\documentclass[review,10pt]{elsarticle}

\usepackage{amsmath} 
\usepackage{amssymb}
\usepackage{amsfonts}
\usepackage{amsthm}
\usepackage{xcolor}
\usepackage{graphicx}
\usepackage{enumerate}
\usepackage{epstopdf} %converting to PDF
\usepackage{array,multirow}
\usepackage{tikz}
\usetikzlibrary{shapes, shadows, arrows,decorations.pathreplacing,calligraphy}
\tikzstyle{block}=[draw,rectangle,fill=blue!5,text width=12 em,text centered, minimum height=12mm, node distance=5 em]
\usepackage{caption}
\usepackage{subcaption}
\tikzstyle{line} = [draw,-latex']
\usetikzlibrary{calc}
\pgfkeys{%
	/polargrid/.cd,
	rmin/.code ={\global\def\rmin {#1}},
	rmax/.code ={\global\def\rmax {#1}},
	amin/.code ={\global\def\amin {#1}},
	amax/.code ={\global\def\amax {#1}},
	rstep/.code={\global\def\rstep{#1}}, 
	astep/.code={\global\def\astep{#1}}}
\def\polargrid{\pgfutil@ifnextchar[{\polar@grid}{\polar@grid[]}}%
\def\polar@grid[#1]{%
	\pgfkeys{/polargrid/.cd,
		rmin ={2},
		rmax ={7},
		amin ={0},
		amax ={360},
		rstep={1}, 
		astep={22.5}}   
	\pgfqkeys{/polargrid}{#1}%
	\pgfmathsetmacro{\addastep}{\amin+\astep} 
	\pgfmathsetmacro{\addrstep}{\rmin+\rstep} 
	\foreach \a in {\amin,\addastep,...,\amax}  \draw[gray] (\a:\rmin) -- (\a:\rmax);  
	\foreach \r in {\rmin,\addrstep,...,\rmax}  \draw[gray] (\amin:\r cm) arc (\amin:\amax:\r cm);    
} 

\usepackage{geometry}
\makeatother   
\usepackage{pgf}

\usepackage{fancyhdr}
\usepackage{bm}
\usepackage{tabularx}
\usepackage{tabularray}

\usepackage{url}
\usepackage{hyperref}

\usepackage{lineno}

\journal{Applied Numerical Mathematics}

\begin{document}
	
	\begin{frontmatter}
		
		
		
		\title{Generalised projective integration scheme in equation-free multiscale modelling}
		
		%%%%%%%%%%%%%%%%%%%%%%%%%%%%%%%%%%%%%%%%%%%
		% Communicated file
		%%%%%%%%%%%%%%%%%%%%%%%%%%%%%%%%%%%%%%%%%%%

		
		\author[add1]{Tanay Kumar Karmakar \corref{cor1}} %% Author name
		\address[add1]{Department of Mathematics, Indian Institute of Technology Guwahati, Assam-781039, India}
		\ead{tanay.kumar@iitg.ac.in}
		\cortext[cor1]{Corresponding author}
		
		\author[add1]{Durga Charan Dalal}
		\ead{durga@iitg.ac.in}
		
		%\affiliation{organization={Department of Mathematics, Indian Institute of Technology Guwahati},%Department and Organization
			%            city={Guwahati},
			%            postcode={781039}, 
			%            state={Assam},
			%            country={India}}
		
		%% Abstract
		\begin{abstract}
			When the spectrum of a system varies significantly over time, fixed choices of macro-, meso-, and micro-time steps, as well as burst lengths, become inadequate, necessitating adaptive and locally informed strategies. To address these challenges, this article proposes a novel and flexible generalised projective integration (\texttt{GPI}) scheme, designed to accommodate time-dependent spectral variation and dynamically evolving scale separation. The proposed framework unifies and extends several existing multiscale methodologies, thereby offering a more general and adaptable computational paradigm. A comprehensive stability analysis of the \texttt{GPI} scheme is carried out, including a detailed investigation of the splitting of the stability region, which forms a central component of this work. Furthermore, problem-dependent strategies for selecting the micro-, meso-, and macro-time steps, as well as the burst length, are developed and their impacts are systematically validated through numerical experiments. To assess the effectiveness of the proposed scheme, three representative problems with distinct types of spectral evolution are considered. The first problem involves a nonlinear stiff system of ordinary differential equations, where the slow eigenvalue remains close to zero along the negative real axis, while the fast eigenvalue increases in magnitude within the negative real plane over time. The second problem examines a linear diffusion equation, in which the entire spectrum evolves dynamically. The third problem considers a highly oscillatory Airy equation, where the eigenvalues always lie on the imaginary axis and move away from the origin as time progresses. The performance of the proposed \texttt{GPI} scheme is evaluated and compared with several existing projective integration methods as well as some widely used stiff solvers, based on (i) number of micro time steps, (ii) accuracy, (iii) computational time, (iv) memory usage and (v) proportion of micro-scale simulations. The results demonstrate that the \texttt{GPI} scheme consistently outperforms the existing methods in terms of performance.
		\end{abstract}
		
		
		%\item A non-rectangular domain or a non-uniform grid can be handled using patch dynamics.
		%\item Patch size and shapes can be non-uniform and non-rectangular.
		%\item A patch dynamics scheme is proposed for generalized curvilinear coordinates.
		%\item The technique is applied to unsteady linear convection-diffusion-reaction equations.
		%\item Patch dynamics scheme can efficiently handle periodic boundary conditions.
		
		%% Keywords
		\begin{keyword}
			multiscale modeling, projective integration, equation-free framework, Numerical methods for \texttt{ODEs} and \texttt{PDEs}, stiff differential equations, highly oscillatory systems
			
		\end{keyword}
		
	\end{frontmatter}
	
%	\linenumbers
	
	%%%%%%%%%%%%%%%%%%%%%%%%%%%%%%%%%%%%%%%%%%%%%%%%%%%%%%%%%%%%%%%%%%%%%%%%%%
	
	\section{Introduction}
Many problems in Science and Engineering are inherently multiscale in nature, involving dynamical processes that evolve over widely separated temporal and spatial scales. Such systems frequently arise from the discretisation of partial differential equations, stochastic differential equations or from detailed microscopic models such as molecular dynamics and kinetic Monte Carlo simulations. A central computational challenge in these settings is the efficient and accurate integration of stiff systems, where rapidly decaying (fast) components coexist with slowly evolving (macroscopic) dynamics of primary interest.

Certain \texttt{ODE} solvers are particularly suited for non-stiff problems, including the explicit Runge--Kutta methods available in \texttt{MATLAB}. For instance, \texttt{ode23} employs the Bogacki–Shampine 3(2) pair \cite{1989_Bogacki_ode23}, while \texttt{ode45} is based on the Dormand–Prince 5(4) pair \cite{1980_Dormand_ode45}. For stiff systems, \texttt{MATLAB}’s \texttt{ode15s} is widely used, relying on backward differentiation formulas (\texttt{BDF}) \cite{1971_Gear_ode15s_BDF} and their modified numerical differentiation formulae (\texttt{NDF}) \cite{1997_Shampine_MATLAB} variants within a variable-step, variable-order framework. The method incorporates local error control and adaptive time stepping to satisfy prescribed tolerances. 
The \texttt{MATLAB} solver \texttt{ode23t} is based on the trapezoidal rule (\texttt{TR}) \cite{1996_Wanner_ode23t}, an implicit one-step method that is unconditionally stable but relatively low in efficiency. In contrast, \texttt{ode23s} employs a modified Rosenbrock-type scheme \cite{1987_Roche_Rosenbrock}, offering improved efficiency for stiff problems within the \texttt{MATLAB ODE} suite. The \texttt{MATLAB} solver \texttt{ode23tb} is based on the \texttt{TR-BDF2} method \cite{1985_Bank_ode23tb,1996_Hosea_ode23tb}, which combines the trapezoidal rule with a second-order backward differentiation formula to achieve improved stability and efficiency. \texttt{Radau IIA} methods \cite{1999_Hairer_Radau_IIA,2026_Mahooti_RadauIIA} are implicit Runge--Kutta schemes based on Radau quadrature nodes, known for their strong stability properties, including A-stability and L-stability. These methods are particularly effective for stiff differential equations, offering high-order accuracy and robust damping of fast transient components.

Classical numerical approaches for stiff systems typically rely on implicit time integration schemes to alleviate severe step size restrictions imposed by stability requirements. However, these methods often entail the repeated solution of large nonlinear systems, making them computationally expensive and difficult to scale, particularly when the underlying models are high-dimensional or available only as legacy simulation codes \cite{2003_Gear_Projective}. This has motivated the development of alternative explicit approaches that can exploit the intrinsic structure of multiscale systems, especially the presence of a spectral gap separating fast and slow modes.

Projective integration (\texttt{PI}) methods \cite{2003_Gear_Projective}, originally introduced for systems with such spectral gaps, provide a framework for accelerating time integration by combining short bursts of fine-scale simulation with extrapolation over larger time intervals. The key idea is to use a stable inner integrator to damp out fast transients, followed by an outer projection step that advances the solution along the slow manifold using polynomial extrapolation. For a class of deterministic multiscale system, Maclean et al. \cite{2014_Maclean_Convergence} presented a convergence analysis of the \texttt{PI} scheme. Givon et al. \cite{2006_Givon_Strong} presented a strong convergence analysis of the \texttt{PI} scheme for singular perturbed stochastic differential systems. The \texttt{PI} is applied on many applications like, liquid crystalline polymers \cite{2003_Siettos_Coarse_Brownian_dynamics}, kinetic Monte Carlo \cite{2004_Rico_Coarse}, stochastic differential systems \cite{2006_Givon_Strong,2021_Przemyslaw_Convergence}, evolving diseases \cite{2004_Jaime_Evolving_Diseases}, movement of organisms and cells \cite{2006_Erban_Equation}, bacterial chemotaxis \cite{2005_Setayeshgar_Bacterial_Chemotaxis}, group-level alignment dynamics of animals moving together \cite{2007_Moon_Heterogeneous_Animal}, heterogeneous cell population dynamics \cite{2007_Bold_Heterogeneous_Cell}, gene regulatory network \cite{2007_Erban_Gene_Regulatory_Network}, disease transmission near eradication \cite{2015_Williams_Disease_Transmission_Near_Eradication}, lattice Boltzmann model \cite{2008_Vandekerckhove_Accuracy}, molecular dynamics \cite{2009_Frederix_Lifting}, materials science \cite{2015_Chuang_coarse}, kinetic theory \cite{2012_Lafitte_Asymptotic}, dynamics of networks \cite{2014_Katherine_equation}, fluid dynamics \cite{2021_Julian_Hyperbolic_Moment}, electric power grid system \cite{2018_Wang_Power_Systems} etc. 

Using multiple projective levels, Gear et al. \cite{2003_Gear_Telescopic} proposed the telescopic projective integration (\texttt{TPI}) method for multiple eigenvalue clusters. Gear et al. \cite{2005_Gear_Projecting} developed a computational framework to initialise dynamical systems on their slow manifold using only a legacy time-stepper, without requiring explicit model equations. The method enforces higher-order derivative conditions to obtain accurate approximations of the missing fast variables, enabling efficient equation-free multiscale computations.
Kavousanakis et al. \cite{2007_Kavousanakis_Projective} enhanced projective and coarse projective integration by exploiting continuous symmetries through a dynamically co-evolving frame, where the effective dynamics become slower and more amenable to extrapolation. By removing translational or scaling effects, the method significantly improves accuracy and enables larger projective time steps for both deterministic and multiscale systems.

Lee et al. \cite{2007_Lee_Second} proposed second-order accurate projective integrators based on Runge--Kutta and Adams--Bashforth formulations, designed as outer schemes for stiff multiscale systems. When combined with telescopic projective integration, these methods yield fully explicit schemes with adaptive time stepping and accuracy comparable to implicit solvers.

Substantial advancements in projective integration have been made, particularly in the context of kinetic theory. For instance, Lafitte et al. \cite{2012_Lafitte_Asymptotic} introduced an asymptotic-preserving projective scheme for kinetic equations in diffusive regimes, ensuring consistency with the corresponding macroscopic limits. Subsequently, Lafitte et al. \cite{2016_Lafitte_High-Order} proposed a high-order projective Runge–Kutta (\texttt{PRK}) method for kinetic equations with linear relaxation, in which a limited number of microsolver steps are employed to estimate time derivatives, followed by high-order extrapolation based on Runge-Kutta scheme. These ideas were further extended by Lafitte et al. \cite{2017_Lafitte_High-order} to nonlinear systems using \texttt{BGK}-type formulations, enabling efficient simulation of multidimensional hyperbolic problems. Maclean et al. \cite{2015_Maclean_Convergence} developed a modified variant of this high-order framework for deterministic multiscale systems with slow--fast structure. Melis et al.~\cite{2018_Melis_Telescopic,2019_Melis_Nonlinear_BGK} introduced telescopic projective integration to address multiscale kinetic equations with multiple relaxation times and subsequently extended this framework to construct high-order explicit projective schemes for nonlinear collisional models such as the \texttt{BGK} and Boltzmann equations. More recent contributions include fully explicit projective methods for multispecies Boltzmann systems~\cite{2022_Rafael_Projective} and for degenerate parabolic equations~\cite{2026_Tenna_Projective}, with stability constraints comparable to classical \texttt{CFL} conditions.

In parallel, efforts have been made to enhance the flexibility of projective integration through adaptivity in space and time. Koellermeier et al. \cite{2022_koellermeier_spatially} introduced spatially adaptive projective integration schemes for stiff hyperbolic balance laws, exploiting spectral gaps arising from spatially varying relaxation times. Their approach applies different time integration strategies across the domain, significantly relaxing stability constraints and improving computational efficiency compared to standard methods. Koellermeier \cite{2025_Koellermeier_Projective} reformulated projective integration methods within a unified Runge--Kutta framework by expressing them through extended Butcher tableaux, enabling systematic analysis of their consistency and order conditions. The study further incorporates spatial and temporal adaptivity via partitioned and embedded Runge--Kutta techniques and rigorously investigates stability, convergence, and error estimation both analytically and numerically.

Recently, George et al. \cite{2026_George_Explicit} introduced explicit time integration schemes based on complex-valued time steps, showing that trajectories in the complex time plane can significantly enlarge stability regions. They demonstrate that such integrators are particularly effective for problems with complex spectra, such as the Schrödinger equation, and can further enhance the efficiency of projective integration methods when applied to stiff systems.

Maclean et al. \cite{2020_Roberts_toolbox,2026_RobertsEquationFree} developed a \texttt{MATLAB}/\texttt{Octave} toolbox implementing equation-free algorithms that enable efficient system-level simulation. The toolbox introduces the coded equation-free functions in an accessible way. Projective integration by second- and fourth-order Runge--Kutta methods is implemented through \texttt{PIRK2} and \texttt{PIRK4}, respectively. These schemes provide accurate approximation of the slow dynamics, provided that the duration of the microsolver bursts remains sufficiently small. The projective integration with a general method (\texttt{PIG}) provides a general formulation of projective integration by allowing the use of any macro-scale time integrator, whether built-in \texttt{MATLAB}/\texttt{Octave} solvers or user-defined. It is particularly effective for highly stiff systems, while for moderately stiff problems it is typically combined with the auxiliary procedure constraint-defined manifold computing (\texttt{cdmc}) \cite{2005_Gear_Projecting}. The \texttt{cdmc} function iteratively applies short microscale simulations together with backward projection steps to drive the fast variables toward the slow manifold without advancing physical time. This correction significantly reduces errors associated with finite microsolver burst lengths and extends the applicability of \texttt{PIG} to a broader class of problems. Based on the macroscale time integrators \texttt{ode23} and \texttt{ode45}, the corresponding \texttt{PIG} schemes are denoted as \texttt{PIG2} and \texttt{PIG4}, respectively, in this article.

%The toolbox includes projective integration schemes such as \texttt{PIRK2} and \texttt{PIRK4} for fixed high-order Runge–Kutta-based macro-integration, as well as the more general \texttt{PIG} framework, which allows the use of arbitrary (including adaptive) macro-integrators through data-driven estimation of slow dynamics.

In multiscale problems, both micro- and macro-level time steps play a fundamental role in accurately capturing system dynamics. The primary objective of this work is to develop a general framework for constructing multiscale, multiphysics methods that offer clear advantages over existing approaches such as systematic upscaling \cite{1977_Brandt_Upscaling,2002_Brandt_Upscaling}, the heterogeneous multiscale method (\texttt{HMM}) \cite{2005_Engquist_Heterogeneous,2012_Abdulle_HMM}, and equation-free \cite{2003_Kevrekidis_Eqn_Free,2026_Karmakar_GPD} techniques. An improved variant of \texttt{HMM}, known as the seamless heterogeneous multiscale method (\texttt{SHMM}), was introduced by Fatkullin et al. \cite{2004_Fatkullin_Computational} and later refined by E \cite{2009_Weinan_General}, where repeated reinitialization of microscale simulations at every macro step is avoided. While most earlier studies focused on two distinct temporal scales—micro and macro—Van Leemput et al. \cite{2008_Van_Mesoscale} first introduced the notion of an intermediate (mesoscopic) time scale within the equation-free framework. Building on this idea, Bunder et al. \cite{2016_Bunder_Accuracy} proposed a modified patch dynamics scheme that exploits mesoscopic scales to reduce communication overhead in large-scale parallel computations. Recently, Karmakar et al. \cite{2026_Karmakar_GPD} proposed the generalised patch dynamics (\texttt{GPD}) scheme, in which several advantages of incorporating a mesoscale time step are discussed. The inclusion of a mesoscale enables the method to capture intermediate dynamics more effectively, leading to improved accuracy with reduced computational cost, and enhances the stability characteristics of the scheme.

To enable proper relaxation of microscale dynamics and their impact on macroscale evolution, a finer mesoscale time step than that used in \texttt{HMM} is employed, giving rise to mesoscale \texttt{HMM} (\texttt{MSHMM}) \cite{2009_Weinan_General}. Vanden-Eijnden \cite{2007_Engquist_Heterogeneous} introduced a framework that avoids explicit identification of slow and fast variables, which was later extended by Tao et al. \cite{2010_Tao_FLAVORS} through the Flow Averaging Integrators (\texttt{FLAVORS}), where the stiff part is alternately turned on and off across micro and meso time scales, respectively. Similarly, boosting algorithms \cite{2009_Weinan_General,2015_Maclean_Note} were developed to effectively reduce stiffness, thereby improving computational efficiency. Lee et al. \cite{2014_Lee_VSHMM} proposed variable step-size \texttt{HMM} (\texttt{VSHMM}), further enhances flexibility by employing adaptive mesoscopic time steps--using finer steps near the boundaries of macro intervals and coarser steps elsewhere to balance accuracy and cost.

%Despite these advancements, increasing the number of interacting time scales often leads to a rapid growth in computational complexity. In contrast, the generalized patch dynamics (GPD) framework provides a unified and flexible structure that naturally accommodates multiple time scales, including micro, meso, and macro levels. By preserving a consistent algorithmic structure, the GPD framework can incorporate and generalize several existing methods—such as HMM, FLAVORS, VSHMM, boosting algorithms, and the UPD scheme—as special cases, depending on the choice of lifting/restriction operators, integrators, and coupling strategies. This highlights the versatility of the GPD approach in addressing complex multiscale phenomena.

\textbf{Aim of this article:}
Within the equation-free framework, the existing schemes were developed primarily for autonomous systems, or for non-autonomous systems in which the eigenvalues exhibit either constant behavior or only mild temporal variation. Consequently, these approaches implicitly assume that the time-scale separation remains nearly uniform throughout the time duration. However, in many real-world applications, the governing systems may not satisfy such restrictive conditions. In particular, the spectral properties of the system may evolve significantly in time, leading to substantial dynamic changes in the degree of time-scale separation.

Such scenarios limit the applicability of existing projective integration schemes to a narrower class of problems. When the spectrum varies strongly with time, fixed choices of macro-, meso- and micro-time steps, as well as burst lengths become inadequate. For such kinds of problems, adaptive and locally informed strategies are required.

To address these challenges, this article proposes a novel and flexible generalised projective integration (\texttt{GPI}) scheme, designed to accommodate time-dependent spectral variation and evolving scale separation. Along with the existing projective integration versions, the proposed framework unifies and extends several other existing multiscale methodologies, including \texttt{HMM}, \texttt{SHMM}, \texttt{MSHMM}, \texttt{FLAVORS}, \texttt{VSHMM} and \texttt{BA} in \texttt{HMM}, thereby providing a broader and more adaptable computational paradigm.

In addition, this work aims to systematically evaluate the performance of the proposed \texttt{GPI} scheme in comparison with classical projective integration methods, such as \texttt{PI}, \texttt{PRK}, \texttt{PIRK} and \texttt{PIG}, as well as widely used stiff solvers such as \texttt{Radau IIA}, \texttt{ode15s}, \texttt{ode23s}, \texttt{ode23t} and \texttt{ode23tb}. A comprehensive stability analysis of the \texttt{GPI} scheme is also carried out, forming a central component of this study. The splitting of the stability region is analysed in detail.

\textbf{Structure of the article:}
This article is organised as follows. In Section~\ref{sec:GPI}, we introduce the generalised projective integration (\texttt{GPI}) scheme formulated across three distinct time scales—macro, meso and micro. In Subsection~\ref{subsec:Features_PI}, we discuss the key features of the proposed \texttt{GPI} scheme and compare them with those of existing projective integration methods. Subsection~\ref{subsec:GPI_Generalised} demonstrates that the \texttt{GPI} scheme provides a unifying framework that generalises several existing multiscale approaches.
In Section~\ref{sec:Stability_GPI}, we present a comprehensive stability analysis of the \texttt{GPI} scheme. Sections~\ref{sec:Relaxation_Time} and \ref{sec:Choice_Micro_Steps} are devoted to the systematic selection of burst length and micro time steps, respectively, based on problem-specific considerations.
Finally, the performance and effectiveness of the proposed \texttt{GPI} scheme are validated through three representative test cases, presented in Subsection~\ref{sec:Results_Discussion}, along with detailed discussions and comparisons with existing methods.

%%%%%%%%%%%%%%%%%%%%%%%%%%%%%%%%%%%%%%%%%%%%%%%%%%%%%%%%

\section{Generalised projective integration (\texttt{GPI}) scheme}\label{sec:GPI}
%	The fundamental idea of projective integration (\texttt{PI}) scheme is to replace the expensive simulation over a long time with sparse computation within small, well separted parts across the time. The \texttt{PI} method is a combination of a few small time-stepper method with a much larger projective time step. The simplest version of the \texttt{PI} scheme is proposed by Gear et al. \cite{2003_Gear_Projective} and is called as ``Projective forward Euler (\texttt{PFE})" method. In this method, to march the solution forrward in time, both the micro simulation in the time-stepper and projective integration are done using the forward Euler scheme, which we are going to discus below. Accuracy of the \texttt{PI} scheme depends on several parameters, such as macro time step ($\Delta T$), relaxation time for the micro-simulation ($K\delta t$) and micro time step ($\delta t$). 

Let the deterministic multiscale systems be represented by:
\begin{equation}\label{eqn:GPI_General_Eqn}
	\frac{du}{dt} = \mathcal{F}\!\left(t,u\right),\hspace{0.5cm} u(t_0)=u_0,
\end{equation}
where $t \in [t_0,T]$ denotes the time variable, $u \in \mathbb{R}^\mathsf{n}$ is the unknown solution vector and the right-hand side function $\mathcal{F} \in \mathbb{R}^\mathsf{n}$ is assumed to be smooth. In equation \eqref{eqn:GPI_General_Eqn}, the right-hand side function might originate from the discretisation of the spatial partial derivatives, and the left-hand side would be a partial time derivative. In most of the variants of projective integrations under the equation-free framework, usually, the system \eqref{eqn:GPI_General_Eqn} is considered as autonomous or non-autonomous with a time-independent Jacobian. However, in this article, the system \eqref{eqn:GPI_General_Eqn} is considered as non-autonomous, that is, the function $\mathcal{F}$ depends on the time variable $t$ as well as the Jacobian may depend on time. Since the system is nonlinear and time-dependent, stiffness is analysed via local linearisation. Let
\begin{equation}
	J(t,u)=\frac{\partial \mathcal{F}}{\partial u}(t,u),
\end{equation}
be the Jacobian. Suppose $\lambda_p(t)$ are the eigenvalues of the Jacobian $J(t,u(t))$ for some finite number of positive integer values of $p$, where the eigenvalues may depend on the local time. The system \eqref{eqn:GPI_General_Eqn} is stiff if the eigenvalues $\lambda_p(t)$ satisfy (i) $Re(\lambda_p(t))<0$ (dissipative behaviour) and (ii) there exists a large separation:
\begin{center}
	$\max_{p}|\lambda_p(t)|\gg\min_{p}|\lambda_p(t)|$.
\end{center}
For the general system \eqref{eqn:GPI_General_Eqn}, the scale separation parameter can be defined locally in time as
\begin{center}
	$\epsilon_\text{stiff}(t)=\frac{\max_{p}|\lambda_p(t)|}{\min_{p}|\lambda_p(t)|}$.
\end{center}
Throughout the following discussion, only real eigenvalues are considered.

%	Let the deterministic multiscale systems be represented by:
%	\begin{equation}\label{eqn:Work3_General_Eqn}
	%		\frac{du}{dt} = \mathcal{F}\!\left(t,u\right),
	%	\end{equation}
%	where $t \in [t_0,T]$ denotes the time variable, $u \in \mathbb{R}^\texttt{n}$ is the unknown solution vector and the right-hand side function $\mathcal{F} \in \mathbb{R}^\texttt{n}$ is assumed to be smooth. In general, the scale-separation parameter may depend on the local time and, in certain cases, also on the local solution state. The proposed \texttt{GPI} scheme is therefore formulated in a general setting based on equation \eqref{eqn:Work3_General_Eqn}.

We discretise the entire time interval $[t_0, T ]$ into macroscopic time levels $\{\, T^n : 0 \le n \le \operatorname{N_t},\ T^0 = t_0 < T^1 < \cdots < T^n < \cdots < T^{\operatorname{N_t}} = T \,\}
$, where $\Delta T^n:=T^{n+1}-T^n$ denotes the macroscopic time step for $n=0,\ldots,N_t-1$.	%A schematic illustration of the discretisation of the interval $[t_0,T]$ into variable macroscopic time steps.  is shown in Figure \ref{fig:GPI_Macro}.

%\begin{figure}%[h!]
%	\centering
%	\begin{tikzpicture}[scale=0.7]
	%		%---------------------------
	%		%		\draw[thick] (-2.5,0) node {\Huge{$(a)$}};
	%		% PI_GPD-I
	%		% Macro_GPI_I		
	%		\draw[thick,->] (0,0) -- (12+1,0);
	%		\draw[thick] (13.5,0) node {$t$};
	%		\draw[very thick, magenta] (0,0-0.6) -- (0,0+0.6);
	%		\draw[thick] (0,-1) node {$T^0=t_0$};
	%		\draw[very thick, magenta] (12,0-0.6) -- (12,0+0.6);
	%		\draw[thick] (12.5,-1) node {$T^\text{Nt}$=$T$};
	%		
	%		\draw[very thick, magenta] (1,0-0.6) -- (1,0+0.6);
	%		\draw[very thick, magenta] (3.5,0-0.6) -- (3.5,0+0.6);
	%		\draw[very thick, magenta] (5,0-0.6) -- (5,0+0.6);
	%		\draw[very thick, magenta] (7,0-0.6) -- (7,0+0.6);
	%		\draw[very thick, magenta] (8,0-0.6) -- (8,0+0.6);
	%		\draw[very thick, magenta] (10.5,0-0.6) -- (10.5,0+0.6);
	%		
	%		\draw[thick] (3.5,-1) node {$T^{n-1}$};
	%		\draw[thick] (5,-1) node {$T^{n}$};
	%		\draw[thick] (7,-1) node {$T^{n+1}$};
	%		
	%		\draw[thick] (5,-1.7-0.3) -- (5,-1.7+0.3);
	%		\draw[thick] (5,-1.7) -- (5.4,-1.7);
	%		\draw[thick] (6,-1.7) node {$\Delta T^{n}$};
	%		\draw[thick] (6.6,-1.7) -- (7,-1.7);
	%		\draw[thick] (7,-1.7-0.3) -- (7,-1.7+0.3);
	%		
	%		\draw[thick,->, magenta] (7.2,-1.7) -- (8,-1.7);
	%		\draw[thick, magenta] (10,-1.7) node {\small{Macro time step}};
	%		
	%	\end{tikzpicture}%\pause
%	\caption{Schematic illustration of the discretisation of the entire time interval $[t_0,T]$ into variable macro time steps $\Delta T^n$ in the \texttt{GPI} scheme.} \label{fig:GPI_Macro}
%\end{figure}

In each macroscopic time step $\Delta T^n$, we introduce a variable number of mesoscopic time steps with non-uniform sizes. Specifically, the macroscopic time interval $[T^n, T^{n+1}]$ is discretised into $l_n$ mesoscopic subintervals in a non-uniform manner. The corresponding mesoscopic time levels are denoted by $\{\, T^{n,m} : 0 \le m \le l_n,\ T^{n,0} = T^n < T^{n,1} < \cdots < T^{n,m} < \cdots < T^{n,l_n} = T^{n+1,0} = T^{n+1} \,\}$. The notation $\partial T^{n,m}:=T^{n,m+1}-T^{n,m}$, $m=0, \ldots, l_n-1$ denotes the $(m+1)^{th}$ mesoscopic time step in $(n+1)^{th}$ macro time step. These mesoscopic time steps satisfy 
\begin{equation}
	\sum_{m=0}^{l_n-1}\partial T^{n,m}=\Delta T^n.
\end{equation}

Due to the presence of time-varying spectrum, the duration of the microscale simulations may vary. Consequently, both the number of micro time steps and their sizes may vary within a meso time step. To maintain the stability of the scheme, the projective extrapolation time step may also need to vary accordingly. Thus, within each mesoscopic time interval $\partial T^{n,m}$, we perform $M_{n,m}+1$ number of micro time steps with variable sizes $\delta T^{n,m,p}:=T^{n,m,p+1}-T^{n,m,p}$, for $p=0,1,\ldots,M_{n,m}$. The corresponding micro time levels are denoted by $T^{n,m,p}$, for $p=0, \ldots, M_{n,m}+1$. Effectively, the projective extrapolation step size is given by
\begin{center}
	$\partial T^{n,m}-\sum_{i=0}^{M_{n,m}}\delta T^{n,m,p}$.
\end{center} 

Starting from the solution $u^{n,m}$ at time $T^{n,m}$, the inner integrator is first applied for $M_{n,m}+1$ small steps:
\begin{equation}\label{eqn:Micro_Discretisation_General}
	u^{n,m,p+1}=\Phi_{\delta T^{n,m,p}}(u^{n,m,p}), \hspace{0.2cm}\forall p=0, \ldots, M_{n,m},
\end{equation}
where $\Phi_{\delta T^{n,m,p}}$ represents the microscopic evolution operator over the inner time step $\delta T^{n,m,p}$.

%The variable micro time steps $\delta T^{n,m,p}$ for $p=0,\ldots,M_{n,m}+1$ are indicated in red solid line in the schematic Figure \ref{fig:GPI_Micro}.
The final micro time step $\delta T^{n,m,M_{n,m}}$ is used to estimate an approximate value of the time derivative of $u$ at $T^{n,m,M_{n,m}+1}$ as 
\begin{equation}\label{eqn:GPI_Slope}
	k_{n,m,M_{n,m}+1}:=\frac{u^{n,m,M_{n,m}+1}-u^{n,m,M_{n,m}}}{\delta T^{n,m,M_{n,m}}}.
\end{equation}
We apply the forward Euler method to march $u$ forward from time $T^{n,m,M_{n,m}+1}$ over a long time step $\partial T^{n,m}-\sum_{i=0}^{M_{n,m}}\delta T^{n,m,p}$ to reach the macroscopic time $T^{n,m+1}$ such that
\begin{equation}\label{eqn:Meso_Projective}
	u^{n,m+1}=u^{n,m,M_{n,m}+1}+\left(\partial T^{n,m}-\sum_{p=0}^{M_{n,m}}\delta T^{n,m,p}\right)k_{n,m,M_{n,m}+1}.
\end{equation}
%A schematic representation of the projective extrapolation step is shown by black solid line in Figure \ref{fig:GPI_Micro}.

We could also write the equation \eqref{eqn:Meso_Projective} as 
\begin{equation}\label{eqn:Meso_Projective_Discretisation}
	u^{n,m+1}=\left(R^{n,m}+1\right)u^{n,m,M_{n,m}+1}-R^{n,m} u^{n,m,M_{n,m}},
\end{equation}
where,
\begin{equation}
	R^{n,m}:=\frac{\partial T^{n,m}-\sum_{p=0}^{M_{n,m}}\delta T^{n,m,p}}{\delta T^{n,m,M_{n,m}}}\in[0,\infty).
\end{equation}

For all microscale simulations presented in this article, we employ the forward Euler method in equation \eqref{eqn:Micro_Discretisation_General} as the inner time integrator:
\begin{equation}\label{eqn:Micro_Discretisation}
	u^{n,m,p+1}=u^{n,m,p}+\delta T^{n,m,p}\mathcal{F}(T^{n,m,p},u^{n,m,p}), \hspace{0.2cm}\forall p=0,\ldots,M_{n,m}.
\end{equation}

While higher-order inner integrators can potentially yield modest improvements in accuracy, such improvements are generally outweighed by the associated increase in computational cost. This observation is consistent with findings reported in the equation-free literature by Karmakar et al. \cite{2024_Karmakar_generalized} and Maclean et al. \cite{2015_Maclean_Convergence}. % by Karmakar et al. \cite{2024_Karmakar_generalized}

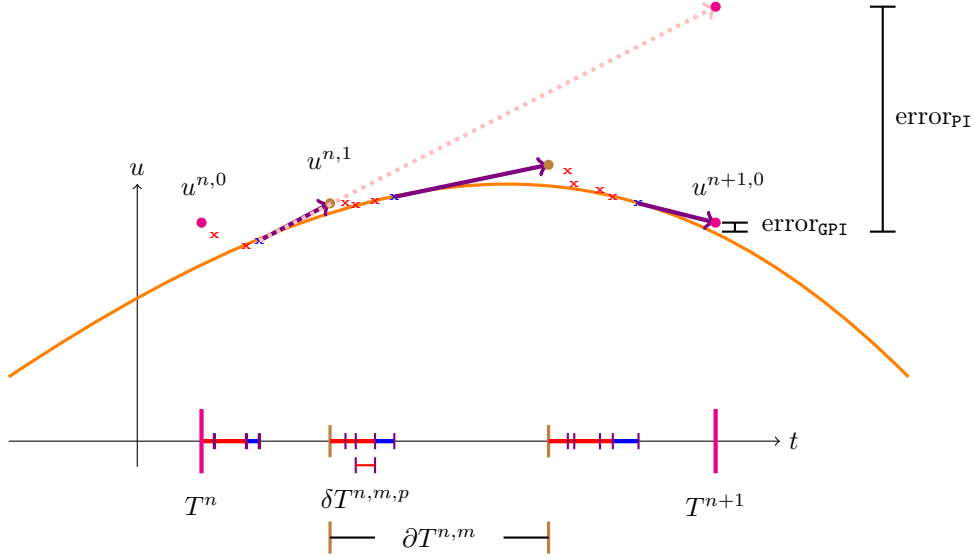
\begin{figure}
	\centering
	\begin{tikzpicture}[scale=1.7]
		% Axis
		\draw[->] (-1,0) -- (5,0) node[right] {$t$};
		\draw[->] (0,-0.2) -- (0,2) node[above] {$u$};
		
		% Curve
		\draw[orange, very thick] (-1,0.5) .. controls (2,2.5) and (4,2.5) .. (6,0.5);
		
		%% First meso step
		% Meso 1 (t-axis)
		\draw[ultra thick,red] (0.5,0) -- (0.85,0);
		\draw[ultra thick,blue] (0.85,0) -- (0.95,0);
		\draw[ultra thick,magenta] (0.5,-0.25) -- (0.5,0.25);
		\draw[ultra thick,magenta] (4.5,-0.25) -- (4.5,0.25);
		\draw[very thick] (0.5,-0.5)  node {$T^{n}$};
		\draw[very thick] (4.5,-0.5)  node {$T^{n+1}$};
		\draw[very thick,violet] (0.60,-0.0625) -- (0.60,0.0625);
		\draw[very thick,violet] (0.85,-0.0625) -- (0.85,0.0625);
		\draw[very thick,violet] (0.95,-0.0625) -- (0.95,0.0625);
		
		% Meso 1 (u-axis)
		\filldraw[magenta] (0.5,1.7) circle (1pt);
		\draw[very thick] (0.52,2)  node {$u^{n,0}$};
		\draw[ultra thick, red] (0.60,1.6)  node {\bf{\tiny x}};
		\draw[ultra thick, red] (0.85,1.52)  node {\bf{\tiny x}};
		\draw[ultra thick, blue] (0.95,1.56)  node {\bf{\tiny x}};
		
		\draw[ultra thick,violet,->] (0.95,1.56) -- (1.5,1.85);
		\filldraw[brown] (1.5,1.85) circle (1pt);
		
		%% Second meso step
		% Meso 2 (t-axis)
		\draw[very thick,brown] (1.5,-0.125) -- (1.5,0.125);
		\draw[ultra thick,red] (1.5,0) -- (1.85,0);
		\draw[ultra thick,blue] (1.85,0) -- (2,0);
		\draw[thick,violet] (1.62,-0.0625) -- (1.62,0.0625);
		\draw[thick,violet] (1.70,-0.0625) -- (1.70,0.0625);
		\draw[thick,violet] (1.85,-0.0625) -- (1.85,0.0625);
		\draw[thick,violet] (2,-0.0625) -- (2,0.0625);
		
		% Meso 2 (u-axis)
		\draw[very thick] (1.5,2.2)  node {$u^{n,1}$};
		\draw[ultra thick, red] (1.62,1.85)  node {\bf{\tiny x}};
		\draw[ultra thick, red] (1.70,1.84)  node {\bf{\tiny x}};
		\draw[ultra thick, red] (1.85,1.87)  node {\bf{\tiny x}};
		\draw[ultra thick, blue] (2,1.9)  node {\bf{\tiny x}};
		
		\draw[ultra thick,violet,->] (2,1.9) -- (3.2,2.15);
		\filldraw[brown] (3.2,2.15) circle (1pt);
		
		%% Last meso step
		% Meso 2 (t-axis)
		\draw[very thick,brown] (3.2,-0.125) -- (3.2,0.125);
		\draw[ultra thick,red] (3.2,0) -- (3.7,0);
		\draw[ultra thick,blue] (3.7,0) -- (3.9,0);
		\draw[thick,violet] (3.35,-0.0625) -- (3.35,0.0625);
		\draw[thick,violet] (3.40,-0.0625) -- (3.40,0.0625);
		\draw[thick,violet] (3.60,-0.0625) -- (3.60,0.0625);
		\draw[thick,violet] (3.70,-0.0625) -- (3.70,0.0625);
		\draw[thick,violet] (3.90,-0.0625) -- (3.90,0.0625);
		
		% Meso 2 (u-axis)
		\draw[ultra thick, red] (3.35,2.1)  node {\bf{\tiny x}};
		\draw[ultra thick, red] (3.4,2)  node {\bf{\tiny x}};
		\draw[ultra thick, red] (3.6,1.95)  node {\bf{\tiny x}};
		\draw[ultra thick, red] (3.7,1.90)  node {\bf{\tiny x}};
		\draw[ultra thick, blue] (3.9,1.85)  node {\bf{\tiny x}};
		
		\draw[ultra thick,violet,->] (3.9,1.85) -- (4.5,1.70);
		\filldraw[magenta] (4.5,1.70) circle (1pt);
		\draw[very thick] (4.6,2)  node {$u^{n+1,0}$};
		
		%% PI
		\draw[ultra thick,dotted,pink,->] (0.95,1.56) -- (4.5,3.38);
		\filldraw[magenta] (4.5,3.38) circle (1pt);
		
		%% Error
		% GPI
		\draw[thick] (4.55,1.70) -- (4.75,1.70);
		\draw[thick] (4.55,1.63) -- (4.75,1.63);
		\draw[thick] (4.65,1.63) -- (4.65,1.70);
		\draw[very thick] (5.2,1.665)  node {$\text{error}_\texttt{GPI}$};
		
		% PI
		\draw[thick] (5.7,3.38) -- (5.9,3.38);
		\draw[thick] (5.7,1.63) -- (5.9,1.63);
		\draw[thick] (5.8,1.63) -- (5.8,3.38);
		\draw[very thick] (6.2,2.5)  node {$\text{error}_\texttt{PI}$};
		
		\draw[thick,violet] (1.70,-2*0.0625) -- (1.70,-4*0.0625);
		\draw[thick,violet] (1.85,-2*0.0625) -- (1.85,-4*0.0625);
		\draw[thick,red] (1.70,-3*0.0625) -- (1.85,-3*0.0625);
		\draw[very thick] (1.775,-0.45)  node {$\delta T^{n,m,p}$};
		
		\draw[very thick,brown] (1.5,-5*0.125) -- (1.5,-7*0.125);
		\draw[very thick,brown] (3.2,-5*0.125) -- (3.2,-7*0.125);
		\draw[thick] (1.5,-6*0.125) -- (1.85,-6*0.125);
		\draw[very thick] (2.35,-6*0.125)  node {$\partial T^{n,m}$};
		\draw[thick] (2.85,-6*0.125) -- (3.2,-6*0.125);
		
	\end{tikzpicture}
	\caption{Schematic representation of the \texttt{GPI} scheme.}
	\label{fig:GPI_Schematic}
\end{figure}

	Figure \ref{fig:GPI_Schematic} presents a schematic representation of the \texttt{GPI} scheme. The orange-coloured curve represents the slow manifold associated with the problem \eqref{eqn:GPI_General_Eqn}. The scheme is applied over a macro time step $\Delta T^n$ to advance the solution from the macro time level $T^n$ to the subsequent macro level $T^{n+1}$ according to the algorithm described above. The macro time levels are represented by magenta vertical lines. Within a single macro time step $\Delta T^n$, three non-uniform meso time steps are employed, whose corresponding meso time levels are indicated by brown vertical lines. Each meso step consists of non-uniform micro time steps (shown by red and blue segments along the $t$-axis), together with non-uniform microsimulations and extrapolation steps. The micro time levels are separated by violet vertical lines.
	
	The solutions at the macro, meso and micro levels are represented by magenta disks, brown disks and red-blue crosses, respectively. The slopes used in the \texttt{GPI} scheme are indicated by solid violet arrows. For the same macro time step $\Delta T^n$ and a comparable micro-burst duration in the first meso step, if the projective integration (\texttt{PI}) scheme \cite{2003_Gear_Projective} is employed, the corresponding slope is represented by a pink dashed arrow. It can be observed from the solution at $T^{n+1}$ that the \texttt{PI} scheme produces a larger error ($\text{error}_{\texttt{PI}}$) compared with the \texttt{GPI} scheme ($\text{error}_{\texttt{GPI}}$). This improved accuracy of the \texttt{GPI} scheme arises because the intermediate meso steps help maintain the solution closer to the slow manifold throughout the integration process.

	%%%%%%%%%%%%%%%%%%%%%%%%%%%%%%%%%%%%%%%%%%%%%%%%%%%%%%%%
	
	In order to express the \texttt{GPI} scheme in the standard Runge--Kutta form, we introduce the following notation:
	\begin{equation}
		\begin{aligned}
			k_{n,m,i}&=\mathcal{F}\left(T^{n,m}+\sum_{j=0}^{i-1}\delta T^{n,m,j}, u^{n,m}+\sum_{j=0}^{i-1}\delta T^{n,m,j}k_{n,m,j}\right), \\
			&=\mathcal{F}\left(T^{n,m}+\sum_{j=0}^{i-1}\delta T^{n,m,j}, u^{n,m}+\partial T^{n,m}\sum_{j=0}^{i-1}\alpha_{n,m,j}k_{n,m,j}\right), \hspace{0.2cm}0\le i\le M_{n,m}
		\end{aligned}
	\end{equation}
	where $\alpha_{n,m,j}:=\frac{\delta T^{n,m,j}}{\partial T^{n,m}}$, for $j=0,\ldots,M_{n,m}$. The solution at the new time step is given by equation \eqref{eqn:Meso_Projective}
	\begin{equation}
		\begin{aligned}
			u^{n,m+1}&=u^{n,m}+\partial T^{n,m}\left(\sum_{j=0}^{M_{n,m}}\frac{\delta T^{n,m,j}}{\partial T^{n,m}}k_{n,m,j}+ \left(1-\sum_{j=0}^{M_{n,m}}\frac{\delta T^{n,m,j}}{\partial T^{n,m}}\right)k_{n,m,M_{n,m}}\right),\\
			&=u^{n,m}+\partial T^{n,m}\left(\sum_{j=0}^{M_{n,m}}\alpha_{n,m,j}k_{n,m,j}+ \left(1-\sum_{j=0}^{M_{n,m}}\alpha_{n,m,j}\right)k_{n,m,M_{n,m}}\right).
		\end{aligned}
	\end{equation}
	
	The \texttt{GPI} method is expressed as a Runge--Kutta method using a block Butcher tableau as follows:

	\begin{equation}\label{eqn:Butcher_Tableau_GPI}
		\renewcommand\arraystretch{1.2}
		\begin{array}{c|cc}
			\boldsymbol{\hat{c}} &  \boldsymbol{\hat{A}}  \\
			\hline
			& \boldsymbol{\hat{b}^T}
		\end{array}
		=
		\begin{array}
			{c|cccccc}
			\hat{c}_1&0\\
			\hat{c}_2&\alpha_{n,m,0}&0\\
			\hat{c}_3&\alpha_{n,m,0}&\alpha_{n,m,1}&0\\
			\vdots&\vdots&\vdots&\vdots\\
			\hat{c}_{M_{n,m}+1}&\alpha_{n,m,0}&\alpha_{n,m,1}&\alpha_{n,m,2}&\hdots&\alpha_{n,m,M_{n,m}-1}&0\\
			\hline
			&\hat{b}_1&\hat{b}_2&\hat{b}_3&\hdots&\hat{b}_{M_{n,m}}&\hat{b}_{M_{n,m}+1} 
		\end{array}
	\end{equation}
	where 
	\begin{center}
		$\boldsymbol{\hat{c}}=\left(\hat{c}_i\right){}^{\top}\in\mathbb{R}^{M_{n,m}+1}$, $\hat{c}_i=\sum_{j=0}^{i-2}\alpha_{n,m,j}$ for $i=1,\ldots,M_{n,m}+1$
	\end{center}
	and 
	\begin{center}
		$\boldsymbol{\hat{b}^T}=\left(\hat{b}_i\right){}^{\top}\in\mathbb{R}^{M_{n,m}+1}$, $\hat{b}_i=\alpha_{n,m,i-1}$, $i=1,\ldots,M_{n,m}$ and $\hat{b}_{M_{n,m}+1}=1-\sum_{j=0}^{M_{n,m}-1}\alpha_{n,m,j}$.
	\end{center}
	
	It is interesting to observe that the Butcher tableau does not depend on the final micro time step of the corresponding microsimulation within a meso time step.

	%\left(\alpha_{n,m,0},\alpha_{n,m,1} \ldots, \alpha_{n,m,M_{n,m}-1}, \left(1-\sum_{j=1}^{M_{n,m}-1}\alpha_{n,m,j}\right)\right)^T
	
	\subsection{Discussion on the features of various projective integration schemes}\label{subsec:Features_PI}
	\begin{table}
		\centering
		\caption{A comparison of the features of various projective integration schemes within the equation-free framework.}
		\label{table:Comparison_Features_PI}
		\begin{tabular}{ccccc}
			\hline
			Scheme & Macro step & Meso step & Micro step & Burst length\\
			\hline
			\texttt{GPI} (present) & variable & variable & variable & variable \\
			\texttt{PI} \cite{2003_Gear_Projective,2012_Lafitte_Asymptotic} & fixed & N/A & fixed & fixed \\
			\texttt{PRK} \cite{2016_Lafitte_High-Order,2017_Lafitte_High-order} & fixed & N/A & fixed & fixed \\
			\texttt{PIRK2} \cite{2021_maclean_toolbox} & fixed & N/A & variable & fixed \\
			\texttt{PIRK4} \cite{2021_maclean_toolbox} & fixed & N/A & variable & fixed \\
			\texttt{PIG2} \cite{2021_maclean_toolbox,2005_Gear_Projecting} & variable & N/A & variable & fixed \\
			\texttt{PIG4} \cite{2021_maclean_toolbox,2005_Gear_Projecting} & variable & N/A & variable & fixed \\
			\hline
		\end{tabular}
	\end{table}
	
	The comparison presented in Table~\ref{table:Comparison_Features_PI} highlights the key differences among various projective integration schemes within the equation-free framework. A primary distinction lies in the flexibility of time stepping across different scales. The \texttt{GPI} scheme represents the most general formulation, allowing variable macro, meso and micro time steps, as well as a variable burst length. This flexibility makes it highly adaptable to complex multiscale systems. 
	
	In comparison with all other schemes listed in Table \ref{table:Comparison_Features_PI}, the \texttt{GPI} scheme is the only one that incorporates meso time steps of variable sizes. Moreover, it allows a variable burst length, whereas the existing schemes do not provide such flexibility.
	
	In contrast, classical projective integration (\texttt{PI}) \cite{2003_Gear_Projective,2012_Lafitte_Asymptotic} and projective Runge--Kutta (\texttt{PRK}) \cite{2016_Lafitte_High-Order,2017_Lafitte_High-order} schemes employ fixed macro and micro time steps along with a fixed burst length and they do not include an explicit meso scale. The projective integration by second and fourth-order Runge--Kutta (\texttt{PIRK} family such as \texttt{PIRK2} and \texttt{PIRK4}) \cite{2021_maclean_toolbox} introduces variable micro time stepping while keeping the macro time step and burst length fixed, thereby partially improving adaptability. These schemes achieve second and fourth-order accuracy, respectively, at the macroscale, but they do not incorporate a meso time scale. 
	
	Similarly, the projective integration via a general macroscale integrator (\texttt{PIG}) schemes \cite{2021_maclean_toolbox,2005_Gear_Projecting} allow variability in both macro and micro time steps, offering greater flexibility than the \texttt{PIRK} methods. However, they do not include an explicit mesoscale or a variable burst length. This limitation may reduce their effectiveness in problems where an intermediate scale plays a significant role.
	
	Overall, the progression from \texttt{PI} to \texttt{GPI} reflects a trade-off between simplicity and flexibility. Compared to the existing \texttt{PI}, \texttt{PRK}, \texttt{PIRK} and \texttt{PIG} families, the \texttt{GPI} scheme offers significantly greater flexibility in selecting macro, meso and micro time steps, as well as the burst length. Due to the lack of such flexibility, existing schemes may struggle to handle problems with time-dependent spectra, whereas the \texttt{GPI} scheme is better suited for such cases, as discussed in Section \ref{sec:Results_Discussion}.
	
	Depending on the requirements of the underlying problem, the full flexibility of the proposed framework may be utilised. However, when the system is autonomous or when the scale-separation parameter exhibits no significant temporal variation, simpler versions of the \texttt{GPI} scheme may be employed, such as uniform choices of the macro, meso and micro time steps, together with uniform burst lengths.
	
	\subsection{The \texttt{GPI} scheme is a generalised version of many other time integration schemes}\label{subsec:GPI_Generalised}
	\begin{table}
		\centering
		\caption{The \texttt{GPI} scheme acts as a general framework for a wide range of existing projective integration and other multiscale methods. }
		\label{table:GPI_Generalised}
		\begin{tabular}{ccc}
			\hline
			Proposed scheme & Values of $l_n$ and $M_{n,m}$& The proposed scheme \\
			&&is equivalent to \\
			\hline
			$\mathsf{(U)(\times)(U)(U)}-\texttt{GPI}$ & $l_n=1$, $M_{n,0}$ is constant & \texttt{PI} \cite{2003_Gear_Projective}\\
			$\mathsf{(U)(\times)(U)(U)}-\texttt{GPI}$ & $l_n=1$, $M_{n,0}=0$& \texttt{PI} version of \texttt{PD} \cite{2006_Samaey_PD_buffers,2020_arbabi_linking}\\
			$\mathsf{(U)(N)(U)(N)}-\texttt{GPI}$ & $l_n\ge1$ is constant, $M_{n,m}$ is variable & \texttt{PI} version of \texttt{GPD-I} \cite{2026_Karmakar_GPD}\\
			$\mathsf{(U)(\times)(U)(U)}-\texttt{GPI}$ & $l_n=1$, $M_{n,0}$ is constant & \texttt{PRK} \cite{2016_Lafitte_High-Order,2017_Lafitte_High-order} of order one\\
			\hline
			$\mathsf{(U)(\times)(U)(U)}-\texttt{GPI}$ & $l_n=1$, $M_{n,0}$ is constant & \texttt{HMM} \cite{2005_Engquist_Heterogeneous,2012_Abdulle_HMM}\\
			$\mathsf{(U)(U)(U)(U)}-\texttt{GPI}$ & $l_n>1$ is constant, $M_{n,m}=1$ & \texttt{FLAVORS} \cite{2010_Tao_FLAVORS}\\
			$\mathsf{(U)(N)(U)(U)}-\texttt{GPI}$ & $l_n>1$ is constant, $M_{n,m}=1$& \texttt{VSHMM} \cite{2014_Lee_VSHMM}\\
			$\mathsf{(U)(U)(U)(U)}-\texttt{GPI}$ & $l_n>1$ is constant, $M_{n,m}=1$& \texttt{BA} in \texttt{HMM} \cite{2015_Maclean_Convergence}\\
			\hline
		\end{tabular}
	\end{table}
	
	Table~\ref{table:GPI_Generalised} illustrates how the \texttt{GPI} scheme acts as a general framework for a wide range of existing projective integration and other multiscale methods. In the notation $\mathsf{(\cdot)(\cdot)(\cdot)(\cdot)}-\texttt{GPI}$, the four brackets respectively denote the characteristics of the macro time step, meso time step, micro time step and burst length, as discussed in Table \ref{table:Comparison_Features_PI}. Here, $\mathsf{U}$ represents uniform (or, fixed), $\mathsf{N}$ represents non-uniform (or, variable) and $\times$ indicates that the corresponding scale is not present (or, not available).
	
	The table shows that several classical schemes can be recovered as special cases of the \texttt{GPI} framework through appropriate choices of $l_n$ and $M_{n,m}$. For instance, when $l_n=1$ and $M_{n,0}$ is constant, the \texttt{GPI} scheme reduces to the classical \texttt{PI} method \cite{2003_Gear_Projective}. If $M_{n,0}=0$, it corresponds to the \texttt{PI} version of patch dynamics (\texttt{PD}) scheme \cite{2006_Samaey_PD_buffers,2020_arbabi_linking}. Similarly, by allowing $l_n \geq 1$ to be constant and $M_{n,m}$ to vary, one obtains the \texttt{PI} version of the generalised patch dynamics (\texttt{GPD}) scheme \cite{2026_Karmakar_GPD} of type-I. The first-order projective Runge--Kutta (\texttt{PRK}) method \cite{2016_Lafitte_High-Order,2017_Lafitte_High-order} is also recovered under the same configuration as the classical \texttt{PI} scheme.
	
	Furthermore, the table indicates that several multiscale methods--such as the heterogeneous multiscale method (\texttt{HMM}) \cite{2005_Engquist_Heterogeneous,2012_Abdulle_HMM}, flow averaging integrators (\texttt{FLAVORS}) \cite{2010_Tao_FLAVORS}, variable step size heterogeneous multiscale methods (\texttt{VSHMM}) \cite{2014_Lee_VSHMM} and the boosting algorithm (\texttt{BA}) within \texttt{HMM} \cite{2015_Maclean_Convergence}--can be interpreted within the \texttt{GPI} framework. These methods retain the same overall structure as the \texttt{GPI} scheme, except for the final micro time step (highlighted in blue in Figure~\ref{fig:GPI_Schematic}).
	
	In this unified perspective, variants of \texttt{HMM} can be formulated through their own appropriate choices of macro to micro operators, micro to macro operators, as well as inner and outer integrator techniques. The primary distinctions among these methods arise from their treatment of the mesoscopic scale and the associated burst length. For instance, \texttt{FLAVORS} and \texttt{VSHMM} correspond to the choice $M_{n,m}=1$ with multiple meso steps $l_n>1$, using uniform and non-uniform step sizes, respectively. In contrast, \texttt{HMM} corresponds to the case where the meso scale is absent, i.e., $l_n=1$, with a fixed value of $M_{n,0}$.
	
	Moreover, by setting $M_{n,m}=1$ and considering a fixed finite number of uniform meso steps in the \texttt{GPI} framework (except for the final micro time step within each meso time step, coloured by blue in Figure \ref{fig:GPI_Schematic}) and by increasing the stiffness parameter, one recovers the \texttt{BA} formulation within \texttt{HMM}. This boosting strategy provides significant computational advantages.
	
	Overall, the table highlights the versatility of the \texttt{GPI} scheme. By systematically varying the structure of time stepping across macro, meso and micro levels, along with the burst length, the \texttt{GPI} framework is capable of reproducing many existing methods. This unified perspective not only clarifies the relationships among different schemes but also provides a flexible foundation for designing new multiscale algorithms.
	
	%%%%%%%%%%%%%%%%%%%%%%%%%%%%%%%%%%%%%%%%%%%%%%%%%%%%%%%%%%%%%%%%%%%%%%%%%%%%%%%%%%%%%%%
	
	\section{Stability analysis of the \texttt{GPI} scheme}\label{sec:Stability_GPI}
	%We assume that the inner integrator is a linear numerical method, meaning that it commutes with any linear transformation of the dependent variables. In other words, applying the method to the original system and then transforming the solution by a nonsingular constant matrix produces the same result as first transforming the differential equation and then applying the numerical method. Most commonly used integration schemes satisfy this property.
	%
	%When such a linear method is applied to a linear system with constant coefficients $u'=Au$, the computation can be decomposed into independent scalar problems corresponding to the eigenvalues of $A$. Therefore, the stability properties of the method can be studied by analysing its behaviour on the Dahlquist test equation 
	We now study the linear stability analysis of the \texttt{GPI} scheme introduced above. To this end, we consider the Dahlquist test equation,
	\begin{equation}\label{eqn:Test_eqn}
		\frac{du}{dt}=\lambda(t) u, \hspace{0.2cm} \lambda(t)<0.
	\end{equation}
	The notation $\lambda(t)$ denotes the eigenvalues of the system at time $t$. 
	
	The one step explicit inner integrator for the problem \eqref{eqn:Test_eqn} over micro step size $\delta T^{n,m,p}$ starting from the time $T^{n,m,p}$ will be
	\begin{equation}\label{eqn:Test_eqn_Inner_Integrator}
		u^{n,m,p+1}=\rho^{n,m,p}u^{n,m,p},
	\end{equation}
	where, $p=0, 1, \ldots, M_{n,m}$; $m=0, 1, \ldots, l_n-1$ and $n=0, 1, \ldots, \operatorname{N_t}-1$. The notation $\rho^{n,m,p}:=\rho(\lambda^{n,m,p}\delta T^{n,m,p})$ denotes the amplification of the method over the micro time step $\delta T^{n,m,p}$ and $\lambda^{n,m,p}:=\lambda(T^{n,m,p})$. For perfect inner integration the amplification factor becomes $\rho(\lambda^{n,m,p}\delta T^{n,m,p})=\exp(\lambda^{n,m,p}\delta T^{m,n,p})$. For general explicit Runge--Kutta microsimulation of order $q$, the amplification factor is $\rho(\lambda^{n,m,p}\delta T^{n,m,p})=\sum_{s=0}^{q}\frac{(\lambda^{n,m,p}\delta T^{m,n,p})^s}{s!}$. In this article, the forward Euler scheme is used as a micro solver. So the amplification factor for a forward Euler step of micro time step $\delta T^{m,n,p}$ is $\rho(\lambda^{n,m,p}\delta T^{n,m,p})=1+\lambda^{n,m,p}\delta T^{m,n,p}$. However, the following stability study is made for a general numerical microsimulation.
	
	$\bullet$ \textbf{Error amplification in micro time steps:}
	
	Suppose the error at time $T^{n,m,p}$ is an eigencomponent corresponding to the current eigenvalue $\lambda^{n,m,p}$ is $\epsilon^{n,m,p}$. After the inner integration step $\delta T^{n,m,p}$, the error is amplified (usually, that decreases) to
	\begin{equation}\label{eqn:Test_eqn_Micro_Error}
		\epsilon^{n,m,p+1}=\rho^{n,m,p}\epsilon^{n,m,p},\hspace{0.2cm} p=0,1,\ldots,M_{n,m}.
	\end{equation}
	
	$\bullet$ \textbf{Error amplification in meso time steps:}
	
	Suppose the error at meso time level $T^{n,m}$ is $\epsilon^{n,m}$ corresponding to the current eigenvalue $\lambda^{n,m}$. After first $M_{n,m}$ inner steps of sizes $\delta T^{n,m,p}$, the error is amplified to
	\begin{equation}\label{eqn:Micro_Error_at_end}
		\epsilon^{n,m,M_{n,m}}=\prod_{p=0}^{M_{n,m}-1}\rho^{n,m,p}\epsilon^{n,m}.
	\end{equation}
	The projective extrapolation \eqref{eqn:Meso_Projective_Discretisation} from the time $T^{n,m,M_{n,m}+1}$ to $T^{n,m+1}$ amplifies the error to
	\begin{equation}
		\begin{aligned}
			\epsilon^{n,m+1}&=\left[\left(R^{n,m}+1\right)\epsilon^{n,m,M_{n,m}+1}-R^{n,m}\epsilon^{n,m,M_{n,m}}\right]\\
			&=\left[\left\{\left(R^{n,m}+1\right)\rho^{n,m,M_{n,m}}-R^{n,m}\right\}\prod_{p=0}^{M_{n,m}-1}\rho^{n,m,p}\right]\epsilon^{n,m},\hspace{0.2cm} (\text{using equations \eqref{eqn:Test_eqn_Micro_Error} and \eqref{eqn:Micro_Error_at_end}})
		\end{aligned}
	\end{equation}
	where, $m=0$, 1,$\ldots$, $l_n-1$.
	
	Hence, the error amplification in the meso time step $\partial T^{n,m}$ is expressed as
	\begin{equation}\label{eqn:Test_eqn_Meso_Error}
		\epsilon^{n,m+1}=\sigma^{n,m}\epsilon^{n,m},
	\end{equation}
	where the amplification factor in meso step $\partial T^{n,m}$ is
	\begin{equation}\label{eqn:Test_eqn_Meso_Amplification_Factor}
		\sigma^{n,m}:=\left\{\left(R^{n,m}+1\right)\rho^{n,m,M_{n,m}}-R^{n,m}\right\}\prod_{p=0}^{M_{n,m}-1}\rho^{n,m,p},
	\end{equation}
	where $m=0$, 1,$\ldots$, $l_n-1$.
	
	$\bullet$ \textbf{Error amplification in macro time steps:}
	
	Suppose the error at macro time level $T^{n}$ is $\epsilon^{n}$ corresponding to the eigenvalue $\lambda^n$. After applying $l_n$ meso steps in the macro time step $\Delta T^n$, the amplified error in macro time step ($\Delta T^{n}$) is expressed as 
	\begin{equation}
		\epsilon^{n+1}=\prod_{m=0}^{l_n-1}\left[\left\{\left(R^{n,m}+1\right)\rho^{n,m,M_{n,m}}-R^{n,m}\right\}\prod_{p=0}^{M_{n,m}-1}\rho^{n,m,p}\right]\epsilon^{n}.
	\end{equation}
	
	Hence, the error amplification in the macro time step $\Delta T^{n}$ is expressed as
	\begin{equation}\label{eqn:Test_eqn_Macro_Error}
		\epsilon^{n+1}=\sigma^{n}\epsilon^{n},
	\end{equation}
	where the amplification factor in macro time step $\Delta T^n$ is
	\begin{equation}\label{eqn:Test_eqn_Macro_Amplification_Factor}
		\sigma^{n}:=\prod_{m=0}^{l_n-1}\left[\left\{\left(R^{n,m}+1\right)\rho^{n,m,M_{n,m}}-R^{n,m}\right\}\prod_{p=0}^{M_{n,m}-1}\rho^{n,m,p}\right].
	\end{equation}
	
	For perfect inner integration, the amplification factor in the macro time step $\Delta T^n$ becomes
	\begin{equation}
		\sigma^n=\prod_{m=0}^{l_n-1}\left[\left\{\left(R^{n,m}+1\right)\exp(\lambda^{n,m,M_{n,m}}\delta T^{n,m,M_{n,m}})-R^{n,m}\right\}\prod_{p=0}^{M_{n,m}-1}\exp(\lambda^{n,m,p}\delta T^{n,m,p})\right].
	\end{equation}
	
	The \texttt{GPI} method is absolutely stable if $|\sigma^n|\le1$, where the absolute stability depends on the values of $\lambda^{n,m,p}\delta T^{n,m,p}$. 
	
	Suppose that $\Theta=\lambda^{n,m,p}\delta T^{n,m,p}$ be fixed for all $n$, $m$ and $p$. The stability region in the $\Theta$-plane is defined as the set of values of $\Theta$ for which $|\sigma^n(\Theta)|\le1$. In this formulation, the amplification factor $\rho$ can be regarded as a function of $\Theta$.
	The amplification factor over a macro time step $\Delta T^n$ is given by
	\begin{equation}\label{eqn:Stability_Gap_Const_lam_delT}
		\begin{aligned}
			\sigma^{n}&=\prod_{m=0}^{l_n-1}\left[\left\{\left(R^{n,m}+1\right)\rho-R^{n,m}\right\}\prod_{p=0}^{M_{n,m}-1}\rho\right],\\
			&=\prod_{m=0}^{l_n-1}\left\{\left(R^{n,m}+1\right)\rho-R^{n,m}\right\}\rho^{\mathfrak{M}_n}
		\end{aligned}
	\end{equation}
	where $\mathfrak{M}_n:=\sum_{m=0}^{l_n-1}M_{n,m}$ and $\rho:=\rho^{n,m,p}$ be fixed for all $n$, $m$ and $p$.

	\begin{figure}
		\centering
		
		\begin{subfigure}{.32\textwidth}
			\centering
			\includegraphics[width=\linewidth]{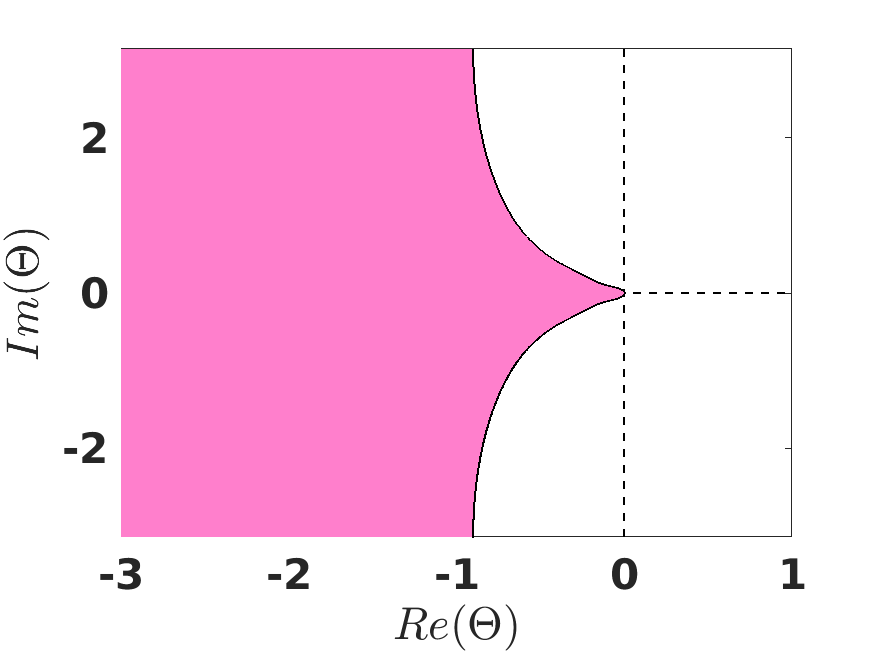}  
			\caption{$R^{n,0}=2$, $R^{n,1}=20$}
			\label{fig:Stability_Perfect_Inner_ln2_mn5_R_2_20}
		\end{subfigure}
		\begin{subfigure}{.32\textwidth}
			\centering
			\includegraphics[width=\linewidth]{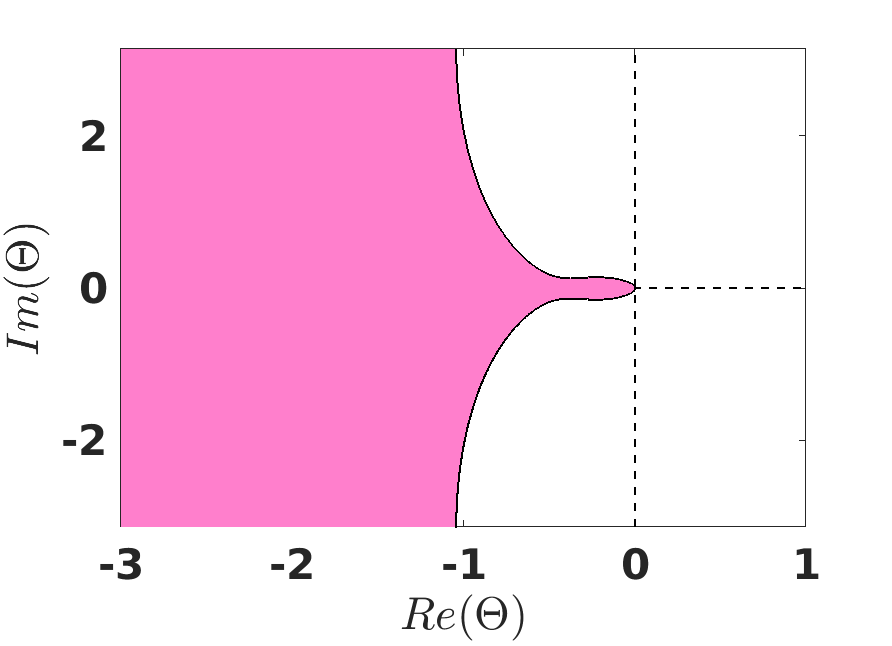}  
			\caption{$R^{n,0}=6$, $R^{n,1}=16$}
			\label{fig:Stability_Perfect_Inner_ln2_mn5_R_6_16}
		\end{subfigure}
		\begin{subfigure}{.32\textwidth}
			\centering
			\includegraphics[width=\linewidth]{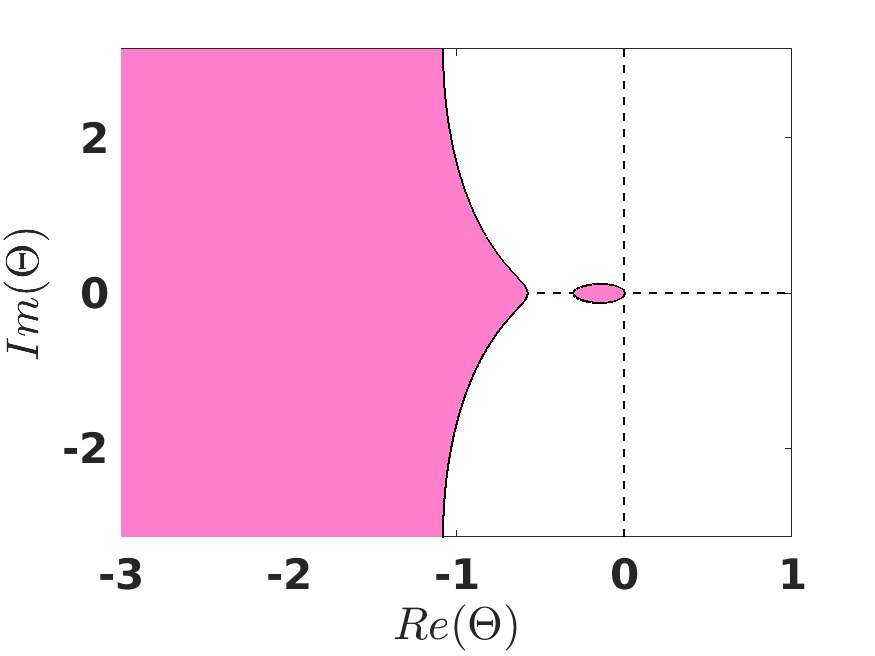}  
			\caption{$R^{n,0}=10$, $R^{n,1}=12$}
			\label{fig:Stability_Perfect_Inner_ln2_mn5_R_10_12}
		\end{subfigure}
		
		\caption{Complex $\Theta$-plane stability for the \texttt{GPI} scheme with perfect inner integrator for $l_n=2$ and $\mathfrak{M}_n=5$. The values of $R^{n,m}$, for $m=0,1$, are considered in three cases: (a) $\{2,20\}$, (b) $\{6,16\}$ and (c) $\{10,12\}$.}
		\label{fig:Stability_Perfect_Inner}
	\end{figure}
	
	The stability region depends on the choice of the inner integrator as well as on the parameters $R^{n,m}$, $M_{n,m}$ or $\mathfrak{M}_n$ (the total number of micro time steps in the macro time step) and $l_n$.
	
	For the case of perfect inner integration, the stability region of the \texttt{GPI} scheme is illustrated in Figure \ref{fig:Stability_Perfect_Inner}. In this setting, two meso steps ($l_n=2$) and a total of $\mathfrak{M}_n=5$ micro steps are considered within macro time steps $\Delta T^n$. The values of $R^{n,m}$ are distributed in three different ways, summing to 22, such as $\{2, 20\}$, $\{6, 16\}$ and $\{10, 12\}$. The corresponding stability regions are depicted within the strip $Im(\Theta)\in[-\pi,\pi]$, which is periodically repeated with period $2\pi$ along the imaginary axis. For the first two distributions, the $\Theta$-plane is divided into two regions. In contrast, for the third distribution, the $\Theta$-plane is partitioned into three regions. In particular, the stability region splits into multiple disconnected components, resulting in a gap along the negative real axis in the $\Theta$-plane where the method remains stable.
	
	%%%%%%%%%%%%%%%%%%%%%%%%%%%%%%%%%%%%%%%%%%%%%%%%%%%%%%
	\begin{figure}
		\centering
		\begin{subfigure}{.32\textwidth}
			\centering
			% include third image
			\includegraphics[width=1\linewidth]{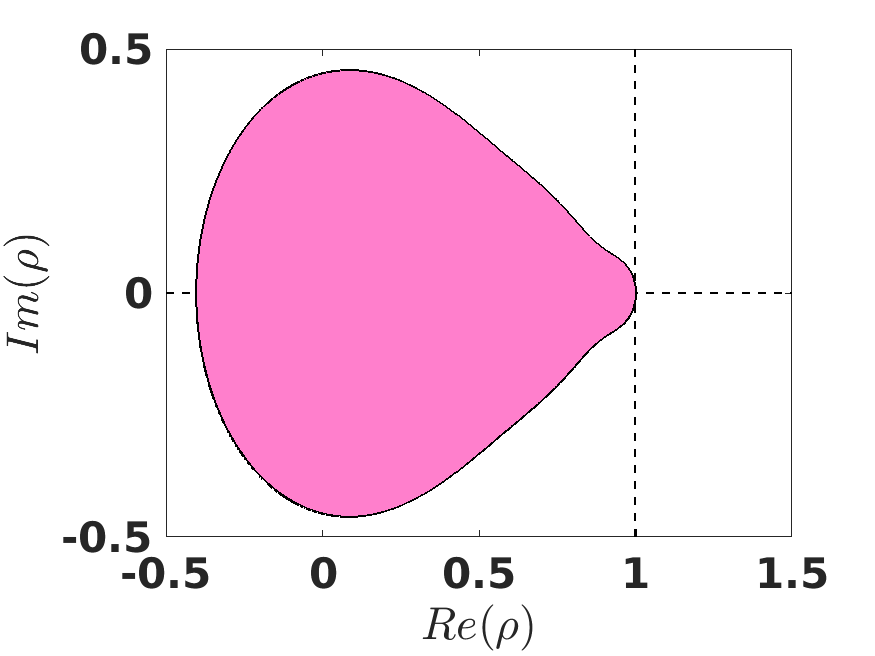}  
			\caption{$R^{n,0}=2$, $R^{n,1}=20$}
			\label{fig:Stability_NonPerfect_Inner_ln2_mn5_R_2_20}
		\end{subfigure}
		\begin{subfigure}{.32\textwidth}
			\centering
			% include fourth image
			\includegraphics[width=1\linewidth]{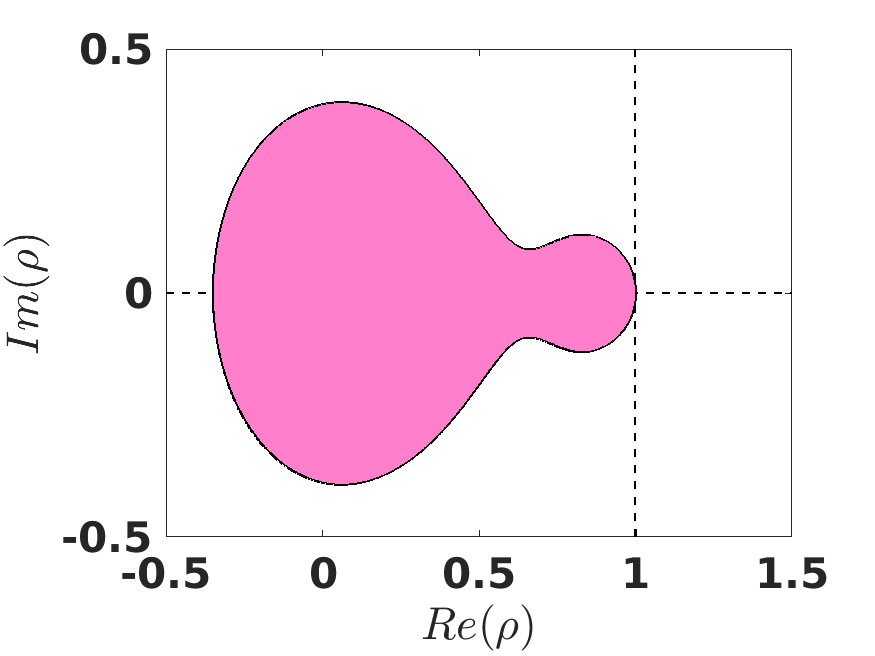}  
			\caption{$R^{n,0}=6$, $R^{n,1}=16$}
			\label{fig:Stability_NonPerfect_Inner_ln2_mn5_R_6_16}
		\end{subfigure}
		\begin{subfigure}{0.32\textwidth}
			\centering
			% include fourth image
			\includegraphics[width=1\linewidth]{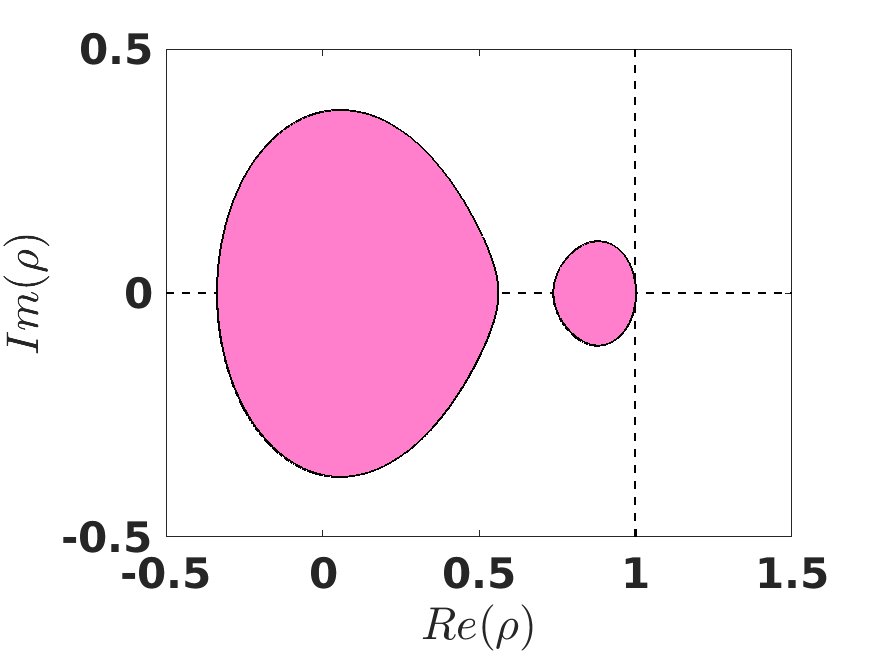}  
			\caption{$R^{n,0}=10$, $R^{n,1}=12$}
			\label{fig:Stability_NonPerfect_Inner_ln2_mn5_R_10_12}
		\end{subfigure}
		\caption{Complex $\rho$-plane stability for the \texttt{GPI} scheme for $l_n=2$ and $\mathfrak{M}_n=5$. The values of $R^{n,m}$, for $m=0$, 1 are considered into three possible ways as (a) $\{2, 20\}$, (b) $\{6, 16\}$ and (c) $\{10, 12\}$.}
		\label{fig:Stability_NonPerfect_Inner}
	\end{figure}
	
	The perfect inner integrator is applicable only to a very limited class of problems. In general, one must employ numerical schemes, for which the amplification factor over the micro time steps provides an approximation to that of the perfect inner integration. For such numerical inner integrators, arbitrarily large values of $\Theta$ can not be realised, except for trivial cases. Consequently, the stability region in the $\Theta$-plane is no longer unbounded.
	
	As illustrated in Figure \ref{fig:Stability_Perfect_Inner}, for the same set of parameter values, the stability region of the \texttt{GPI} scheme corresponding to numerical (non-perfect) inner integration is shown in Figure \ref{fig:Stability_NonPerfect_Inner}. In this case, the stability region becomes finite, since numerical integrators are unstable for sufficiently large values of $\Theta$. Nevertheless, a similar pattern in the splitting of the stability region is observed. 
	
	To determine the stability region for different types of inner integrators, one must map the stability region from the $\rho$-plane to the $\Theta$-plane using the specific form of $\rho(\Theta)$ associated with the chosen inner integrator. In this sense, Figure \ref{fig:Stability_Perfect_Inner} can be interpreted as the logarithmic mapping of Figure \ref{fig:Stability_NonPerfect_Inner}.
	
	%%%%%%%%%%%%%%%%%%%%%%%%%%%%%%%%%%%%%%%%%%%%%%%%%%%%
	
	\begin{figure}
		\centering
		\begin{subfigure}{.32\textwidth}
			\centering
			% include third image
			\includegraphics[width=1\linewidth]{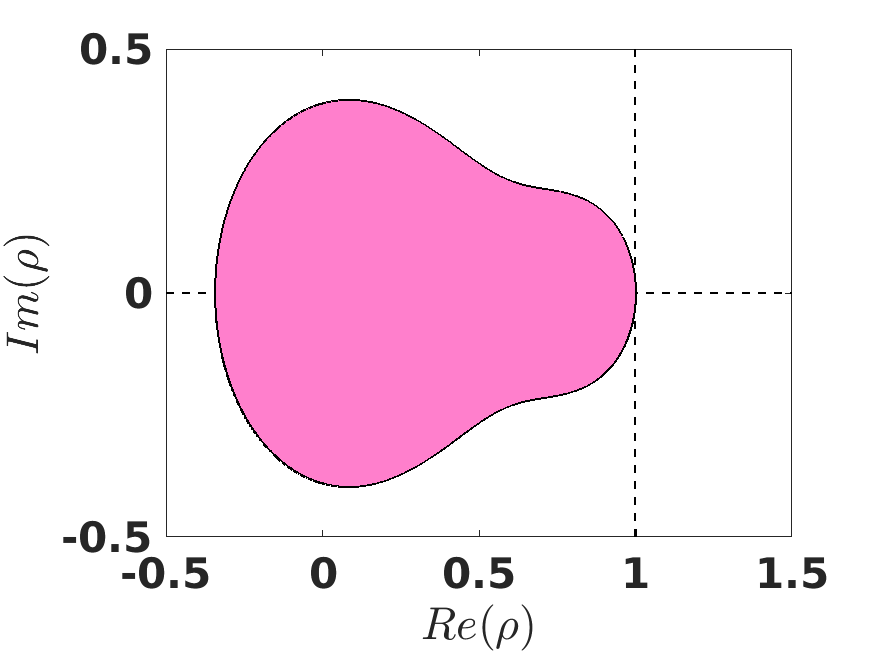}  
			\caption{uniform $R^{n,m}=6$}
			\label{fig:Stability_lm4_Mn8_Uniform_R6}
		\end{subfigure}
		\begin{subfigure}{.32\textwidth}
			\centering
			% include fourth image
			\includegraphics[width=1\linewidth]{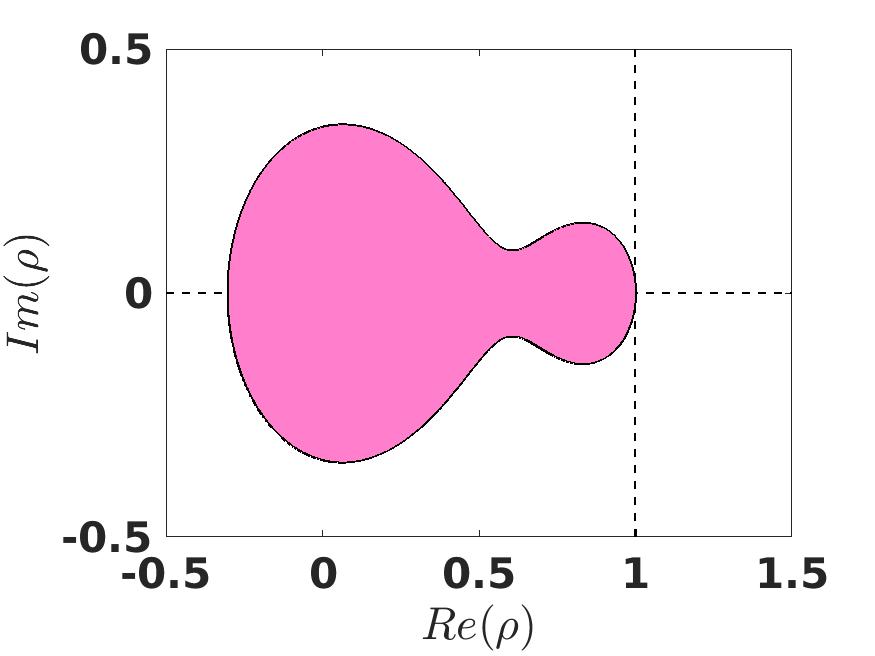}  
			\caption{Uniform $R^n=8$}
			\label{fig:Stability_lm4_Mn8_Uniform_R8}
		\end{subfigure}
		\begin{subfigure}{0.32\textwidth}
			\centering
			% include fourth image
			\includegraphics[width=1\linewidth]{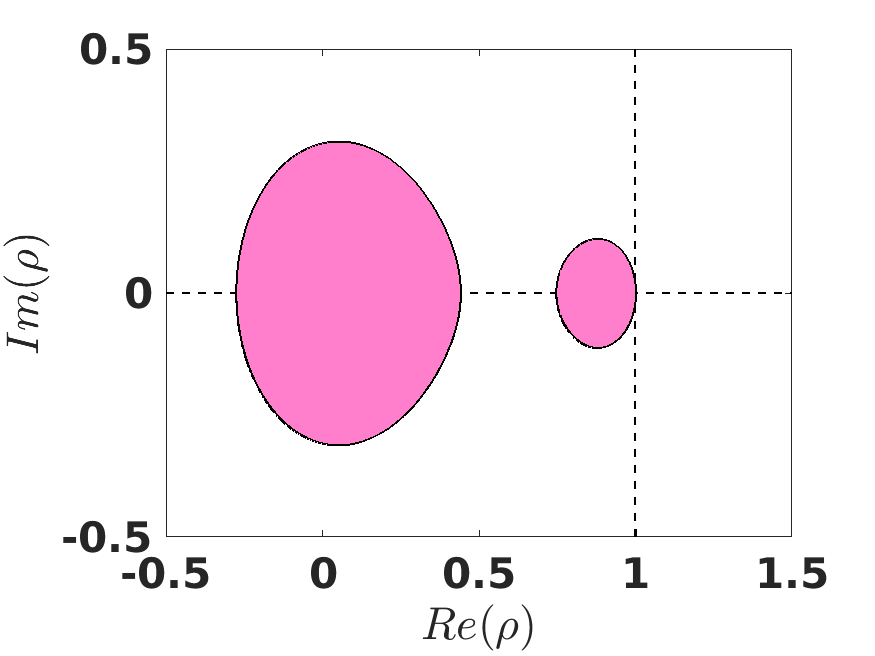}  
			\caption{Uniform $R^{n,m}=10$}
			\label{fig:Stability_lm4_Mn8_Uniform_R10}
		\end{subfigure}
		\caption{Complex $\rho$-plane stability for the \texttt{GPI} scheme for $l_n=4$, $\mathfrak{M}_n=8$. The values of uniform $R^{n,m}$, for $m=0$, 1, 2 and 3 are considered into three possible ways as (a) 6, (b) 8 and (c) 10.}
		\label{fig:Stability_lm4_Mn8_Uniform_R}
	\end{figure}
	
	Figure \ref{fig:Stability_lm4_Mn8_Uniform_R} illustrates the stability regions as the parameter $R^{n,m}$ increases. In this case, uniform values of $R^{n,m}$ are considered with $l_n=4$ and $\mathfrak{M}_n=8$. The uniform values of $R^{n,m}$ used in Subfigures \ref{fig:Stability_lm4_Mn8_Uniform_R6}, \ref{fig:Stability_lm4_Mn8_Uniform_R8} and \ref{fig:Stability_lm4_Mn8_Uniform_R10} are 6, 8 and 10, respectively. The results indicate that, as $R^{n,m}$ increases, the initially connected stability region gradually shrinks. Beyond a certain threshold value of $R^{n,m}$, this single connected region splits into two disconnected components.
	
	%%%%%%%%%%%%%%%%%%%%%%%%%%%%%%%%%%%%%%%%%%%%%%%%%%%
	
	\begin{figure}
		\begin{subfigure}{.32\textwidth}
			\centering
			% include third image
			\includegraphics[width=1\linewidth]{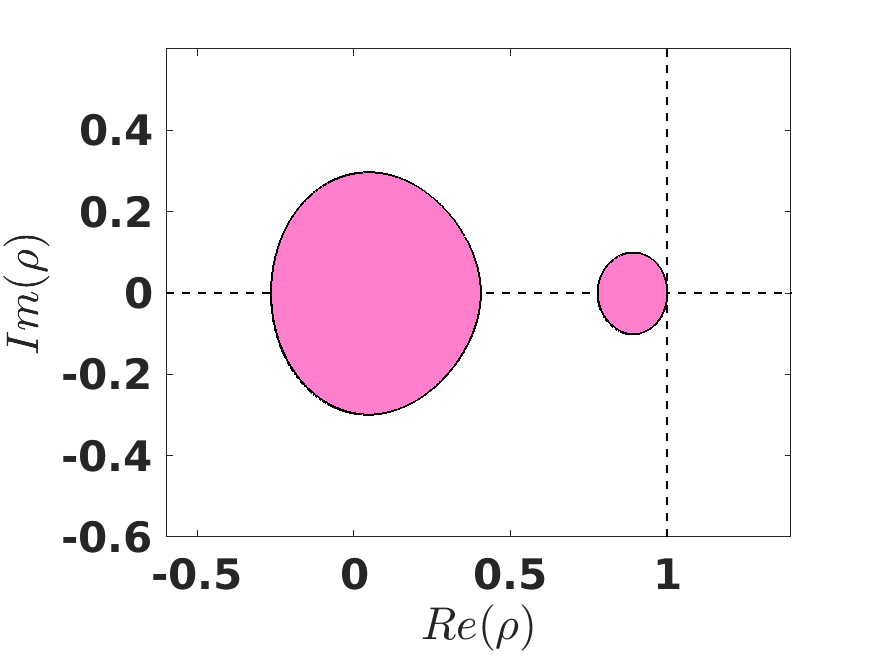}  
			\caption{$\mathfrak{M}_n=4$}
			\label{fig:Stability_lm2_R_10_12_Uniform_Mn4}
		\end{subfigure}
		\begin{subfigure}{.32\textwidth}
			\centering
			% include fourth image
			\includegraphics[width=1\linewidth]{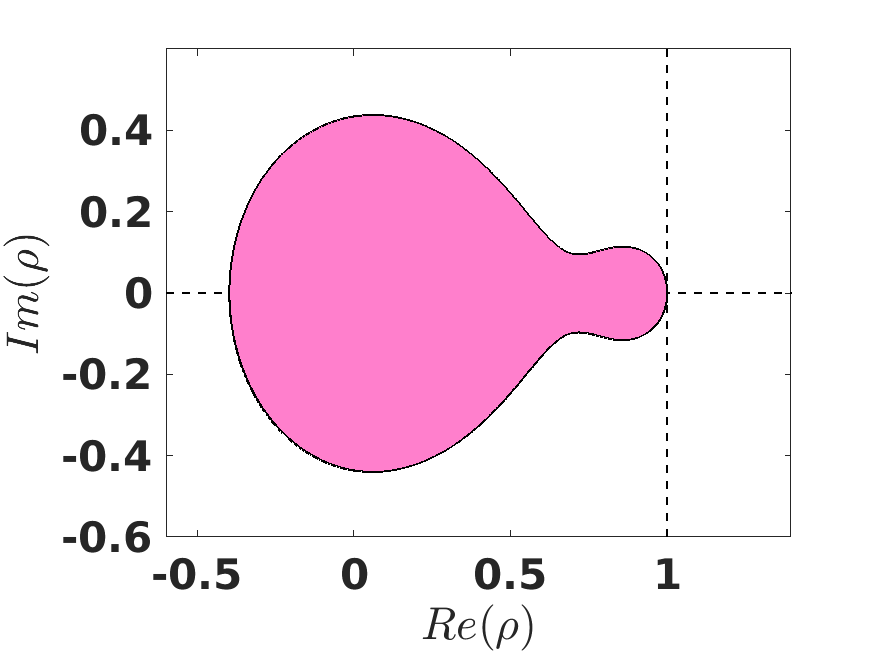}  
			\caption{$\mathfrak{M}_n=6$}
			\label{fig:Stability_lm2_R_10_12_Uniform_Mn6}
		\end{subfigure}
		\begin{subfigure}{0.32\textwidth}
			\centering
			% include fourth image
			\includegraphics[width=1\linewidth]{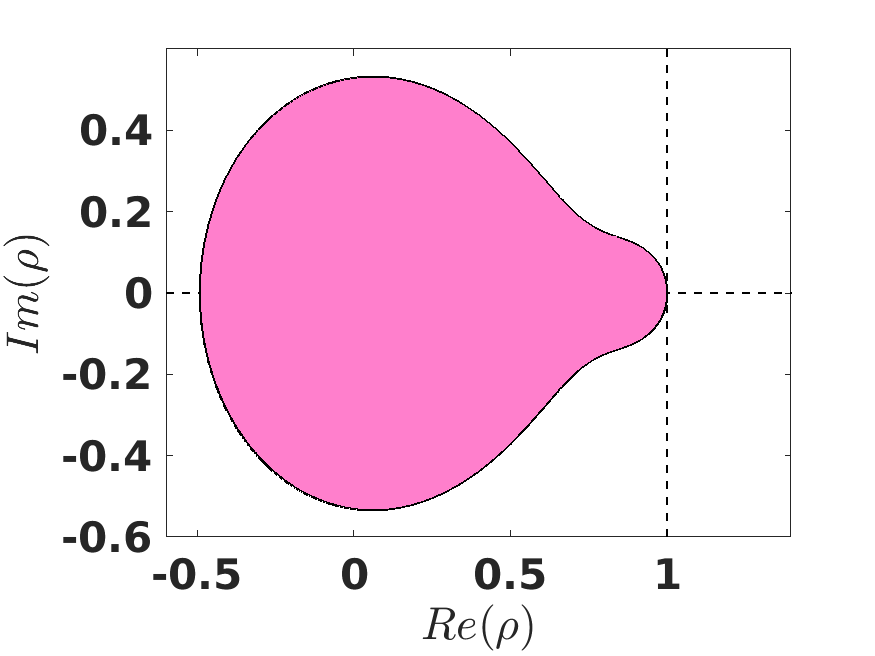}  
			\caption{$\mathfrak{M}_n=8$}
			\label{fig:Stability_lm2_R_10_12_Uniform_Mn8}
		\end{subfigure}
		\caption{Complex $\rho$-plane stability for the \texttt{GPI} scheme for $l_n=2$, $R^{n,0}=10$, $R^{n,1}=12$. The values of $\mathfrak{M}_n$ are considered in three possible ways: (a) 4, (b) 6 and (c) 8.}
		\label{fig:Stability_FE_Inner_Vari_Mn}
	\end{figure}
	
	Figure \ref{fig:Stability_FE_Inner_Vari_Mn} illustrates the stability regions for varying values of $\mathfrak{M}_n$. In this case, the parameters are fixed as $l_n=2$, $R^{n,0}=10$ and $R^{n,1}=12$, while $\mathfrak{M}_n$ takes the values 4, 6 and 8 in Subfigures \ref{fig:Stability_lm2_R_10_12_Uniform_Mn4}, \ref{fig:Stability_lm2_R_10_12_Uniform_Mn6} and \ref{fig:Stability_lm2_R_10_12_Uniform_Mn8}, respectively. The results show that, as the number of micro time steps within each macro time step increases, the stability region expands, transitioning from multiple disconnected components to a single connected region. 
	
	Furthermore, by comparing Figures \ref{fig:Stability_NonPerfect_Inner}, \ref{fig:Stability_lm4_Mn8_Uniform_R} and \ref{fig:Stability_FE_Inner_Vari_Mn}, it is evident that the splitting of the stability region depends explicitly on both $R^{n,m}$ and $\mathfrak{M}_n$. The tendency for the stability region to split increases with $R^{n,m}$ and decreases with $\mathfrak{M}_n$. That is, as $R^{n,m}$ is larger than $\gamma\mathfrak{M}_n$ would cause the stability region to break into two pieces, where $\gamma$ represents a splitting parameter.

	%\begin{tikzpicture}
	%	
	%	\node (plot) at (0,0)
	%	{\includegraphics[width=6cm]{L2_K8_M_10_12_FE_Kvari.eps}};
	%	
	%	\draw[->,thick] (-4,0) -- (-1,0);
	%	\node at (-4.5,0) {Input};
	%	
	%	\draw[->,thick] (3,0) -- (5,0);
	%	\node at (5.5,0) {Output};
	%	
	%\end{tikzpicture}
	
	%%%%%%%%%%%%%%%%%%%%%%%%%%%%%%%%%%%%%%%%%%%%%%%%%
	
	\subsection{Stability region gap in \texttt{GPI} scheme}
	In this subsection, the objective is to analyse the splitting of the stability region. This can be investigated by examining \eqref{eqn:Stability_Gap_Const_lam_delT} and determining when the locus of $\rho$ intersects the real axis as $\sigma^n$ traverses the unit circle. Such an intersection occurs only when $\sigma^n$ is real, that is, when $\sigma^n=\pm 1$. 
	Consequently, equation \eqref{eqn:Stability_Gap_Const_lam_delT} admits real roots $\rho$ only for these two equations. 
	
	To observe the splitting of the stability region, $\sigma^n=\pm 1$ must attain at least four distinct real roots. The analysis is carried out separately for uniform and non-uniform choices of $R^{n,m}$ in Subsections \ref{subsubsec:Uni_R_nm} and \ref{subsubsec:Non_Uni_R_nm}, respectively.
	
	\subsubsection{The values of $R^{n,m}$ are uniform in macro time step}\label{subsubsec:Uni_R_nm}
	
	Let the values of $R^{n,m}$ are uniform in macro time step $\Delta T^n$ such that $R^{n,m}=R^n$, $\forall m=0, 1, \ldots, l_n-1$. That is 
	\begin{center}
		extrapolation step of $m^{th}$ meso step $\propto$ last micro time step of $m^{th}$ meso step,
	\end{center}
	for all $m$. The amplification factor becomes
	\begin{equation}\label{eqn:Stability_Uniform_Amplification_Factor}
		\sigma^{n}=\left\{\left(R^n+1\right)\rho-R^n\right\}^{l_n}\rho^{\mathfrak{M}_n}.
	\end{equation}
	To find the multiple real roots of $\sigma^n=\pm 1$ depends on the number of meso time steps $l_n$. That is obtained when 
	\begin{equation}\label{eqn:Stability_Amplification_Boundary}
		\sigma^n=\pm (-1)^{l_n}.
	\end{equation}
	That is 
	\begin{equation}\label{eqn:Stability_Uniform_Root_Finding_Equations_1}
		\left\{\left(R^n+1\right)\rho-R^n\right\}^{l_n}\rho^{\mathfrak{M}_n}=(-1)^{l_n},
	\end{equation}
	and
	\begin{equation}\label{eqn:Stability_Uniform_Root_Finding_Equations_2}
		\left\{\left(R^n+1\right)\rho-R^n\right\}^{l_n}\rho^{\mathfrak{M}_n}=(-1)^{l_n+1}.
	\end{equation}
	
	Equations \eqref{eqn:Stability_Uniform_Root_Finding_Equations_1} and \eqref{eqn:Stability_Uniform_Root_Finding_Equations_2} can be rewritten as
	\begin{equation}\label{eqn:Stability_Uniform_Root_Finding_Equations_LR1}
		\frac{1}{\rho^{\mathfrak{M}_n}}=(-1)^{l_n}\left\{\left(R^n+1\right)\rho-R^n\right\}^{l_n},
	\end{equation}
	and
	\begin{equation}\label{eqn:Stability_Uniform_Root_Finding_Equations_LR2}
		\frac{1}{\rho^{\mathfrak{M}_n}}=(-1)^{l_n+1}\left\{\left(R^n+1\right)\rho-R^n\right\}^{l_n}.
	\end{equation}

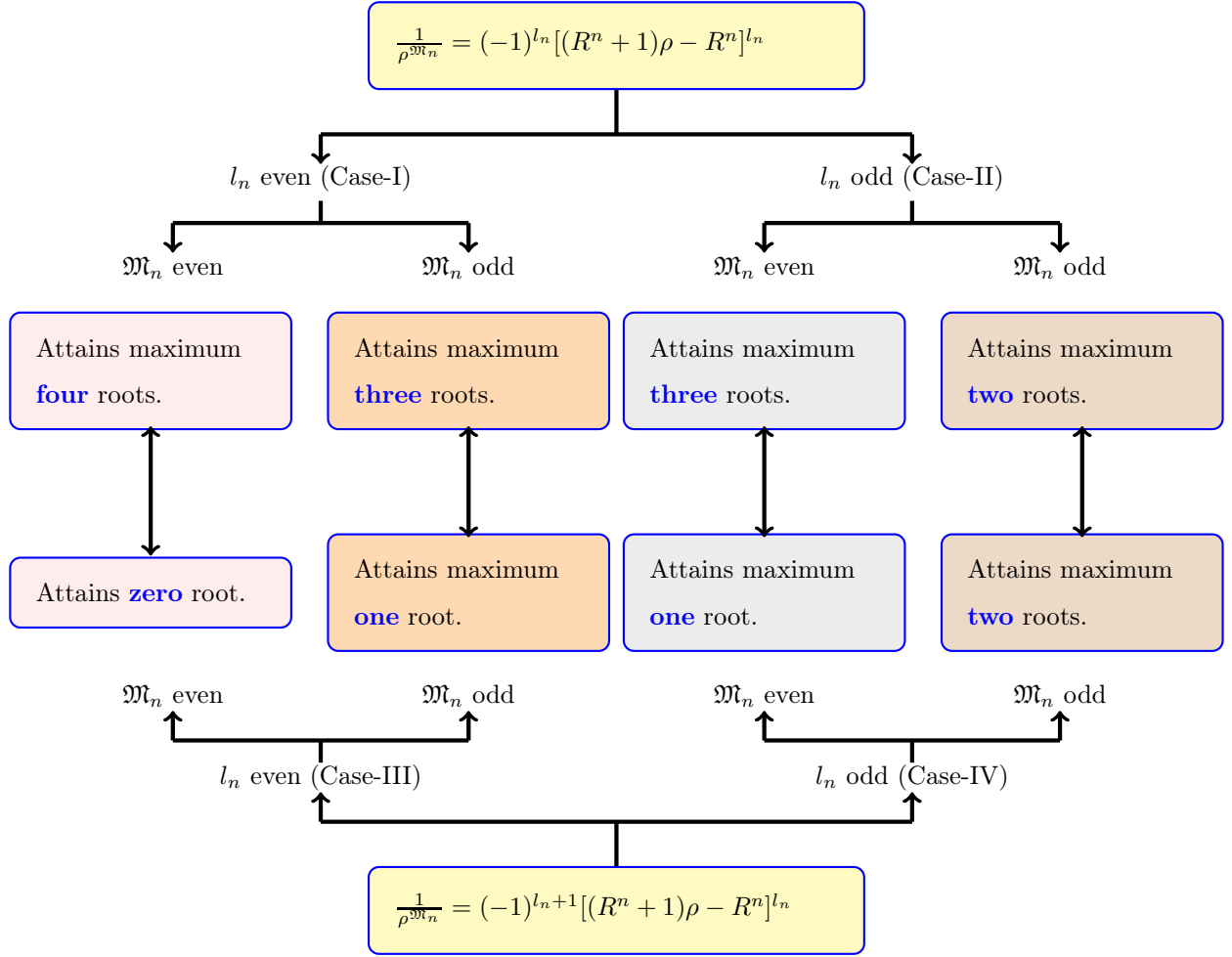
\begin{figure}
		\centering
		\begin{tikzpicture}[scale=1]
			\node at (0,5.4) [
			draw=blue,
			fill=yellow!30,
			thick,
			rounded corners,
			inner sep=10pt,
			text width=6cm,
			align=left
			]{$\frac{1}{\rho^{\mathfrak{M}_n}}=(-1)^{l_n}[(R^n+1)\rho-R^n]^{l_n}$
			};
			
			% First step
			\draw[ultra thick,] (0,4.8) -- (0,4.2);
			\draw[ultra thick,] (-4,4.2) -- (4,4.2);
			
			\draw[ultra thick,->] (-4,4.2) -- (-4,3.8);
			\draw[ultra thick,->] (4,4.2) -- (4,3.8);
			
			\draw[very thick] (-4,3.6)  node {$l_n$ even (Case-I)};
			\draw[very thick] (4,3.6)  node {$l_n$ odd (Case-II)};
			
			% Second step
			\draw[ultra thick,] (-4,3.3) -- (-4,3);
			\draw[ultra thick,] (-6,3) -- (-2,3);
			\draw[ultra thick,->] (-6,3) -- (-6,2.6);
			\draw[ultra thick,->] (-2,3) -- (-2,2.6);
			\draw[very thick] (-6,2.4)  node {$\mathfrak{M}_n$ even};
			\draw[very thick] (-2,2.4)  node {$\mathfrak{M}_n$ odd};

			\draw[ultra thick,] (4,3.3) -- (4,3);
			\draw[ultra thick,] (2,3) -- (6,3);
			\draw[ultra thick,->] (6,3) -- (6,2.6);
			\draw[ultra thick,->] (2,3) -- (2,2.6);
			\draw[very thick] (6,2.4)  node {$\mathfrak{M}_n$ odd};
			\draw[very thick] (2,2.4)  node {$\mathfrak{M}_n$ even};
			
			% Third step
			\node at (-6.3,1) [
			draw=blue,
			fill=pink!30,
			thick,
			rounded corners,
			inner sep=10pt,
			text width=3.1cm,
			align=left
			]{Attains maximum \textbf{\textcolor{blue}{four}} roots.
			};

			\node at (-2,1) [
			draw=blue,
			fill=orange!30,
			thick,
			rounded corners,
			inner sep=10pt,
			text width=3.1cm,
			align=left
			]{Attains maximum \textbf{\textcolor{blue}{three}} roots.
			};
			
			\node at (2,1) [
			draw=blue,
			fill=lightgray!30,
			thick,
			rounded corners,
			inner sep=10pt,
			text width=3.1cm,
			align=left
			]{Attains maximum \textbf{\textcolor{blue}{three}} roots.
			};
			
			\node at (6.3,1) [
			draw=blue,
			fill=brown!30,
			thick,
			rounded corners,
			inner sep=10pt,
			text width=3.1cm,
			align=left
			]{Attains maximum \textbf{\textcolor{blue}{two}} roots.
			};

			%%%%%%%%%%%%%%%%%%%%%%%%%%%%%%%%%%%%%%%%
			% Next part
			
			% Third step from bottom
			\node at (-6.3,-2) [
			draw=blue,
			fill=pink!30,
			thick,
			rounded corners,
			inner sep=10pt,
			text width=3.1cm,
			align=left
			]{Attains \textbf{\textcolor{blue}{zero}} root.
			};

			\node at (-2,-2) [
			draw=blue,
			fill=orange!30,
			thick,
			rounded corners,
			inner sep=10pt,
			text width=3.1cm,
			align=left
			]{Attains maximum \textbf{\textcolor{blue}{one}} root.
			};
			
			\node at (2,-2) [
			draw=blue,
			fill=lightgray!30,
			thick,
			rounded corners,
			inner sep=10pt,
			text width=3.1cm,
			align=left
			]{Attains maximum \textbf{\textcolor{blue}{one}} root.
			};
			
			\node at (6.3,-2) [
			draw=blue,
			fill=brown!30,
			thick,
			rounded corners,
			inner sep=10pt,
			text width=3.1cm,
			align=left
			]{Attains maximum \textbf{\textcolor{blue}{two}} roots.
			};
			
			% Second step
			\draw[ultra thick,] (-4,-4.3) -- (-4,-4);
			\draw[ultra thick,] (-6,-4) -- (-2,-4);
			\draw[ultra thick,->] (-6,-4) -- (-6,-3.6);
			\draw[ultra thick,->] (-2,-4) -- (-2,-3.6);
			\draw[very thick] (-6,-3.4)  node {$\mathfrak{M}_n$ even};
			\draw[very thick] (-2,-3.4)  node {$\mathfrak{M}_n$ odd};

			\draw[ultra thick,] (4,-4.3) -- (4,-4);
			\draw[ultra thick,] (2,-4) -- (6,-4);
			\draw[ultra thick,->] (6,-4) -- (6,-3.6);
			\draw[ultra thick,->] (2,-4) -- (2,-3.6);
			\draw[very thick] (6,-3.4)  node {$\mathfrak{M}_n$ odd};
			\draw[very thick] (2,-3.4)  node {$\mathfrak{M}_n$ even};
			
			% First step
			\draw[ultra thick,] (0,-5.7) -- (0,-5.1);
			\draw[ultra thick,] (-4,-5.1) -- (4,-5.1);
			
			\draw[ultra thick,->] (-4,-5.1) -- (-4,-4.7);
			\draw[ultra thick,->] (4,-5.1) -- (4,-4.7);
			
			\draw[very thick] (-4,-4.5)  node {$l_n$ even (Case-III)};
			\draw[very thick] (4,-4.5)  node {$l_n$ odd (Case-IV)};
			
			\node at (0,-6.3) [
			draw=blue,
			fill=yellow!30,
			thick,
			rounded corners,
			inner sep=10pt,
			text width=6cm,
			align=left
			]{$\frac{1}{\rho^{\mathfrak{M}_n}}=(-1)^{l_n+1}[(R^n+1)\rho-R^n]^{l_n}$
			};
			
			% Connect the boxes
			\draw[ultra thick,<->] (-6.3,-1.5) -- (-6.3,0.25);
			\draw[ultra thick,<->] (-2,-1.25) -- (-2,0.25); %6.3,1
			\draw[ultra thick,<->] (2,-1.25) -- (2,0.25);
			\draw[ultra thick,<->] (6.3,-1.25) -- (6.3,0.25);
			
		\end{tikzpicture}
		%	\label{fig:Schematic_Stability_Uniform}
		\caption{A schematic illustrating the maximum number of distinct real roots of equations \eqref{eqn:Stability_Uniform_Root_Finding_Equations_LR1} and \eqref{eqn:Stability_Uniform_Root_Finding_Equations_LR2} for various combinations of $l_n$ and $\mathfrak{M}_n$.}
		\label{fig:Schematic_Stability_Chart_Roots}
	\end{figure}

	For even $\mathfrak{M}_n$, the curve $\frac{1}{\rho^{\mathfrak{M}_n}}$ lies in the first and second quadrants, whereas for odd $\mathfrak{M}_n$, it lies in the first \& third quadrants. In both cases, the curve has a horizontal asymptote at $y=0$ and a vertical asymptote at $x=0$. Here, we need to find the possible maximum number of real, distinct roots of both the equations \eqref{eqn:Stability_Uniform_Root_Finding_Equations_LR1} and \eqref{eqn:Stability_Uniform_Root_Finding_Equations_LR2}. When both side curves intersect for each individual equations, these are the possible roots.
	\begin{itemize}
		\item \textbf{Case-I: $l_n$ is even in equation \eqref{eqn:Stability_Uniform_Root_Finding_Equations_LR1}}
		
		Suppose,
		\begin{center}
			$y=[(R^n+1)\rho-R^n]^{l_n}$,
		\end{center}
		which decreases for $\rho<\frac{R^n}{R^n+1}$ and increases for $\rho>\frac{R^n}{R^n+1}$. That is $y$ attains its minimum value zero at $\rho=\frac{R^n}{R^n+1}$ and $y>0$ for all $\rho$. For $\rho\rightarrow\infty$, the solution tends to $y\rightarrow\infty$. For any $R^n>0$ and $\mathfrak{M}_n\in\mathbb{N}$, $\rho=1$ is a root of the equation \eqref{eqn:Stability_Uniform_Root_Finding_Equations_LR1}. As $R^n$ increases from zero, minimum point $\frac{R^n}{R^n+1}$ moves toward $\rho=1$ and becomes stepper \& sharper. As $R^n$ increases, for certain $R^n=R^n_\text{critical}$, $y$ touches $\frac{1}{\rho^{\mathfrak{M}_n}}$ in the first quadrant for any $\mathfrak{M}_n$. For $R<R_\text{critical}$, there are no positive real roots for equation \eqref{eqn:Stability_Uniform_Root_Finding_Equations_LR1}. For $R^n>R^n_\text{critical}$, there are two positive real roots for equation \eqref{eqn:Stability_Uniform_Root_Finding_Equations_LR1}. There is one negative real root for any $R^n>0$ if $\mathfrak{M}_n$ is even and none for $\mathfrak{M}_n$ is odd. So, for even and odd $\mathfrak{M}_n$, in total there are a maximum of four and three real roots available, respectively, for the equation \eqref{eqn:Stability_Uniform_Root_Finding_Equations_LR1}. A schematic representation is shown in Figure \ref{fig:Schematic_Stability_Chart_Roots}, with pink and orange coloured boxes in the upper part.
		
		\item \textbf{Case-II: $l_n$ is odd in equation \eqref{eqn:Stability_Uniform_Root_Finding_Equations_LR1}}
		
		Suppose,
		\begin{center}
			$y=-[(R^n+1)\rho-R^n]^{l_n}$,
		\end{center}
		which decreases for all $\rho$. The function $y$ is positive for $\rho<\frac{R^n}{R^n+1}$, negative for $\rho>\frac{R}{R+1}$. As $\rho$ approaches $\pm\infty$, $y$ tends to $\mp\infty$. As $R^n$ increases from zero, for certain $R^n=R^n_\text{critical}$, $y$ touches $\frac{1}{\rho^{\mathfrak{M}_n}}$ in the first quadrant for any $\mathfrak{M}_n\in\mathbb{N}$. For $R^n<R^n_\text{critical}$, there are no positive real roots for equation \eqref{eqn:Stability_Uniform_Root_Finding_Equations_LR1}. For $R^n>R^n_\text{critical}$, there are two positive real roots of \eqref{eqn:Stability_Uniform_Root_Finding_Equations_LR1}. There is one negative real root for any $R^n>0$ if $\mathfrak{M}_n$ is even and none for $\mathfrak{M}_n$ is odd. So, for even and odd $\mathfrak{M}_n$, there are a maximum of three and two real roots available, respectively, for the equation \eqref{eqn:Stability_Uniform_Root_Finding_Equations_LR1}. A schematic representation is shown in Figure \ref{fig:Schematic_Stability_Chart_Roots}, with light-grey and brown coloured boxes in the upper part.

		\item \textbf{Case-III: $l_n$ is even in equation \eqref{eqn:Stability_Uniform_Root_Finding_Equations_LR2}}
		
		Suppose,
		\begin{center}
			$y=-[(R^n+1)\rho-R^n]^{l_n}$,
		\end{center}
		is a downward-opening, even-power curve, symmetric about the vertical line $\rho=\frac{R^n}{R^n+1}$, where the curve achieves maximum value zero. As $|\rho|$ approaches to $\infty$, $y$ tends to $-\infty$. So $y$ always lies in the lower half plane. For even $\mathfrak{M}_n$, $\frac{1}{\rho^{\mathfrak{M}_n}}$ always lies on the upper half plane, so it never intersects $y$. So for even $\mathfrak{M}_n$, there are no real roots present of the equation \eqref{eqn:Stability_Uniform_Root_Finding_Equations_LR2}. However, for odd $\mathfrak{M}_n$, $\frac{1}{\rho^{\mathfrak{M}_n}}$ and $y$ intersect once in the third quadrant. So there is one real negative root. For even and odd $\mathfrak{M}_n$, there are maximum zero and one real roots available, respectively, for the equation \eqref{eqn:Stability_Uniform_Root_Finding_Equations_LR2}. A schematic representation is shown in Figure \ref{fig:Schematic_Stability_Chart_Roots}, with light-grey and brown coloured boxes in the lower part.
		
		\item \textbf{Case-IV: $l_n$ is odd in equation \eqref{eqn:Stability_Uniform_Root_Finding_Equations_LR2}}
		
		Suppose,
		\begin{center}
			$y=[(R^n+1)\rho-R^n]^{l_n}$,
		\end{center}
		which is a monotonically increasing polynomial, becomes flat near $\rho=\frac{R^n}{R^n+1}$ for $l_n>1$ and a rapid growth is seen for large $|\rho|$. For $\rho>\frac{R^n}{R^n+1}$, $y$ lies in upper half plane, whereas for $\rho<\frac{R^n}{R^n+1}$, $y$ lies in lower half plane. As $\rho$ approches to $\pm\infty$, $y$ also tends to $\pm\infty$. For any $\mathfrak{M}_n\in\mathbb{N}$, $y$ intersect $\frac{1}{\rho^{\mathfrak{M}_n}}$ at (1,1) in the first quadrant. For even $\mathfrak{M}_n$, (1,1) point is the only root of the equation \eqref{eqn:Stability_Uniform_Root_Finding_Equations_LR2}, whereas for odd $\mathfrak{M}_n$, $\frac{1}{\rho^{\mathfrak{M}_n}}$ and $y$ intersect once in the third quadrant. Therefore, for even and odd $\mathfrak{M}_n$, there are maximum one and two real roots available, respectively, for the equation \eqref{eqn:Stability_Uniform_Root_Finding_Equations_LR2}. A schematic short representation is shown in Figure \ref{fig:Schematic_Stability_Chart_Roots}, coloured by light grey and brown in the lower part.
	\end{itemize}
	The number of real individual roots of the equations \eqref{eqn:Stability_Uniform_Root_Finding_Equations_LR1} and \eqref{eqn:Stability_Uniform_Root_Finding_Equations_LR2} are presented individually for the even and odd cases of both $l_n$ and $\mathfrak{M}_n$ in Figure \ref{fig:Schematic_Stability_Chart_Roots}. The root counts are indicated using distinct colours corresponding to identical combinations of $l_n$ and $\mathfrak{M}_n$ for both the equations \eqref{eqn:Stability_Uniform_Root_Finding_Equations_LR1} and \eqref{eqn:Stability_Uniform_Root_Finding_Equations_LR2}. By summing the number of roots associated with the same cases (i.e., the same coloured entries), we observe that the total number of distinct real roots of equations \eqref{eqn:Stability_Uniform_Root_Finding_Equations_LR1} and \eqref{eqn:Stability_Uniform_Root_Finding_Equations_LR2} is at most four. This implies that, for uniform $R^{n,m}$, a single connected stability region can split into at most two disconnected components.
	
	We have observed in Cases I and II that a critical value $R^n_\text{critical}$ arises from the equation \eqref{eqn:Stability_Uniform_Root_Finding_Equations_LR1}. As $R^n$ increases from zero and for $R^n<R^n_\text{critical}$, equation \eqref{eqn:Stability_Uniform_Root_Finding_Equations_LR1} admits at most two and one real roots for even and odd $\mathfrak{M}_n$, respectively, in the Case-I and at most one and no real roots, respectively, in the Case-II. At the critical value $R^n=R^n_\text{critical}$, the curves $\frac{1}{\rho^{\mathfrak{M}_n}}$ and $(-1)^{l_n}[(R^n+1)\rho-R^n]^{l_n}$ are tangent to each other, yielding a repeated real root. For $R^n>R^n_\text{critical}$, two distinct real roots emerge, which are responsible for the separation of the stability region. Hence, the onset of the splitting of the stability region occurs precisely at $R^n=R^n_\text{critical}$, where the two curves $\frac{1}{\rho^{\mathfrak{M}_n}}$ and $(-1)^{l_n}[(R^n+1)\rho-R^n]^{l_n}$ touch each other. At this point, their tangents coincide. Therefore, differentiating equation \eqref{eqn:Stability_Uniform_Root_Finding_Equations_LR1} and equating the derivaties, we obtain
	\begin{equation}\label{eqn:Stability_Uniform_Root_Finding_Equations_LR1_Tangents}
		\frac{-\mathfrak{M}_n}{\rho^{\mathfrak{M}_n+1}}=(-1)^{l_n}l_n(R^n+1)[(R^n+1)\rho-R^n]^{l_n-1}.
	\end{equation}
	Dividing equation \eqref{eqn:Stability_Uniform_Root_Finding_Equations_LR1_Tangents} by equation \eqref{eqn:Stability_Uniform_Root_Finding_Equations_LR1}, we obtain
	\begin{equation}\label{eqn:Stability_Uniform_Root_Value_rho}
		\frac{\mathfrak{M}_n}{\rho}+\frac{l_n(R^n+1)}{(R^n+1)\rho-R^n}=0,
	\end{equation}
	which yields
	\begin{equation}\label{eqn:Stability_Uniform_Root_Value_rho_final}
		\rho=\frac{\mathfrak{M}_n}{\mathfrak{M}_n+l_n}\frac{R^n}{R^n+1}.
	\end{equation}
	
	Substituting this expression for $\rho$ into equation \eqref{eqn:Stability_Uniform_Root_Finding_Equations_LR1}, we obtain
	\begin{equation}\label{eqn:M_gamma_k}
		\left(1+\frac{1}{R^n}\right)^{\mathfrak{M}_n}=(l_n)^{l_n}\left(\frac{R^n}{\mathfrak{M}_n+l_n}\right)^{l_n}\left(\frac{\mathfrak{M}_n}{\mathfrak{M}_n+l_n}\right)^{\mathfrak{M}_n}.
	\end{equation}
	
	To characterise the onset of splitting, assume that $R^n$ larger than about $\gamma \mathfrak{M}_n$ would cause the stability region to break into two pieces. So, we put $R^n=\gamma \mathfrak{M}_n$ in \eqref{eqn:M_gamma_k}, we obtain 
	\begin{equation}
		\left(1+\frac{1}{\gamma \mathfrak{M}_n}\right)^{\mathfrak{M}_n}=(l_n)^{l_n}\gamma^{l_n}\left(\frac{\mathfrak{M}_n}{\mathfrak{M}_n+l_n}\right)^{l_n}\left(\frac{\mathfrak{M}_n}{\mathfrak{M}_n+l_n}\right)^{\mathfrak{M}_n}.
	\end{equation}
	
	Taking the limit as $\mathfrak{M}_n\to\infty$, we obtain the transidental equation for $\gamma$: 
	\begin{equation}\label{eqn:M_Critical_Transidental}
		(l_n)^{l_n}\gamma^{l_n}=\exp\left(l_n+\frac{1}{\gamma}\right).
	\end{equation}
	
	%\textcolor{red}{This can be written as}
	%\begin{equation}
	%	l\gamma=e^{1+\frac{1}{l\gamma}}
	%\end{equation}

	\begin{table}[h]
		\centering
		\caption{The values of the splitting parameter $\gamma$ are presented for various numbers of meso time steps $l_n$.}
		\begin{tabular}{p{1cm}p{1.5cm}|p{1cm}p{1.5cm}|p{1cm}p{1.5cm}}
			\hline
			$l_n$& $\gamma$&$l_n$&$\gamma$&$l_n$&$\gamma$\\
			\hline
			1&3.5911& 9&0.3990& 17&0.2112 \\
			2&1.7956& 10&0.3591& 18&0.1995 \\
			3&1.1970& 11&0.3265& 19&0.1890 \\
			4&0.8978& 12&0.2993& 20&0.1796 \\
			5&0.7182& 13&0.2762& 21&0.1710\\
			6&0.5985& 14&0.2565& 22&0.1632 \\
			7&0.5130& 15&0.2394& 23&0.1561 \\
			8&0.4489& 16&0.2245& 24&0.1496\\
			\hline
		\end{tabular}
		\label{table:Splitting_Parameters}
	\end{table}
	
	Table \ref{table:Splitting_Parameters} presents the splitting parameters $\gamma$ for various numbers of meso time steps $l_n$. The roots of the transcendental equation \eqref{eqn:M_Critical_Transidental} for various values of $l_n$ are calculated using the Newton-Raphson method. It can be observed that, as the number of meso steps within a macro time step increases, the value of the splitting parameter $\gamma$ decreases. This behaviour can be explained as follows: for a fixed total number of micro steps $\mathfrak{M}_n$ within macro time step $\Delta T^n$, increasing $l_n$ reduces the effective extrapolation step size. Consequently, the relation $R^n=\gamma\mathfrak{M}_n$ indicates that a smaller value of $\gamma$ is required to trigger the splitting of the stability region.
	
	%\textbf{\textcolor{red}{$\gamma$ values are in GP series. \{For any arbitrary values of $k_{n,j}\le k$ this series does not change, but\} As the number of meso time steps are increasing the splitting parameter is decreasing, which is clear from the equation (29). The splitting parameter proposed by Gear et al. is a fundamental splitting parameter.}}
	
	\begin{table}
		\centering
		\caption{For fixed values of $\mathfrak{M}_n$ and $l_n$, numerical approximations of $R^n$ are shown, where $R^n_s$ and $R^n_d$ denote the values of $R^n$ for which the stability region is single and double (i.e., split into two components), respectively. The parameter $\gamma_\text{numerical}$ is then defined as the numerical splitting parameter, given by $\frac{R^n_{\text{approx}}}{\mathfrak{M}_n}$. }
		\begin{tabular}{|cccc|} %{|p{1.2cm}|p{3cm}|p{3.3cm}|p{1.7cm}| }
			\hline
			$\mathfrak{M}_n$&$R^n_s$&$R^n_d$&$\gamma_\text{numerical}$\\
			\hline
			\multicolumn{4}{|c|}{$l_n$=1}\\
			\hline
			20&73.1&73.2&3.6575\\
			200&719.5&719.6&3.5977\\
			2000&7183.5&7183.6&3.5918\\
			20000&71823&71824&3.5912\\
			200000&718224&718225&3.5911\\
			2000000&7182201&7182202&3.5911\\
			\hline
			\multicolumn{4}{|c|}{$l_n$=2}\\
			\hline
			20&37.19&37.20&1.8598\\
			200&360.40&360.41&1.8020\\
			2000&3592.41&3592.42&1.7962\\
			20000&35912&35913&1.7956\\
			200000&359113&359114&1.7956\\
			\hline
			\multicolumn{4}{|c|}{$l_n$=3}\\
			\hline
			20&25.22&25.23&1.2612\\
			200&240.70&240.71&1.2035\\
			2000&2395.37&2395.38&1.1977\\
			20000&23942.10&23942.11&1.1971\\
			200000&239409.41&239409.42&1.1970\\
			2000000&2394082.19&2394082.20&1.1970\\
			\hline
		\end{tabular}
		\label{table:Stability_Numerical_Validation}
	\end{table}
	
	Table \ref{table:Splitting_Parameters} presents the analytically obtained splitting parameters. We now validate the agreement between the analytical predictions and the numerical results through direct numerical visualisation. For stability, the amplification factor in equation \eqref{eqn:Stability_Uniform_Amplification_Factor} must satisfy
	\begin{equation}\label{inequ:Stability_Numerical}
		\left|[(R^n+1)\rho-R^n]^{l_n}\rho^{\mathfrak{M}_n}\right|\le1,
	\end{equation}
	which defines the stability region over the macro time step $\Delta T^n$. For fixed values of $\mathfrak{M}_n$ and $l_n$, we plot the inequality \eqref{inequ:Stability_Numerical}. In doing so, we identify two nearby values of $R^n$ such that a slight increase from the smaller value causes the stability region to transition from a single connected component to two disconnected components. These numerical values are denoted by $R^n_s$ and $R^n_d$ in Table \ref{table:Stability_Numerical_Validation}, where the subscripts ``$s$" and ``$d$" represent ``single" and ``double", respectively. 
	
	For the numerical experiments, we consider $l_n=1$, 2 and 3, although other values may also be examined. The numerical splitting parameter is computed using the relation $\gamma=\frac{R^n_\text{approx}}{\mathfrak{M}_n}$, where $R^n_\text{approx}$ is approximated by $\frac{R^n_s+R^n_d}{2}$. As $\mathfrak{M}_n$ increases, the corresponding value of $R^n$ required to observe the splitting also increases. Consequently, the numerical values of $\gamma$ form a sequence that converges to the analytical splitting parameter reported in Table \ref{table:Stability_Numerical_Validation}. For example, when $l_n=2$, the sequence of numerical values $\gamma_\text{numerical}$ is 1.8598, 1.8020, 1.7962 and 1.7956 and rest of the terms are the repetition of the last term. Thus, the numerical splitting parameter converges to 1.7956, which agrees with the analytical value reported in Table \ref{table:Splitting_Parameters} for $l_n=2$. Therefore, we conclude that the numerically computed splitting parameter $\gamma_\text{numerical}$ is in excellent agreement with the analytically derived splitting parameter $\gamma$.

	%\textbf{\textcolor{red}{The index of the macro 'n' should be  removed.}}    
	
	\subsubsection{The values of $R^{n,m}$ are non-uniform in macro time step}\label{subsubsec:Non_Uni_R_nm}
	
	For non-uniform values of $R^{n,m}$ within a macro time step $\Delta T^n$, the amplification factor is given by 
	\begin{equation}\label{eqn:Stability_Gap_Vari_R}
		\sigma^{n}=\prod_{m=0}^{l_n-1}\left\{\left(R^{n,m}+1\right)\rho-R^{n,m}\right\}\rho^{\mathfrak{M}_n}.
	\end{equation}
	To investigate the splitting of the stability region, we examine the roots of $\sigma^n=\pm1$, where the sign depends on the number of meso time steps, as discussed in equation \eqref{eqn:Stability_Amplification_Boundary}. This leads to the following equations: 
	\begin{equation}\label{eqn:Stability_NonUniform_Root_Finding_Equations_LR1}
		\frac{1}{\rho^{\mathfrak{M}_n}}=(-1)^{l_n}\prod_{m=0}^{l_n-1}\left\{\left(R^{n,m}+1\right)\rho-R^{n,m}\right\},
	\end{equation}
	and
	\begin{equation}\label{eqn:Stability_NonUniform_Root_Finding_Equations_LR2}
		\frac{1}{\rho^{\mathfrak{M}_n}}=(-1)^{l_n+1}\prod_{m=0}^{l_n-1}\left\{\left(R^{n,m}+1\right)\rho-R^{n,m}\right\}.
	\end{equation}
	
	The difference between the amplification factors \eqref{eqn:Stability_Uniform_Amplification_Factor} (uniform case) and \eqref{eqn:Stability_Gap_Vari_R} (non-uniform case) lies in the first product term. For uniform $R^{n,m}$, the expression $[(R^n+1)\rho-R^n]^{l_n}$ intersects the $\rho$-axis at a single point $\rho=\frac{R^n}{R^n+1}$. In contrast, for non-uniform $R^{n,m}$, the product term $\prod_{m=0}^{l_n-1}\left\{\left(R^{n,m}+1\right)\rho-R^{n,m}\right\}$ intersects the $\rho$-axis at multiple points $\rho=\frac{R^{n,m}}{R^{n,m}+1}$, for $m=0$, 1,$\ldots$, $l_n-1$. We define,
	\begin{center}
		$\displaystyle
		\hat{\rho}_{\min}=\min_{m}\frac{R^{n,m}}{R^{n,m}+1}
		\quad \text{and} \quad
		\hat{\rho}_{\max}=\max_{m}\frac{R^{n,m}}{R^{n,m}+1}.
		$
	\end{center}
	Thus, for non-uniform $R^{n,m}$, the product term $\prod_{m=0}^{l_n-1}\left\{\left(R^{n,m}+1\right)\rho-R^{n,m}\right\}$ intersects the $\rho$-axis at multiple points within the interval $[\rho_{\min},\rho_{\max}]$. In the uniform case, the equations $\sigma^n=\pm(-1)^{l_n}$ admit at most four distinct real roots. However, in the non-uniform case, $\sigma^n=\pm(-1)^{l_n}$ equations may admit more than four distinct real roots. Outside the interval $[\rho_{\min},\rho_{\max}]$, denoted by $[\rho_{\min},\rho_{\max}]^c$, both the uniform and non-uniform cases exhibit similar behaviour, as discussed in Subsection \ref{subsubsec:Uni_R_nm}. In this region, the number of distinct real roots remains at most four. However, within the interval $[\rho_{\min},\rho_{\max}]$, equations \eqref{eqn:Stability_NonUniform_Root_Finding_Equations_LR1} and \eqref{eqn:Stability_NonUniform_Root_Finding_Equations_LR2} may admit additional real roots beyond these four. Consequently, for non-uniform $R^{n,m}$, the total number of distinct real roots can exceed four, leading to the possibility of more than two disconnected stability regions.
	
	Furthermore, in the region $[\rho_{\min},\rho_{\max}]^c$, bahaviours of the equations \eqref{eqn:Stability_NonUniform_Root_Finding_Equations_LR1} and \eqref{eqn:Stability_NonUniform_Root_Finding_Equations_LR2} remain analogous to the uniform case. Therefore, a critical value $R_\text{critical}^{n,m}$ may arise for equation \eqref{eqn:Stability_NonUniform_Root_Finding_Equations_LR1}. However, if such a critical value occurs within the interval $[\rho_{\min},\rho_{\max}]$, a more detailed investigation is required. At $R_\text{critical}^{n,m}$, the tangents to the left and right sides of \eqref{eqn:Stability_NonUniform_Root_Finding_Equations_LR1} are identical. Thus, we obtain 
	\begin{equation}\label{eqn:Stability_NonUniform_Root_Finding_Equations_LR1_Tangents}
		\frac{-\mathfrak{M}_n}{\rho^{\mathfrak{M}_n+1}}=(-1)^{l_n}\sum_{j=0}^{l_n-1}(R^{n,j}+1)\prod_{\substack{m=0 \\ m \ne j}}^{l_n-1}[(R^{n,m}+1)\rho-R^{n,m}].
	\end{equation}
	Dividing equation \eqref{eqn:Stability_NonUniform_Root_Finding_Equations_LR1_Tangents} by \eqref{eqn:Stability_NonUniform_Root_Finding_Equations_LR1}, we obtain
	\begin{equation}\label{eqn:Stability_General_Polynomial}
		\begin{aligned}
			\frac{\mathfrak{M}_n}{\rho}+\sum_{m=0}^{l_n-1}\frac{R^{n,m}+1}{(R^{n,m}+1)\rho-R^{n,m}}&=0,\\
			\frac{\mathfrak{M}_n}{\rho}+\sum_{m=0}^{l_n-1}\frac{1}{\rho-\hat{\rho}^{n,m}}&=0,
		\end{aligned}
	\end{equation}
	where, $\hat{\rho}^{n,m}=\frac{R^{n,m}}{R^{n,m}+1}$. This provides a general formulation for determining the values of $\rho$ in the case of variable micro, meso and macro time steps with time-dependent spectra, under the assumption that $\Theta=\lambda^{n,m,p}\delta T^{n,m,p}$ remains fixed.
	
	\begin{enumerate}
		\item For uniform $R^{n,m}=R^n$, equation \eqref{eqn:Stability_Uniform_Root_Value_rho} is a special case of the general equation \eqref{eqn:Stability_General_Polynomial}.
		\item \textbf{Only non-uniform macro-micro scales for time-dependent spectra:} Suppose that no meso scale is present, i.e., $l_n=1$, while the macro and micro time steps are non-uniform and the spectrum is time-dependent. This setting corresponds to a particular case of the general framework discussed above. For this configuration of the \texttt{GPI} scheme, the splitting parameter is $\gamma=3.5911$, as reported in Table \ref{table:Splitting_Parameters} and numerically verified in Table \ref{table:Stability_Numerical_Validation}.
		\item \textbf{Only uniform macro-micro scales for time-independent spectra:} The projective integration method proposed by Gear et al. \cite{2003_Gear_Projective} considers uniform macro and micro time steps with time-independent spectra, implying that $\Theta$ remains constant. In addition, no mesoscale is present, i.e., $l_n=1$ and the parameter $R^n$ is uniform. Therefore, this setting is a particular case of the scenario described above. In this case, the splitting parameter is reported as $\gamma=3.6$ in the article \cite{2003_Gear_Projective}, which is in close agreement with the more precise value $\gamma=3.5911$ obtained in this present work. 
	\end{enumerate}

	%That is we may write
	%\begin{center}
	%	\textit{Gear et al.}~\cite{2003_Gear_Projective} corresponds to a special case of $l_n = 1$ with uniform $R^{n,m}$ (see Subsection~\ref{subsubsec:Uni_R_nm}), which is a particular instance of the more general framework with uniform $R^{n,m}$, which is further extended to the case of non-uniform $R^{n,m}$ (see Subsection~\ref{subsubsec:Non_Uni_R_nm}).
	%\end{center}
	
	Since analysing the splitting of the stability region for the general equation \eqref{eqn:Stability_General_Polynomial} is highly challenging, in this work we restrict our attention to the case $l_n=2$. The analysis for larger numbers of meso steps needs future investigation.
	
	\vspace{0.2cm}

	$\bullet$ \textbf{Splitting of the stability region for the general equation \eqref{eqn:Stability_General_Polynomial} when $l_n=2$:}
	
	For $l_n=2$, Equation \eqref{eqn:Stability_General_Polynomial} reduces to
	\begin{equation}
		\frac{\mathfrak{M}_n}{\rho}+\sum_{m=0}^{1}\frac{R^{n,m}+1}{(R^{n,m}+1)\rho-R^{n,m}}=0
	\end{equation}
	which leads to a quadratic equation in $\rho$. This equation can be written as 
	\begin{equation}
		(\mathfrak{M}_n+2)(R^{n,0}+1)(R^{n,1}+1)\rho^2-(\mathfrak{M}_n+1)\left\{R^{n,0}(R^{n,1}+1)+R^{n,1}(R^{n,0}+1)\right\}\rho+\mathfrak{M}_nR^{n,0}R^{n,1}=0.
	\end{equation}
	
	The roots of this quadratic equation are given by 
	\begin{equation}\label{eqn:M_Critical_Roots_l2}
		%\resizebox{\textwidth}{!}{$
			\begin{aligned}
				\rho = \frac{(\mathfrak{M}_n+1)S \pm \sqrt{(\mathfrak{M}_n+1)^2 S^2 - 4\mathfrak{M}_n(\mathfrak{M}_n+2)\mathcal{B}}}{2(\mathfrak{M}_n+2)(R^{n,0} + 1)(R^{n,1} + 1)},
			\end{aligned}
			%$}
	\end{equation}
	where
	\[
	\begin{aligned}
		S &:= R^{n,0}(R^{n,1} + 1) + R^{n,1}(R^{n,0} + 1),\\
		\mathcal{A} &:= R^{n,0}- R^{n,1},\\
		\mathcal{B} &:=R^{n,0} R^{n,1} (R^{n,0} + 1)(R^{n,1} + 1).
	\end{aligned}
	\]
	
	The discriminant of this quadratic polynomial is
	\begin{equation}
		\begin{split}
			&(\mathfrak{M}_n+1)^2S^2-4\mathfrak{M}_n(\mathfrak{M}_n+2)\mathcal{B}\\
			&=(\mathfrak{M}_n+1)^2\{R^{n,0}(R^{n,1}+1)-R^{n,1}(R^{n,0}+1)\}^2+4R^{n,0}R^{n,1}(R^{n,0}+1)(R^{n,1}+2)\\
			&=(\mathfrak{M}_n+1)^2\mathcal{A}^2+4\mathcal{B}>0.
		\end{split}
	\end{equation}
	Hence, both roots are real and distinct.
	
	Finally, we introduce the scaling 
	\begin{equation}\label{eqn:Relation_M_k}
		R^{n,m}=\gamma c^{n,m}\mathfrak{M}_n, \text{where}\hspace{0.2cm} c^{n,m}>0,\hspace{0.2cm} m=0, 1.
	\end{equation} 
	%\textbf{\textcolor{red}{Discuss about the $\frac{R^{n,j}}{k_{n,j}}=\gamma$ and for uniform case $c_{n,j}=1$. Discuss about our restriction.}}
	
	We obtain
	\begin{equation}
		\begin{split}
			(R^{n,0}+1)\rho-R^{n,0}&=\frac{(\mathfrak{M}_n+1)R^{n,1}(R^{n,0}+1)-(\mathfrak{M}_n+3)R^{n,0}(R^{n,1}+1)\pm\sqrt{(\mathfrak{M}_n+1)^2\mathcal{A}^2+4\mathcal{B}}}{2(\mathfrak{M}_n+2)(R^{n,1}+1)},\\
			&=\gamma \mathfrak{M}_n\frac{(c^{n,1}-c^{n,0}-2\gamma c^{n,0}c^{n,1})\mathfrak{M}_n+(c^{n,1}-3c^{n,0})\pm\sqrt{(\mathfrak{M}_n+1)^2\mathbb{A}^2+4\mathbb{B}}}{2(\mathfrak{M}_n+2)(\gamma c^{n,1}\mathfrak{M}_n+1)},
		\end{split}
	\end{equation}
	where 
	\[
	\mathbb{A} := c^{n,0} - c^{n,1}\hspace{0.5cm} \&\hspace{0.5cm} \mathbb{B}:=c^{n,0}c^{n,1}(\gamma c^{n,0}\mathfrak{M}_n+1)(\gamma c^{n,1}\mathfrak{M}_n+1).
	\]
	
	Similarly,
	\begin{equation}
		\begin{split}
			(R^{n,1}+1)\rho-R^{n,1}&=\frac{(\mathfrak{M}_n+1)R^{n,0}(R^{n,1}+1)-(\mathfrak{M}_n+3)R^{n,1}(R^{n,0}+1)\pm\sqrt{(\mathfrak{M}_n+1)^2\mathcal{A}^2+4\mathcal{B}}}{2(\mathfrak{M}_n+2)(R^{n,0}+1)},\\
			&=\gamma \mathfrak{M}_n\frac{(c^{n,0}-c^{n,1}-2\gamma c^{n,0}c^{n,1})\mathfrak{M}_n+(c^{n,0}-3c^{n,1})\pm\sqrt{(\mathfrak{M}_n+1)^2\mathbb{A}^2+4\mathbb{B}}}{2(\mathfrak{M}_n+2)(\gamma c^{n,0}\mathfrak{M}_n+1)}.
		\end{split}
	\end{equation}

	Eliminating $\rho$ from equation \eqref{eqn:Stability_NonUniform_Root_Finding_Equations_LR1} and using the relation \eqref{eqn:Relation_M_k}, we obtain
	\begin{equation}\label{eqn:M_Critical_rho_Eli_Rela_M_k}
		\begin{split}
			&\{(R^{n,0}+1)\rho-R^{n,0}\}\{(R^{n,1}+1)\rho-R^{n,1}\}\rho^{\mathfrak{M}_n}=1\\
			i.e.,&\left[\gamma \mathfrak{M}_n\frac{(c^{n,1}-c^{n,0}-2\gamma c^{n,0}c^{n,1})\mathfrak{M}_n+(c^{n,1}-3c^{n,0})\pm\sqrt{(\mathfrak{M}_n+1)^2\mathbb{A}^2+4\mathbb{B}}}{2(\mathfrak{M}_n+2)(\gamma c^{n,1}\mathfrak{M}_n+1)}\right]\\
			\times &\left[\gamma \mathfrak{M}_n\frac{(c^{n,0}-c^{n,1}-2\gamma c^{n,0}c^{n,1})\mathfrak{M}_n+(c^{n,0}-3c^{n,1})\pm\sqrt{(\mathfrak{M}_n+1)^2\mathbb{A}^2+4\mathbb{B}}}{2(\mathfrak{M}_n+2)(\gamma c^{n,0}\mathfrak{M}_n+1)}\right]\\
			\times&\left[\gamma \mathfrak{M}_n\frac{(\mathfrak{M}_n+1)(2\gamma c^{n,0}c^{n,1}\mathfrak{M}_n+c^{n,0}+c^{n,1})\pm\sqrt{(\mathfrak{M}_n+1)^2\mathbb{A}^2+4\mathbb{B}}}{2(\mathfrak{M}_n+2)(\gamma c^{n,0}\mathfrak{M}_n+1)(\gamma c^{n,1}\mathfrak{M}_n+1)}\right]^{\mathfrak{M}_n}=1.
		\end{split}
	\end{equation}
	
	Substituting the above expressions, we obtain an explicit equation in terms of $\gamma$, $\mathfrak{M}_n$ and $c^{n,m}$. Taking the limit as $\mathfrak{M}_n\rightarrow\infty$, equation \eqref{eqn:M_Critical_rho_Eli_Rela_M_k} reduces to 
	\begin{equation}\label{eqn:M_Critical_Limit_Trans}
		\begin{split}
			&\left[\frac{c^{n,1}-c^{n,0}-2\gamma c^{n,0}c^{n,1}\pm\sqrt{(c^{n,0}-c^{n,1})^2+4\gamma^2(c^{n,0})^2(c^{n,1})^2}}{2c^{n,1}}\right]\\
			\times &\left[\frac{c^{n,0}-c^{n,1}-2\gamma c^{n,0}c^{n,1}\pm\sqrt{(c^{n,0}-c^{n,1})^2+4\gamma^2(c^{n,0})^2(c^{n,1})^2}}{2c^{n,0}}\right]\\
			\times&\exp\left[\frac{-(2\gamma c^{n,0}c^{n,1}+c^{n,0}+c^{n,1})\pm\sqrt{(c^{n,0}-c^{n,1})^2+4\gamma^2(c^{n,0})^2(c^{n,1})^2}}{2\gamma c^{n,0}c^{n,1}}\right]=1.
		\end{split}
	\end{equation}
	
	This leads to the transcendental equation
	\begin{equation}
		2\gamma^2c^{n,0}c^{n,1}\mp\gamma\sqrt{(c^{n,0}-c^{n,1})^2+4\gamma^2(c^{n,0})^2(c^{n,1})^2}=\exp\left[1+\frac{c^{n,0}+c^{n,1}\mp\sqrt{(c^{n,0}-c^{n,1})^2+4\gamma^2(c^{n,0})^2(c^{n,1})^2}}{2\gamma c^{n,0}c^{n,1}}\right]
	\end{equation}
	corresponding to the roots in \eqref{eqn:M_Critical_Roots_l2}. 
	
	For the root
	$
	\begin{aligned}
		\rho = \frac{(\mathfrak{M}_n+1)S + \sqrt{(\mathfrak{M}_n+1)^2 S^2 - 4\mathfrak{M}_n(\mathfrak{M}_n+2)\mathcal{B}}}{2(\mathfrak{M}_n+2)(R^{n,0} + 1)(R^{n,1} + 1)},
	\end{aligned}
	$
	we obtain the transidental equation 
	\begin{equation}\label{eqn:Not_root_l2}
		2\gamma^2c^{n,0}c^{n,1}-\gamma\sqrt{(c^{n,0}-c^{n,1})^2+4\gamma^2(c^{n,0})^2(c^{n,1})^2}=\exp\left[1+\frac{c^{n,0}+c^{n,1}-\sqrt{(c^{n,0}-c^{n,1})^2+4\gamma^2(c^{n,0})^2(c^{n,1})^2}}{2\gamma c^{n,0}c^{n,1}}\right].
	\end{equation}
	If we consider the uniform case $R^{n,m}$ (as discussed in Subsection \ref{subsubsec:Uni_R_nm}), then $c^{n,0}=c^{n,1}=1$. In this case, the left-hand side of equation \eqref{eqn:Not_root_l2} becomes zero, whereas the right-hand side evaluates to $\exp\left(\frac{1}{\gamma}\right)$, leading to a contradiction. Therefore, this root is not admissible and must be discarded. 
	
	%(\textbf{\textcolor{red}{This below comparison is not connected with the proposed scheme. Instead of $M=\gamma k$, take $lM=\gamma k$ in subsection 6.1.1 and then compare.}}) 
	
	However, for the root 
	$
	\begin{aligned}
		\rho = \frac{(\mathfrak{M}_n+1)S - \sqrt{(\mathfrak{M}_n+1)^2 S^2 - 4\mathfrak{M}_n(\mathfrak{M}_n+2)\mathcal{B}}}{2(\mathfrak{M}_n+2)(R^{n,0} + 1)(R^{n,1} + 1)},
	\end{aligned}
	$
	the resulting transidental equation becomes 
	\begin{equation}\label{eqn:Yes_root_l2}
		2\gamma^2c^{n,0}c^{n,1}+\gamma\sqrt{(c^{n,0}-c^{n,1})^2+4\gamma^2(c^{n,0})^2(c^{n,1})^2}=\exp\left[1+\frac{c^{n,0}+c^{n,1}+\sqrt{(c^{n,0}-c^{n,1})^2+4\gamma^2(c^{n,0})^2(c^{n,1})^2}}{2\gamma c^{n,0}c^{n,1}}\right].
	\end{equation}
	
	$\bullet$ For uniform $R^{n,m}$, as discussed in Subsection \ref{subsubsec:Uni_R_nm}, we set $c^{n,0}=c^{n,1}=1$ in equation \eqref{eqn:Yes_root_l2}. This yields
	\begin{equation}
		2^2(\gamma)^2=\exp\left(2+\frac{1}{\gamma}\right).
	\end{equation}
	which coincides with equation \eqref{eqn:M_Critical_Transidental} for $l_n=2$.
	
	\begin{table}
		\centering
		\caption{For fixed values of $\mathfrak{M}_n$ and $l_n=2$, numerical approximations of $R^{n,m}$ are computed for $c^{n,0}=0.4$ and $c^{n,1}=0.6$. The quantities $R^{n,m}_s$ and $R^{n,m}_d$ denote the values of $R^{n,m}$ for which the stability region is single and double (i.e., split into two components), respectively. The parameter $\gamma_\text{numerical}$ is then defined as the numerical splitting parameter, given by $\frac{R^{n,m}_{\text{approx}}}{\mathfrak{M}_n}$.}
		\begin{tabular}{cccc}
			\hline
			$\mathfrak{M}_n$& $R_s^n$&$R_d^n$&$\gamma_\text{numerical}$\\
			\hline
			20&76.3380&76.3383&3.8169\\
			200&739.8685&739.8687&3.6993\\
			2000&7375&7375.0015&3.6875\\
			20000&73726.4105&73726.4107&3.6863\\
			200000&737240.5565&737240.5567&3.6862\\
			2000000&7372382.0205&7372382.0208&3.6862\\
			\hline
		\end{tabular}
		\label{table:Stability_Numerical_Validation_Non_Uni_l2_Example1}
	\end{table}
	
	\begin{table}
		\centering
		\caption{For fixed values of $\mathfrak{M}_n$ and $l_n=2$, numerical approximations of $R^{n,m}$ are computed for $c^{n,0}=0.2$ and $c^{n,1}=0.8$. The quantities $R^{n,m}_s$ and $R^{n,m}_d$ denote the values of $R^{n,m}$ for which the stability region is single and double (i.e., split into two components), respectively. The parameter $\gamma_\text{numerical}$ is then defined as the numerical splitting parameter, given by $\frac{R^{n,m}_{\text{approx}}}{\mathfrak{M}_n}$.}
		\begin{tabular}{cccc}
			\hline
			$\mathfrak{M}_n$& $R_s^n$&$R_d^n$&$\gamma_\text{numerical}$\\
			\hline
			20&19.7919&19.7920&4.9480\\
			200&191.6891&191.6892&4.7922\\
			2000&1918.8133&1918.8134&4.7970\\
			20000&19183.2129&19183.2130&4.7958\\
			200000&191828.1197&191828.1198&4.7957\\
			2000000&1918275.6872&1918275.6873&4.7957\\
			\hline
		\end{tabular}
		\label{table:Stability_Numerical_Validation_Non_Uni_l2_Example2}
	\end{table}
	
	To verify the analytically obtained splitting parameter, we plot the stability region defined by 
	\begin{equation}
		\left|\prod_{m=0}^{l_n-1}\left\{\left(R^{n,m}+1\right)\rho-R^{n,m}\right\}\rho^{\mathfrak{M}_n}\right|\le1.
	\end{equation}   
	In particular, we consider two numerical experiments. In the first case, we set $R^{n,1}=\frac{3}{2}R^{n,0}$, which corresponds to $c^{n,0}=0.4$ and $c^{n,1}=0.6$. In the second case, we set $R^{n,1}=4R^{n,0}$, which corresponds to $c^{n,0}=0.2$ and $c^{n,1}=0.8$. The transcendental equation \eqref{eqn:Yes_root_l2} is then solved using the Newton--Raphson method, yielding $\gamma=3.6862$ and $\gamma=4.7957$ for the first and second cases, respectively.
	
	Following the same procedure as in Table \ref{table:Stability_Numerical_Validation}, we compute the stability regions and identify the values of $R^{n,m}$ at which the stability region transitions from a single connected component to two disconnected components. From these values, we estimate the numerical splitting parameter $\gamma_\text{numerical}$, which forms a sequence as $\mathfrak{M}_n$ increases. The numerical values $\gamma_\text{numerical}$ converge to the corresponding analytical splitting parameters in both examples, as shown in Tables \ref{table:Stability_Numerical_Validation_Non_Uni_l2_Example1} and \ref{table:Stability_Numerical_Validation_Non_Uni_l2_Example2}. In the limit $\mathfrak{M}_n\rightarrow\infty$, the numerical splitting parameter agrees with the analytical value $\gamma$.

	%\textcolor{red}{For $l_n$ meso scales, it can have maximum $2(l_n+1)$ roots, effectively it can have maximum $l_n+1$ splitted stability regions.}
	
	\begin{figure}
		\begin{subfigure}{.32\textwidth}
			\centering
			% include third image
			\includegraphics[width=1\linewidth]{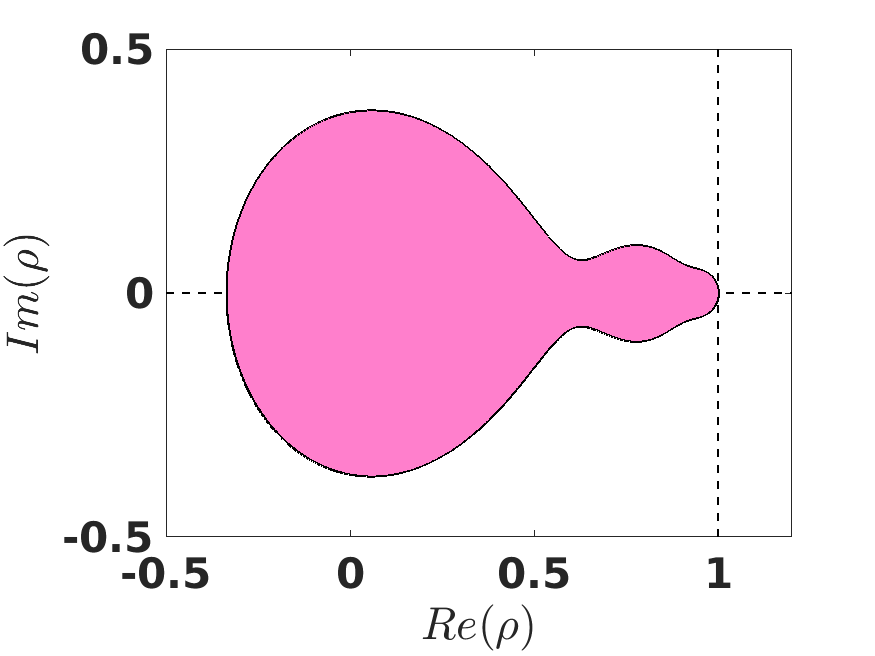}  
			\caption{$R^{n,0}=30$}
			\label{fig:One_Regions_k5_R_30_4}
		\end{subfigure}
		\begin{subfigure}{.32\textwidth}
			\centering
			% include fourth image
			\includegraphics[width=1\linewidth]{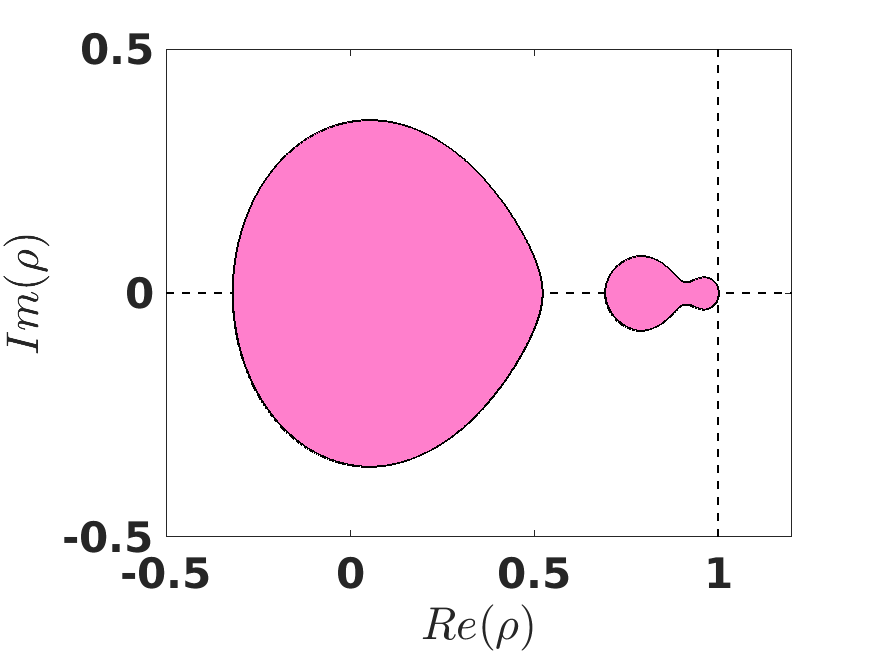}  
			\caption{$R^{n,0}=40$}
			\label{fig:Two_Regions_k5_R_40_4}
		\end{subfigure}
		\begin{subfigure}{0.32\textwidth}
			\centering
			% include fourth image
			\includegraphics[width=1\linewidth]{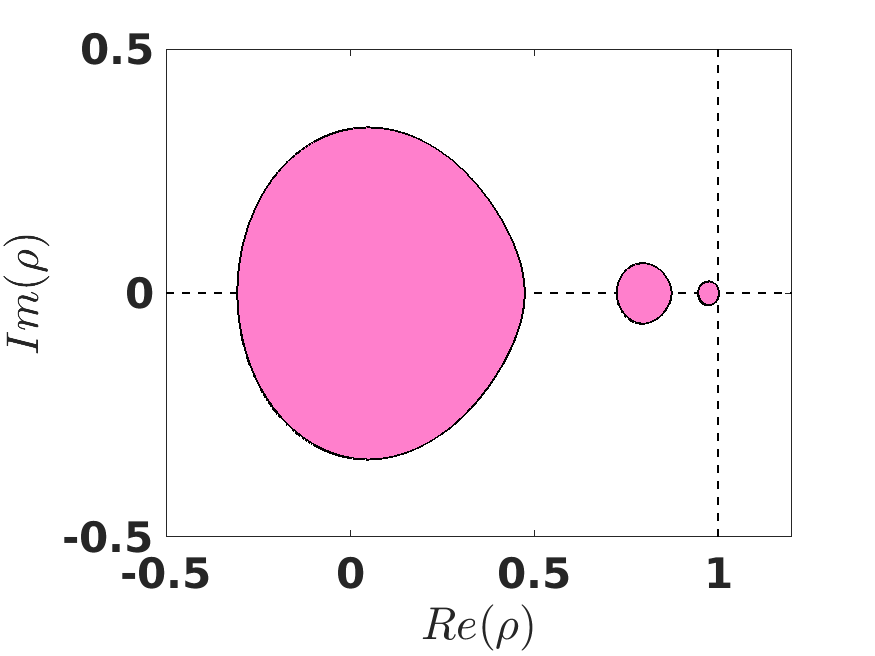}  
			\caption{$R^{n,0}=50$}
			\label{fig:Three_Regions_k5_R_50_4}
		\end{subfigure}
		\caption{Splitting of the stability region in \texttt{GPI} scheme due to variations in the parameter $R^{n,m}$. Here, $l_n=2$, $\mathfrak{M}_n=5$ and $R^{n,1}=4$ are set with (a) $R^{n,0}=30$, (b) $R^{n,0}=40$ and (c) $R^{n,0}=50$.}
		\label{fig:One_To_Three_Regions_k5}
	\end{figure}
	
	For two non-uniform meso steps ($l_n=2$), as discussed earlier, equations \eqref{eqn:Stability_NonUniform_Root_Finding_Equations_LR1} and \eqref{eqn:Stability_NonUniform_Root_Finding_Equations_LR2} may admit more than four distinct real roots. The term $\prod_{m=0}^{1}\left\{\left(R^{n,m}+1\right)\rho-R^{n,m}\right\}$ intersects the real $\rho$-axis at two points, namely $\rho=\frac{R^{n,m}}{R^{n,m}+1}$ for $m=0$ and 1. Therefore, in addition to at most four distinct real roots in $[\rho_{\min},\rho_{\max}]^c$, the equations \eqref{eqn:Stability_NonUniform_Root_Finding_Equations_LR1} and \eqref{eqn:Stability_NonUniform_Root_Finding_Equations_LR2} may admit up to two additional distinct real roots within the interval $[\rho_{\min},\rho_{\max}]$. Consequently, for two meso steps $l_n=2$, there can be at most six distinct real roots of these equations. This implies that the stability region may split into at most three disconnected components. 
	
	To illustrate this behaviour, we consider the example with $\mathfrak{M}_n=5$, $R^{n,1}=4$ and $R^{n,0}=30$, 40, 50. The corresponding stability regions are shown in Figures \ref{fig:One_Regions_k5_R_30_4}, \ref{fig:Two_Regions_k5_R_40_4} and \ref{fig:Three_Regions_k5_R_50_4}, respectively. For $R^{n,0}=30$, the stability region is single and connected. As $R^{n,0}$ increases to 40, the stability region splits into two disconnected components. With a further increase to $R^{n,0}=50$, the stability region splits further into three disconnected components.
	
	\section{An adaptive strategy for selecting the micro-burst length}\label{sec:Relaxation_Time}
	
	All variants of projective integration schemes require a balance between the stability of the fast modes and the accurate resolution of the slow dynamics. During the microsimulation stage, the fast modes decay, whereas during the projective extrapolation stage, they may regrow. The interplay between these two mechanisms determines the overall stability and efficiency of the method.
	
	From the test equation \eqref{eqn:Test_eqn}, we consider
	\begin{equation}\label{eqn:Test_Eqn_Max_Eigen}
		\frac{du}{dt}=-\lambda_*(t)u,
		\qquad 
		\lambda_*(t):=\max_{p}|\lambda_p(t)|\gg0,
	\end{equation}
	where $\lambda_p(t)<0$ denote the eigenvalues of the Jacobian corresponding to the system \eqref{eqn:GPI_General_Eqn}. The exact solution of \eqref{eqn:Test_Eqn_Max_Eigen} is given by
	\begin{equation}\label{eqn:GPI_Test_Eqn_Exact}
		u(t+\tau)
		=
		\exp\left(
		-\int_{t}^{t+\tau}\lambda_*(s)\,ds
		\right)
		u(t).
	\end{equation}
	
	Therefore, over a micro-burst interval of length $\tau$, the fast modes decay by the factor
	\[
	\exp\left(
	-\int_{t}^{t+\tau}\lambda_*(s)\,ds
	\right).
	\]
	This quantity measures the damping of the fast components during the inner microsimulation and characterises the relaxation of the solution toward the slow manifold.
	
	In the \texttt{GPI} scheme, suppose that the microsimulation is performed over a duration $\tau^{n,m}$ and is followed by a projective extrapolation over the remaining interval 
	\[
	\bar{\partial T}^{n,m}:=\partial T^{n,m}-\tau^{n,m},
	\]
	using the forward Euler scheme within the meso time step $\partial T^{n,m}$. Then,
	\begin{equation}
		u^{n,m+1}
		\approx
		\left(
		1-\lambda_*(\bar{t})\bar{\partial T}^{n,m}
		\right)
		u(\bar{t}),
	\end{equation}
	where 
	\[
	\bar{t}=T^{n,m}+\tau^{n,m}.
	\]
	
	Substituting the exact solution \eqref{eqn:GPI_Test_Eqn_Exact} into the above expression yields
	\begin{equation}
		u^{n,m+1}
		\approx
		\left(
		1-\lambda_*(\bar{t})\bar{\partial T}^{n,m}
		\right)
		\exp\left(
		-\int_{T^{n,m}}^{T^{n,m}+\tau^{n,m}}
		\lambda_*(s)\,ds
		\right)
		u^{n,m}.
	\end{equation}
	
	The first factor represents the numerical amplification introduced by the projective extrapolation, whereas the exponential term corresponds to the physical damping of the fast modes during the microsimulation stage. Therefore, the overall amplification factor over the meso time step $\partial T^{n,m}$ is given by
	\begin{equation}
		\tilde{\sigma}^{n,m}
		=
		\left(
		1-\lambda_*(\bar{t})\bar{\partial T}^{n,m}
		\right)
		\exp\left(
		-\int_{T^{n,m}}^{T^{n,m}+\tau^{n,m}}
		\lambda_*(s)\,ds
		\right).
	\end{equation}
	
	Since $\lambda_*(t)\gg0$ is assumed to be large, we approximate
	\[
	1-\lambda_*(\bar{t})\bar{\partial T}^{n,m}
	\approx
	-\lambda_*(\bar{t})\bar{\partial T}^{n,m}.
	\]
	Consequently, the stability condition becomes
	\begin{equation}\label{eqn:Relax_stability_Inequality}
		\lambda_*(\bar{t})\bar{\partial T}^{n,m}
		\exp\left(
		-\int_{T^{n,m}}^{T^{n,m}+\tau^{n,m}}
		\lambda_*(s)\,ds
		\right)
		\le1.
	\end{equation}
	
	In the \texttt{GPI} framework,
	\[
	\tau^{n,m}
	=
	\sum_{i=0}^{M_{n,m}}
	\delta T^{n,m,i},
	\]
	which denotes the total microsimulation duration over which the fast variables relax toward the slow manifold.

	\paragraph{Case I: Time-independent $\lambda_*$}
	\begin{enumerate}
		\item \textbf{Non-uniform macro, meso and micro scales:} 
		If all three scales---macro, meso and micro---are present and $\lambda_*$ is time-independent, then the stability condition \eqref{eqn:Relax_stability_Inequality} reduces to
		\begin{equation}
			\lambda_*\bar{\partial T}^{n,m}
			\exp(-\lambda_*\tau^{n,m})
			\le1,
		\end{equation}
		which implies
		\begin{equation}
			\tau^{n,m}
			\ge 
			\frac{1}{\lambda_*}
			\log\left(
			\lambda_*\bar{\partial T}^{n,m}
			\right).
		\end{equation}
		
		Suppose that the meso time step $\partial T^{n,m}$ is prescribed and the corresponding micro-burst length is to be determined. Since the projective extrapolation interval $\bar{\partial T}^{n,m}$ depends on the burst length, it is generally not known a priori. However, using the relation
		\[
		\partial T^{n,m}\ge\bar{\partial T}^{n,m},
		\]
		we obtain
		\begin{equation}\label{eqn:GPI_Dt_greater_Dt_bar}
			\frac{1}{\lambda_*}
			\log\left(
			\lambda_*\partial T^{n,m}
			\right)
			\ge
			\frac{1}{\lambda_*}
			\log\left(
			\lambda_*\bar{\partial T}^{n,m}
			\right).
		\end{equation}
		
		Therefore, a practical lower bound for the relaxation time corresponding to each meso time step $\partial T^{n,m}$ is given by
		\begin{equation}
			\tau^{n,m}
			\ge
			\frac{1}{\lambda_*}
			\log\left(
			\lambda_*\partial T^{n,m}
			\right).
		\end{equation}
		
		\item \textbf{Non-uniform macro and micro scales only:} 
		In the \texttt{GPI} framework, if only macro and micro scales are present, that is, $l_n=1$ for all $n$, then
		\[
		\partial T^{n,0}=\Delta T^n.
		\]
		
		For time-independent $\lambda_*$, the micro-burst length in each stage of the \texttt{GPI} scheme should satisfy
		\begin{equation}\label{eqn:GPRK_Burst_Time_Independent}
			\tau^{n}
			\ge
			\frac{1}{\lambda_*}
			\log\left(
			\lambda_*\Delta T^n
			\right),
			\qquad 
			n=0,1,\ldots,N_t-1.
		\end{equation}
		
		This condition is applicable to general projective integration schemes with non-uniform macro and micro scales under time-independent $\lambda_*$ throughout the simulation interval.
		
		\item \textbf{Uniform macro and non-uniform micro scales only:} 
		If the macro time steps are uniform, as commonly assumed in existing projective integration schemes such as \texttt{PI} \cite{2003_Gear_Projective,2012_Lafitte_Asymptotic}, \texttt{PRK} \cite{2016_Lafitte_High-Order,2017_Lafitte_High-order}, and \texttt{PIRK} \cite{2020_Roberts_toolbox,2021_maclean_toolbox}, then the relaxation time condition simplifies to
		\begin{equation}\label{eqn:Formula_Burst_Uniform_Macro_Time_Ind_Spectra}
			\tau
			\ge
			\frac{1}{\lambda_*}
			\log\left(
			\lambda_*\Delta T
			\right).
		\end{equation}
		
		This expression is consistent with the relaxation-time condition reported by Maclean et al.~\cite{2021_maclean_toolbox}. Here, the index `$n$' is omitted since the macro time step is uniform throughout the simulation interval.
		
	\end{enumerate}

	\paragraph{Case II: Time-dependent $\lambda_*$}
	\begin{enumerate}
		\item \textbf{General time-dependent $\lambda_*(t)$:}
		
		Assume that $\lambda_*(t)$ attains its minimum value at some time 
		\[
		t_{\min}\in[T^{n,m},T^{n,m}+\tau^{n,m}],
		\]
		such that
		\begin{equation}
			\lambda_*^{\min}
			=
			\lambda_*(t_{\min})
			:=
			\min_{t\in[T^{n,m},T^{n,m}+\tau^{n,m}]}
			\lambda_*(t).
		\end{equation}
		
		Then,
		\[
		\exp\left(
		-\int_{T^{n,m}}^{T^{n,m}+\tau^{n,m}}
		\lambda_*(s)\,ds
		\right)
		\le
		\exp\left(
		-\lambda_*^{\min}\tau^{n,m}
		\right).
		\]
		
		Therefore, from the amplification factor, we obtain
		\begin{equation}
			\begin{aligned}
				|\tilde{\sigma}^{n,m}|
				&\le
				\left|
				\left(
				1-\lambda_*(\bar{t})\bar{\partial T}^{n,m}
				\right)
				\exp\left(
				-\lambda_*^{\min}\tau^{n,m}
				\right)
				\right| \\
				&\approx
				\lambda_*(\bar{t})
				\bar{\partial T}^{n,m}
				\exp\left(
				-\lambda_*^{\min}\tau^{n,m}
				\right),
			\end{aligned}
		\end{equation}
		where the approximation holds for sufficiently large values of $\lambda_*(t)$.
		
		To ensure the stability condition
		\[
		|\tilde{\sigma}^{n,m}|\le1,
		\]
		it is sufficient to impose
		\[
		\lambda_*(\bar{t})
		\bar{\partial T}^{n,m}
		\exp\left(
		-\lambda_*^{\min}\tau^{n,m}
		\right)
		\le1.
		\]
		
		Hence,
		\begin{equation}
			\bar{\partial T}^{n,m}
			\exp\left(
			-\lambda_*^{\min}\tau^{n,m}
			\right)
			\le
			\frac{1}{\lambda_*(\bar{t})}
			\le
			\frac{1}{\lambda_*^{\min}},
		\end{equation}
		which implies
		\begin{equation}
			\tau^{n,m}
			\ge
			\frac{1}{\lambda_*^{\min}}
			\log\left(
			\bar{\partial T}^{n,m}\lambda_*^{\min}
			\right).
		\end{equation}
		
		Using the same argument as in \eqref{eqn:GPI_Dt_greater_Dt_bar}, we obtain the practical relaxation-time bound
		\begin{equation}\label{eqn:Relaxation_Bound_Stiffness_General}
			\tau^{n,m}
			\ge
			\frac{1}{\lambda_*^{\min}}
			\log\left(
			\partial T^{n,m}\lambda_*^{\min}
			\right).
		\end{equation}
		
		\item \textbf{Strictly increasing $\lambda_*(t)$ over the meso time step $\partial T^{n,m}$:}
		
		In this case, $\lambda_*(t)$ attains its minimum value at the beginning of the meso interval, namely
		\[
		\lambda_*^{\min}
		=
		\lambda_*(T^{n,m}).
		\]
		
		Therefore, from \eqref{eqn:Relaxation_Bound_Stiffness_General}, the relaxation-time condition becomes
		\begin{equation}\label{eqn:Relaxation_Bound_Increasing_Stiffness}
			\tau^{n,m}
			\ge
			\frac{1}{\lambda_*(T^{n,m})}
			\log\left(
			\partial T^{n,m}\lambda_*(T^{n,m})
			\right).
		\end{equation}
		
		If the meso scale is absent, that is, $l_n=1$, then the micro-burst length should satisfy
		\begin{equation}\label{eqn:Relaxation_Bound_Macro_Increasing_Stiffness}
			\tau^{n}
			\ge
			\frac{1}{\lambda_*(T^{n})}
			\log\left(
			\Delta T^{n}\lambda_*(T^{n})
			\right).
		\end{equation}
		
		%\item \textbf{$\lambda_*(t)$ is strictly decreasing in the meso time step $\partial T^{n,m}$:} For this case, $\lambda_*(t)$ attains its minimum value $\lambda_*^{\min}=\lambda_*(T^{n,m+1})$ at $t_{\min}=T^{n,m+1}$ in the meso time step $\partial T^{n,m}$. Therefore from \eqref{eqn:Relaxation_Bound_Stiffness_General}, the relaxation time limit becomes
		%\begin{equation}\label{eqn:Relaxation_Bound_Decreasing_Stiffness}
		%	\bar{\delta}^{n,m}\ge \frac{1}{\lambda_*(T^{n,m+1})}\log\left(\bar{\partial T}^{n,m}\lambda_*(T^{n,m+1})\right).
		%\end{equation} 
		
	\end{enumerate}
	
	%Since the burst length is not known a priori and $\log\left(\bar{\partial T}^{n,m}\lambda_*^{\min}\right) < \log\left(\partial T^{n,m}\lambda_*^{\min}\right)$, it is therefore customary to use $\partial T^{n,m}$ in place of $\bar{\partial T}^{n,m}$. A similar approach is also adopted by Maclean et al. \cite{2021_maclean_toolbox}.

	%%%%%%%%%%%%%%%%%%%%%%%%%%%%%%%%%%%%%%%%%%%%%%%%%%%%%%%%%%%%%
	
	\section{A strategy for selecting non-uniform micro time steps}\label{sec:Choice_Micro_Steps}
	
	To ensure stability of the forward Euler microsolver, the micro time steps must satisfy the condition
	\[
	0<\delta T^{n,m,i}\le\frac{2}{\lambda_*(T^{n,m,i})},
	\]
	where
	\[
	\lambda_*(t):=\max_{p}|\lambda_p(t)|
	\]
	denotes the maximum magnitude of the eigenvalues of the local spectrum at time $t$.
	
	Motivated by this stability requirement, we select non-uniform micro time steps according to the local spectral properties of the system. In particular, within a meso time step $\partial T^{n,m}$, the micro time steps are chosen as
	\begin{equation}\label{eqn:Variable_Inner_Steps_Formula}
		\begin{aligned}
			&h^{n,m}_0
			=
			\frac{1}{\lambda_*(T^{n,m})},\\
			&\delta T^{n,m,j}
			=
			\frac{1}{\lambda_*\!\left(
				T^{n,m,j}+h^{n,m}_0
				\right)},\\
			&T^{n,m,j+1}
			=
			T^{n,m,j}
			+
			\delta T^{n,m,j},
		\end{aligned}
	\end{equation}
	for
	\[
	j=0,1,\ldots,M_{n,m},
	\qquad
	m=0, 1,\ldots,l_n-1,
	\qquad
	n=0,1,\ldots,N_t-1.
	\]
	
	This provides one possible adaptive strategy for selecting non-uniform micro time steps based on the evolving spectral properties of the system, as employed throughout this work. Depending on the nature of the problem and the desired numerical properties, alternative adaptive strategies may also be constructed within the proposed framework.

	%%%%%%%%%%%%%%%%%%%%%%%%%%%%%%%%%%%%%%%%%%%%%%%%%%%%%%%%%%%%%%%%
	
	\section{Results and discussion}\label{sec:Results_Discussion}
	The proposed generalised projective integration (\texttt{GPI}) scheme is validated through three representative test cases. The first case considers a nonlinear system of stiff ordinary differential equations (\texttt{ODEs}), in which the slow eigenvalue remains close to zero along the negative real axis, while the magnitude of the fast eigenvalue increases over time within the negative real plane (see Subsection \ref{subsec:Non_Linear_Stiff_variable_Stiffness}). In the second case, a linear diffusion problem is examined, where the entire spectrum evolves over time, as discussed in Subsection \ref{subsec:Pure_Diff}. Finally, a highly oscillatory Airy equation is studied in Subsection \ref{subsec:Airy}. All the numerical computations presented in this work were carried out using \texttt{MATLAB R2024a} on a desktop system equipped with an \texttt{Intel Core i5-7500} processor (3.40\,\texttt{GHz}) and 8\,\texttt{GB} RAM, running \texttt{Ubuntu 22.04 LTS}. 
	
	%\textbf{\textcolor{red}{Problem 1: Varying spectra as well as scale separation. Problem 2: Varying spectra but not scale separation.}}
	
	\subsection{System of nonlinear stiff \texttt{ODEs} with time-dependent spectra and scale separation}\label{subsec:Non_Linear_Stiff_variable_Stiffness}
	We consider a nonlinear system of stiff \texttt{ODEs} with time-dependent spectra and scale separation, given by
	\begin{equation}\label{eqn:NonLinear_Stiff_ODE_Increase_Stiffness}
		\begin{aligned}
			&\frac{du_1}{dt}=-u_1u_2-\alpha u_1^2, & u_1(t_0)=1,\\
			&\frac{du_2}{dt}=f(t)(-u_2+\sin^2(u_1)), & u_2(t_0)=0,\\
		\end{aligned}
	\end{equation}
	where $\alpha=0.2$ and $f(t)\gg1$ for all $t\in[t_0,T]$. The system exhibits a clear separation of time scales. One eigenvalue behaves as $\lambda_\text{fast}\approx -f(t)$, which attains large negative values as $f(t)$ increases, corresponding to a rapidly decaying (fast) mode. The other eigenvalue, $\lambda_\text{slow}=-(u_2+2\alpha u_1+u_1\sin(2u_1))$ remains of order
	$\mathcal{O}(1)$ and governs the slow evolution of the system. Since $\alpha$, $u_1(t)$ and $u_2(t)$ remain within the interval $[0,1]$ for $t\in[t_0,T]$, the system is stable. In particular, $\lambda_\text{slow}$ stays close to zero, whereas $\lambda_\text{fast}$ rapidly shifts deeper into the negative real axis as time progresses. This behaviour demonstrates a time-dependent separation of scales. 
	
	The corresponding reduced slow system associated with \eqref{eqn:NonLinear_Stiff_ODE_Increase_Stiffness} is given by
	\begin{equation}\label{eqn:Nonlinear_Effective_ODE_Increase_Stiffness}
		\begin{aligned}
			&\frac{dU}{dt}=-U\sin^2(U)-\alpha U^2, & U(t_0)=1.
		\end{aligned}
	\end{equation}
	
	From a physical perspective, the stiffness parameter $f(t)$ may vary with time. In the following subsection, we consider the case where $f(t)$ is an increasing function, implying that the stiffness of the system intensifies as time evolves.
	
	\subsubsection{Increasing stiffness: choice of $f(t)$}\label{subsubsec:Non_Linear_Increasing_Stiffness}
	
	Let $f(t)=10^t$, which is a strictly increasing function of time and we set the initial time as $t_0=2$. As time progresses, the fast dynamics become increasingly stiff due to the growth of $f(t)$. The problem is solved using the proposed \texttt{GPI} scheme. Its performance is compared with several existing equation-free multiscale methods, namely \texttt{PI}, \texttt{PRK2}, \texttt{PRK4}, \texttt{PIRK2}, \texttt{PIRK4}, \texttt{PIG2} and \texttt{PIG4}. These methods are primarily designed for systems with time-independent spectra or for problems where the separation scales do not vary significantly over time and they are known to perform efficiently in such settings. However, when the spectrum evolves in time--leading to a dynamically changing separation of scales--it becomes particularly important to assess how these existing schemes perform in comparison with the proposed \texttt{GPI} scheme.
	
	\paragraph{Projective integration (\texttt{PI}) scheme:}
	
	The projective integration (\texttt{PI}) scheme \cite{2003_Gear_Projective,2012_Lafitte_Asymptotic} is developed using uniform macro and micro time steps, along with a fixed burst length. To solve the problem \eqref{eqn:NonLinear_Stiff_ODE_Increase_Stiffness} using the \texttt{PI} scheme, the micro time step must satisfy the stability condition of the forward Euler method, such as
	\begin{center}
		$\delta T\le\frac{2}{\sup_{t \in [t_0, T]} \max_p |\lambda_p(t)|}$,
	\end{center}
	where $\lambda_p(t)$ denotes the time-dependent eigenvalues of the system. For the present problem, the dominant eigenvalue behaves approximately as $-f(t)=-10^t$. Accordingly, a conservative choice of the micro time step is taken as $\delta T=1\mathrm{e}{-T}$. The initial time is fixed at $t_0=2$, while the final time $T$ is varied starting from $3$ in order to examine performance. %Effectively, the micro time step size decreases, effectively the number of micro time steps increases to satisfy the relaxation time limit. 
	
	A uniform macro time step $\Delta T=5\mathrm{e}{-2}$ is used. For $T=3$, at least $M_n = 40$ micro steps of size $\delta T = 1\mathrm{e}{-3}$ are required within the first macro step to ensure relaxation of the fast variable towards the slow manifold. If fewer micro steps are used, the fast dynamics exhibit oscillatory behaviour during the initial macro step. 
	
	A key limitation of the \texttt{PI} scheme in this setting is that the micro time step depends explicitly on the final time $T$. As $T$ increases by one unit, $\delta T$ decreases by a factor of $10$. Consequently, to maintain the fixed relaxation time, the number of micro steps must increase by a factor of $10$ within each macro step. Thus, the number of micro steps per macro step becomes $40 \times 10^{T-3}$, with uniform step size $1\mathrm{e}{-T}$. This implies that the \texttt{PI} scheme spends a significant portion, approximately $80\%$ of the entire time interval $[t_0,T]$. The relaxation time for the micro simulation in the \texttt{PI} scheme is given by $\frac{1}{\lambda_*}\log\left(\Delta T\,\lambda_*\right)$ for constant $\lambda_*$, as mentioned in equation \eqref{eqn:Formula_Burst_Uniform_Macro_Time_Ind_Spectra}. However, in the present problem, $\lambda_*$ increases with time. As a result, the effective burst length decreases once $\lambda_*>\frac{\exp(1)}{\Delta T}\approx54.35$. This leads to a reduction in the required relaxation time as time progresses. Nevertheless, since the \texttt{PI} scheme enforces a uniform burst length, it performs excessive micro simulations, resulting in computational inefficiency.
	
\begin{figure}[h!]
		\begin{subfigure}{.5\textwidth}
			\centering
			% include third image
			\includegraphics[width=1\linewidth]{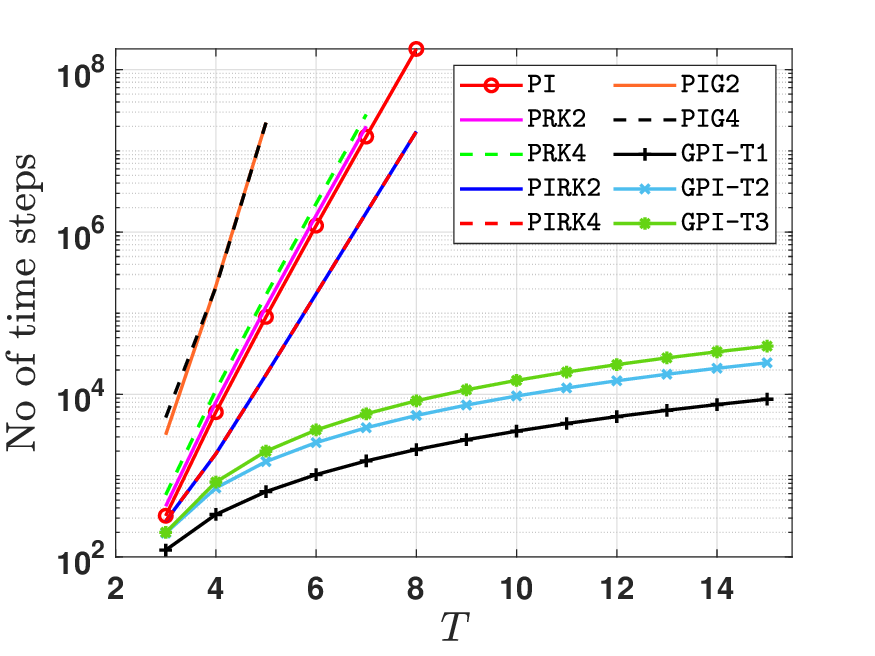}  
			\caption{Micro time-step count comparison}
			\label{fig:O1_Steps}
		\end{subfigure}
		\begin{subfigure}{.5\textwidth}
			\centering
			% include fourth image
			\includegraphics[width=1\linewidth]{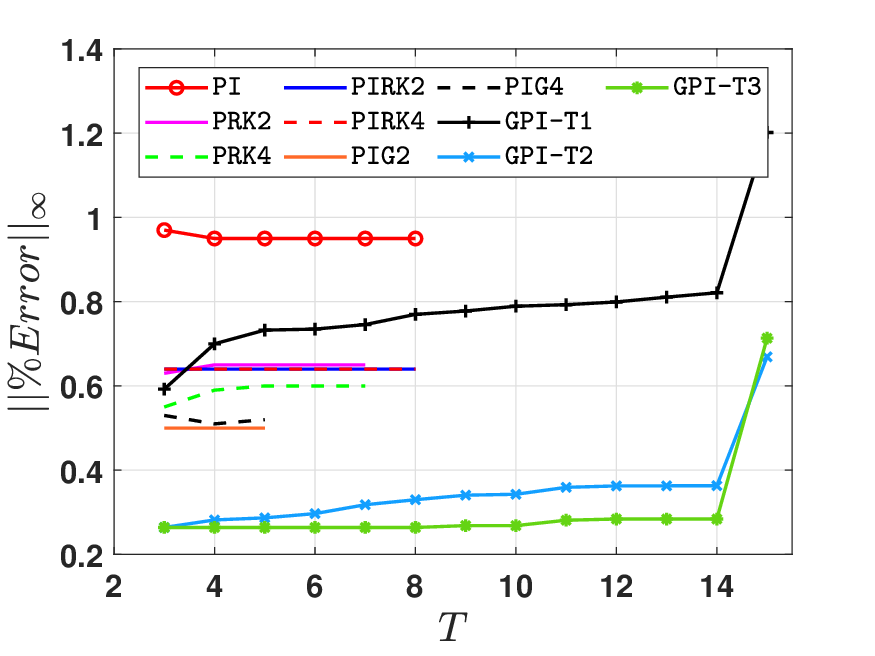}  
			\caption{Accuracy comparison}
			\label{fig:O1_Per_Err}
		\end{subfigure}
		
		\begin{subfigure}{.5\textwidth}
			\centering
			% include third image
			\includegraphics[width=1\linewidth]{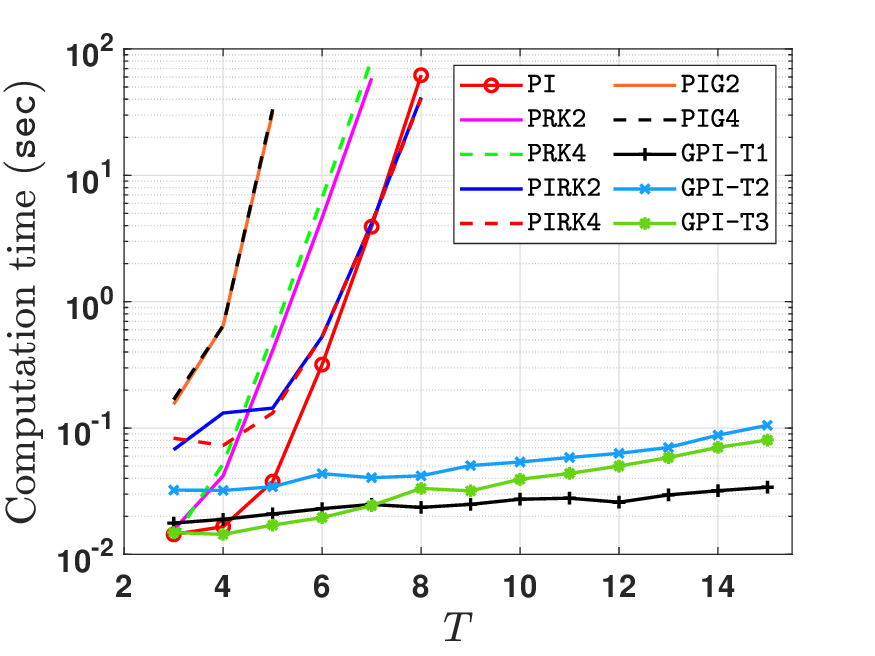}  
			\caption{Computational time comparison}
			\label{fig:O1_Time}
		\end{subfigure}
		\begin{subfigure}{.5\textwidth}
			\centering
			% include fourth image
			\includegraphics[width=1\linewidth]{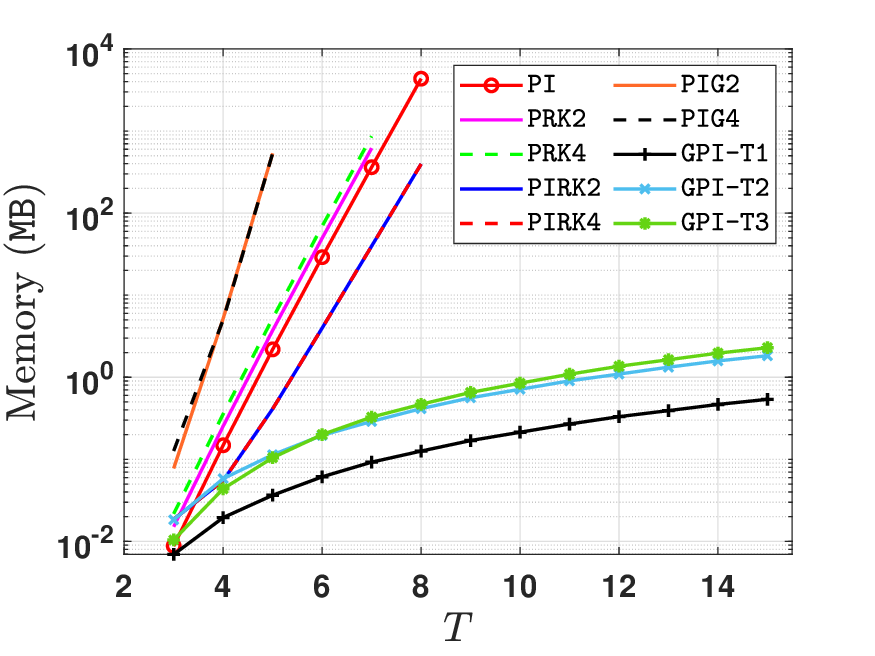}  
			\caption{Memory usage comparison}
			\label{fig:O1_Memory}
		\end{subfigure}
		\caption{Performance comparison of \texttt{PI}, \texttt{PRK2}, \texttt{PRK4}, \texttt{PIRK2}, \texttt{PIRK4}, \texttt{PIG2}, \texttt{PIG4}, \texttt{GPI-T1}, \texttt{GPI-T2} and \texttt{GPI-T3}.}
		\label{fig:Comparison_O1}
	\end{figure}

	The performance of the \texttt{PI} scheme is illustrated in Figure \ref{fig:Comparison_O1} in terms of total number of micro steps, accuracy, computational time and memory usage. The total number of micro time steps over the interval $[t_0=2,T]$ is shown in Figure \ref{fig:O1_Steps} for $T=3,4,5,6,7$ and $8$. As $T$ increases, the micro time step decreases to maintain stability, while the number of micro steps increases substantially to preserve the burst length. This results in high memory consumption and increased computational time, as depicted in Figures \ref{fig:O1_Memory} and \ref{fig:O1_Time}. The maximum percentage error, shown in Figure \ref{fig:O1_Per_Err}, remains approximately $0.95\%$ across all values of $T$. However, for larger final times $T>8$, the simulations cannot be completed using our computer due to excessive memory requirements.
	
	To address problems with time-dependent spectra using the \texttt{PI} scheme, a relatively long microsimulation is typically required. A similar observation was reported by Gear et al. \cite{2003_Gear_Projective}, where it was noted that the Euler–Lagrange formulation of such systems leads to operators whose spectra evolve over time. Although the separation of time scales may vary only mildly, a sufficiently large number of micro time steps is still necessary to effectively damp the fast components prior to extrapolation. This requirement increases the overall computational cost of the method.
	
	\paragraph{Projective Runge--Kutta schemes of order two (\texttt{PRK2}) and four (\texttt{PRK4}):}
	
	The projective Runge--Kutta schemes (\texttt{PRK}) \cite{2016_Lafitte_High-Order,2017_Lafitte_High-order} extend the \texttt{PI} framework to higher-order accuracy while retaining uniform macro--micro time steps and burst lengths. For the present problem \eqref{eqn:NonLinear_Stiff_ODE_Increase_Stiffness}, both \texttt{PRK2} and \texttt{PRK4} use the same macro and micro time steps as the \texttt{PI} scheme, namely $\Delta T=0.05$ and $\delta T=1\mathrm{e}{-T}$ for $t\in[t_0=2,T]$. Due to their higher-order structure, the fast variables relax over shorter burst lengths compared to the \texttt{PI} scheme. Specifically, \texttt{PRK2} and \texttt{PRK4} require approximately $10$ and $7$ micro steps per stage, respectively, to achieve relaxation. This corresponds to approximately $40\%$ and $56\%$ microsimulation of the macro time step or the entire interval $[t_0=2, T]$. The underlying Runge--Kutta methods used in \texttt{PRK2} and \texttt{PRK4} are defined by the following Butcher tableaux:
	\[
	\begin{array}{c|cc}
		0 & 0 & 0 \\
		1 & 1 & 0 \\
		\hline
		& \frac{1}{2} & \frac{1}{2}
	\end{array}
	\qquad\text{and}\qquad
	\begin{array}{c|cccc}
		0 & 0 & 0 & 0 & 0 \\
		\frac{1}{2} & \frac{1}{2} & 0 & 0 & 0 \\
		\frac{1}{2} & 0 & \frac{1}{2} & 0 & 0 \\
		1 & 0 & 0 & 1 & 0 \\
		\hline
		& \frac{1}{6} & \frac{1}{3} & \frac{1}{3} & \frac{1}{6}
	\end{array}
	\]
	
	Figure \ref{fig:Comparison_O1} illustrates the performance of the \texttt{PRK2} and \texttt{PRK4} schemes for $T=3,4,5,6,$ and $7$. Similar to the \texttt{PI} scheme, both methods require a large number of micro steps, as shown in Figure \ref{fig:O1_Steps}, leading to high computational cost in terms of time and memory (see Figures \ref{fig:O1_Time} and \ref{fig:O1_Memory}).
	The percentage errors for \texttt{PRK2} and \texttt{PRK4} are approximately $0.65\%$ and $0.60\%$, respectively, across all tested values of $T$, as shown in Figure \ref{fig:O1_Per_Err}. However, similar to the \texttt{PI} scheme, simulations fail for larger final times ($T>7$) due to excessive memory demands.

	\paragraph{Projective integration by second and fourth-order Runge--Kutta (\texttt{PIRK2} and \texttt{PIRK4}):}
	
	The \texttt{PIRK} schemes of second and fourth order are built using uniform macro time steps and non-uniform (adaptive) micro time steps, while maintaining a uniform micro burst length. In general, these schemes are not well-suited for solving the problem \eqref{eqn:NonLinear_Stiff_ODE_Increase_Stiffness}. In particular, \texttt{PIRK2} and \texttt{PIRK4} require the micro burst length to exceed $\frac{\Delta T}{2}$ and $\frac{\Delta T}{4}$, respectively, in order to obtain convergent solution, where $\Delta T$ denotes the uniform macro time step. This requirement implies that more than $100\%$ of each macro time step is effectively spent on microsimulation. This implies \texttt{PIRK} schemes do more microsimulation compared to the single scale explicit full time microsimulation schemes, making them inefficient for problems with time-dependent spectra such as \eqref{eqn:NonLinear_Stiff_ODE_Increase_Stiffness}.
	
	Nevertheless, we implement these schemes to assess their practical performance. A uniform macro time step $\Delta T = 0.05$ is used. The micro burst lengths are chosen as $0.0251$ and $0.0126$ for the \texttt{PIRK2} and \texttt{PIRK4} schemes, respectively. These choices correspond to approximately $100.4\%$ and $100.8\%$ micro-scale simulation within each macro time step over the entire time interval.
	
	The performance of the \texttt{PIRK} schemes is illustrated in Figure \ref{fig:Comparison_O1}. A notable feature of these schemes, absent in \texttt{PI} and \texttt{PRK}, is the ability to employ adaptive micro time stepping, for instance through the second-order \texttt{ode23} solver. This adaptability reduces the number of micro steps compared to \texttt{PI} and \texttt{PRK}, despite the excessive total micro simulation. As a result, \texttt{PIRK} schemes exhibit relatively lower memory usage.
	
	However, similar to \texttt{PI} and \texttt{PRK}, the total number of micro steps required by \texttt{PIRK} grow rapidly with the final time $T$, as shown in Figure \ref{fig:O1_Steps}. This leads to high computational cost in terms of both time and memory (see Figures \ref{fig:O1_Time} and \ref{fig:O1_Memory}). The percentage error remains approximately $0.64\%$ for $T=3,4,5,6,7,$ and $8$, as shown in Figure \ref{fig:O1_Per_Err}. For larger final times $T>8$, simulations cannot be completed due to excessive memory requirements.
	
	\paragraph{Projective integration via a general macroscale integrator of order two (\texttt{PIG2}) and four (\texttt{PIG4}):}
	
	The \texttt{PIG2} and \texttt{PIG4} schemes \cite{2021_maclean_toolbox,2005_Gear_Projecting} are higher-order methods that allow adaptive selection of both macro and micro time steps. However, the micro burst length is kept uniform across macro steps. In this study, it is fixed at $0.02$ for both the \texttt{PIG} schemes.
	
	Figure \ref{fig:Comparison_O1} presents the performance of the \texttt{PIG2} and \texttt{PIG4} schemes for final times $T=3,4,$ and $5$. Compared to \texttt{PI}, \texttt{PRK} and \texttt{PIRK}, these methods offer greater flexibility in selecting macro--micro time steps. However, they require the largest number of micro time steps among all the considered schemes, as shown in Figure \ref{fig:O1_Steps}. Consequently, they incur the highest computational cost in terms of both computational time and memory usage. In terms of accuracy, the \texttt{PIG} schemes outperform the other methods. The \texttt{PIG2} and \texttt{PIG4} schemes achieve errors of approximately $0.50\%$ and $0.52\%$, respectively, for the chosen burst length. Although increasing the burst length can further improve accuracy, it also significantly increases computational cost (in terms of micro steps, runtime and memory), making such choices impractical. For larger final times $T>5$, simulations with \texttt{PIG} schemes cannot be completed due to excessive memory demands.
	
	\vspace{0.5cm}
	
	For problems with time-dependent spectra, the performance of the existing projective integration schemes is neither satisfactory nor computationally efficient. In general, the application of classical single-scale numerical methods to multiscale problems becomes impractical due to their prohibitive computational cost. This highlights the need for robust multiscale schemes that can effectively balance computational efficiency and accuracy.
	
	From the preceding discussion, it is evident that none of the existing multiscale approaches--namely \texttt{PI}, \texttt{PRK2}, \texttt{PRK4}, \texttt{PIRK2}, \texttt{PIRK4}, \texttt{PIG2} and \texttt{PIG4}--achieve such a balance when applied to problems with time-dependent spectra along with scale separation. These methods lack certain adaptive features that are essential for handling dynamically evolving spectra and scale separation.
	
	\begin{table}
		\centering 
		\caption{Comparison of the features of three types of \texttt{GPI} schemes.}
		\begin{tabular}{ccccc}
			\hline
			\textbf{Scheme} & \textbf{Macro step} & \textbf{Meso step} & \textbf{Micro step} & \textbf{Burst length} \\
			\hline
			\texttt{GPI-T1} & fixed & N/A & variable & variable \\
			\texttt{GPI-T2} & fixed & variable & variable & variable \\
			\texttt{GPI-T3} & variable & N/A & variable & variable \\
			\hline
		\end{tabular}
		\label{table:Three_Types_GPI}
	\end{table}
	
	In the following, problem \eqref{eqn:NonLinear_Stiff_ODE_Increase_Stiffness} is solved using the proposed \texttt{GPI} scheme in three different configurations. These are referred to as type-1, type-2 and type-3 variants, denoted by \texttt{GPI-T1}, \texttt{GPI-T2} and \texttt{GPI-T3}, respectively. The key features of these variants are summarised in Table \ref{table:Three_Types_GPI}.
	
	For all three variants, variable micro time steps and adaptive burst lengths are employed to accommodate the time-dependent stiffness of the problem. The \texttt{GPI-T2} scheme additionally incorporates a variable meso time scale, bridging the macro and micro levels. While \texttt{GPI-T1} and \texttt{GPI-T2} use fixed macro time steps, the \texttt{GPI-T3} scheme allows for adaptive macro time stepping, providing greater flexibility in resolving the evolving dynamics.
	
	\paragraph{\texttt{GPI} scheme of type-1 (\texttt{GPI-T1}):}
	
	In the \texttt{GPI-T1} scheme, non-uniform micro time steps are employed following the strategy described in Subsection \ref{sec:Choice_Micro_Steps}, while a uniform macro time step $\Delta T = 0.05$ is used, consistent with the settings in the \texttt{PI}, \texttt{PRK} and \texttt{PIRK} schemes. A key feature of the \texttt{GPI} framework is the ability to adapt the micro burst length based on local problem characteristics. In this variant, the burst length is chosen as
	\begin{center}
		$\frac{\beta}{10^{T^{n}}}\log\!\left(\Delta T \, 10^{T^{n}}\right)$,
	\end{center}
	with burst control parameter $\beta=2$, motivated by the relaxation bound in \eqref{eqn:Relaxation_Bound_Increasing_Stiffness}. 
	
	\begin{table}
		\centering 
		\caption{Maximum percentage error of \texttt{GPI-T1} solution for different values of the burst control parameter $\beta$ over the interval $[2,14]$.}
		\begin{tabular}{ccccc}
			\hline
			$\beta$ & 2 & 2.2 & 2.4 & 2.6 \\
			\textbf{Maximum \%Error} & 0.82 & 0.79 & 0.75 & 0.72\\
			\hline
		\end{tabular}
		\label{table:GPI_T1_Relaxed}
	\end{table}
	
	Table \ref{table:GPI_T1_Relaxed} examines the sensitivity of the \texttt{GPI-T1} scheme to the choice of burst length, controlled by the parameter $\beta$. The results are obtained by varying $\beta$ while keeping all other parameters fixed over the interval $[2,14]$. Increasing $\beta$ (and hence the burst length) leads to only marginal improvements in accuracy. This indicates that even for $\beta = 2$, the micro solver sufficiently relaxes the fast dynamics onto the slow manifold. Therefore, the proposed burst length selection is both efficient and consistent with the theoretical relaxation estimate.
	
	This adaptive choice significantly reduces the number of required micro steps, as illustrated in Figure \ref{fig:O1_Steps}. Consequently, \texttt{GPI-T1} achieves substantially lower memory usage and computational time compared to existing projective integration schemes (see Figures \ref{fig:O1_Memory} and \ref{fig:O1_Time}). Over the range $T \in [3,15]$, the maximum computational time is $0.034$ second and the peak memory usage is $0.54$ MB, demonstrating a marked improvement in efficiency.
	
	As time progresses, the burst length decreases under fixed macro time stepping, which effectively increases the projective extrapolation size. Since the accuracy of projective integration methods depends on the extrapolation length, a gradual increase in the error is observed with increasing $T$. Nevertheless, the method maintains good accuracy: the maximum percentage error remains within the interval $[0.61, 0.75]$ for $T = 3,4,\ldots,14$.

	This loss of accuracy, induced by the growth of the extrapolation interval, can be mitigated in two ways. One approach is to introduce an intermediate (meso) time scale in regions where the extrapolation interval becomes large, thereby reducing the effective projection length. Alternatively, the extrapolation interval can be explicitly controlled by fixing its size throughout the computation. These two strategies are incorporated in the \texttt{GPI-T2} and \texttt{GPI-T3} schemes, respectively.
	
	\paragraph{\texttt{GPI} scheme of type-2 (\texttt{GPI-T2}):}
	
	The \texttt{GPI-T2} scheme incorporates three time scales: macro, meso and micro. The macro time steps are uniform, except for the first step, while both meso and micro time steps are chosen adaptively. The micro time steps follow the strategy described in Subsection \ref{sec:Choice_Micro_Steps}. The micro burst length within each meso step is given by
	\begin{center}
		$\frac{\beta}{10^{T^{n,m}}}\log\!\left(\partial T^{n,m} \, 10^{T^{n,m}}\right)$,
	\end{center}
	with burst control parameter $\beta = 2$, based on the relaxation bound \eqref{eqn:Relaxation_Bound_Increasing_Stiffness}. With this choice, the micro time steps in the first macro step are approximately $9.77\mathrm{e}{-3}$, $9.56\mathrm{e}{-3}$, $9.35\mathrm{e}{-3}$ and $9.15\mathrm{e}{-3}$. The first projective extrapolation step is fixed at $0.05$. Thereafter, uniform macro time stepping is employed, with $20$ macro steps per unit time, consistent with the previous schemes. Consequently, the total number of macro steps in \texttt{GPI-T2} matches that of \texttt{GPI-T1}.
	
	Meso time stepping is introduced adaptively to control the extrapolation interval. Specifically, no meso steps are used in the first seven macro steps. In macro steps 8--15, two meso steps of equal size are used, while for subsequent macro steps, three meso steps are employed of equal size.

	%We set the first projective extrapolation step size is 0.05 in the first macro time step. For the rest of the time, uniform sizes of macro time steps are considered. As the previously discussed schemes, in each unit time step, 20 macro time steps are considered. The total number of macro time steps in \texttt{GPI-T1} and \texttt{GPI-T2} are identical. In the first seven macro time steps, meso time steps are not being used. However, for eighth to fiftinth macro time step, two meso time steps are used each of size half of its macro time step. For the rest of the macro time steps, three meso time steps are used in each macro time step of size one-third of the macro time step. 
	
	The inclusion of meso time steps leads to a moderate increase in the number of micro steps compared to \texttt{GPI-T1} (see Figure \ref{fig:O1_Steps}), resulting in slightly higher memory usage and computational time (Figures \ref{fig:O1_Memory} and \ref{fig:O1_Time}). However, this additional cost is offset by a significant improvement in accuracy. The maximum percentage error remains within the interval $[0.26, 0.36]$ for $T = 3,4,\ldots,14$. In addition, the maximum computational time and peak memory consumption are only $0.088$ second and $1.58$ MB, respectively, highlighting the computational efficiency of the proposed method. In particular, \texttt{GPI-T2} achieves higher accuracy than \texttt{GPI-T1} and the other existing schemes. Moreover, the use of meso time steps effectively controls the growth of the extrapolation interval. As a result, the extrapolation step increases at a slower rate compared to \texttt{GPI-T1}, leading to only a mild increase in the percentage error as time progresses.
	
	\begin{table}
		\centering 
		\caption{Maximum percentage error of \texttt{GPI-T2} for different values of the burst control parameter $\beta$ over the interval $[2,14]$.}
		\begin{tabular}{ccccc}
			\hline
			$\beta$ & 2 & 2.2 & 2.4 & 2.6 \\
			\textbf{Maximum \%Error} & 0.36 & 0.33 & 0.33 & 0.30\\
			\hline
		\end{tabular}
		\label{table:GPI_T2_Relaxed}
	\end{table}
	
	Table \ref{table:GPI_T2_Relaxed} examines the sensitivity of the \texttt{GPI-T2} scheme to the burst control parameter $\beta$. The results are obtained by varying $\beta$ while keeping all other parameters fixed over the interval $[2,14]$. Increasing the burst length yields only marginal improvements in accuracy, indicating that the fast dynamics are already sufficiently relaxed onto the slow manifold for $\beta = 2$. This confirms the effectiveness of the proposed burst length selection.
	
	\paragraph{\texttt{GPI} scheme of type-3 (\texttt{GPI-T3}):}
	
	The \texttt{GPI-T3} scheme employs two time scales--macro and micro--to solve problem \eqref{eqn:NonLinear_Stiff_ODE_Increase_Stiffness}. The first macro step is treated in the same manner as in the \texttt{GPI-T2} scheme. Thereafter, both macro and micro time steps are selected adaptively. In particular, the macro step sizes, micro step sizes and micro burst lengths are all allowed to vary throughout the simulation. The micro time steps are determined using \eqref{eqn:Variable_Inner_Steps_Formula}, while the burst length is chosen as
	\begin{center}
		$\frac{\beta}{10^{T^{n}}}\log\!\left(\Delta T^{n-1} \, 10^{T^{n}}\right)$,
	\end{center}
	with burst control parameter $\beta = 2$. Since the current macro step $\Delta T^n$ is not known a priori, the previously computed value $\Delta T^{n-1}$ is used for $n \geq 2$. Following the micro simulation, projective extrapolation is performed with a fixed step size of $0.01$ at each macro step (except the first), similar to the \texttt{PI} scheme.
	
	Compared to \texttt{GPI-T1} and \texttt{GPI-T2}, the \texttt{GPI-T3} scheme requires a larger number of micro steps. However, its computational time and memory usage remain comparable to those of \texttt{GPI-T2}. The maximum computational time and peak memory consumption are only $0.070$ second and $1.97$ MB, respectively, highlighting the computational efficiency of the proposed method. The maximum percentage error remains within the interval $[0.26, 0.28]$ for $T = 3,4,\ldots,14$.  In terms of accuracy, \texttt{GPI-T3} outperforms both \texttt{GPI-T1}, \texttt{GPI-T2}, as well as the other existing schemes. This improvement is primarily due to the use of a fixed extrapolation step size, which prevents the growth of the projection interval over time. Consequently, nearly uniform accuracy is maintained across the range $T \in [3,14]$.
	
	In contrast to \texttt{GPI-T1} and \texttt{GPI-T2}, where the extrapolation interval increases with time, the \texttt{GPI-T3} scheme maintains a constant projection length (except in the first step), resulting in stable error behaviour.
	
	\begin{table}
		\centering 
		\caption{Maximum percentage error of \texttt{GPI-T3} for different values of the burst control parameter $\beta$ over the interval $[2,14]$.}
		\begin{tabular}{ccccc}
			\hline
			$\beta$ & 2 & 2.2 & 2.4 & 2.6 \\
			\textbf{Maximum \%Error} & 0.28 & 0.28 & 0.28 & 0.27\\
			\hline
		\end{tabular}
		\label{table:GPI_T3_Relaxed}
	\end{table}
	
	Table \ref{table:GPI_T3_Relaxed} evaluates the sensitivity of the \texttt{GPI-T3} scheme to the burst control parameter $\beta$. The results are obtained by varying $\beta$ while keeping all other parameters fixed over the interval $[2,14]$. Increase in the burst length produces negligible improvement in accuracy, indicating that the fast dynamics are already sufficiently relaxed onto the slow manifold for $\beta = 2$. This further validates the effectiveness of the proposed burst length selection.
	
	\vspace{0.5cm}
	
	For large final times $T \leq 14$, none of the three \texttt{GPI} variants exhibits a noticeable increase in error. In particular, when $T=15$, the maximum percentage error rises to approximately $1.20\%$, $0.67\%$ and $0.71\%$ for \texttt{GPI-T1}, \texttt{GPI-T2} and \texttt{GPI-T3}, respectively. This behaviour can be attributed to finite-precision limitations. At $T=15$, the fast eigenvalue is approximately $-10^{15}$, causing the micro time steps to fall below $10^{-15}$. Since double-precision arithmetic (as used in \texttt{MATLAB}) provides only about 15--16 digits of accuracy, round-off errors become dominant at this scale, leading to a degradation in numerical accuracy.
	
	\paragraph{Overall comparison:}
	
\begin{figure}
		\centering
		\includegraphics[width=0.6\linewidth]{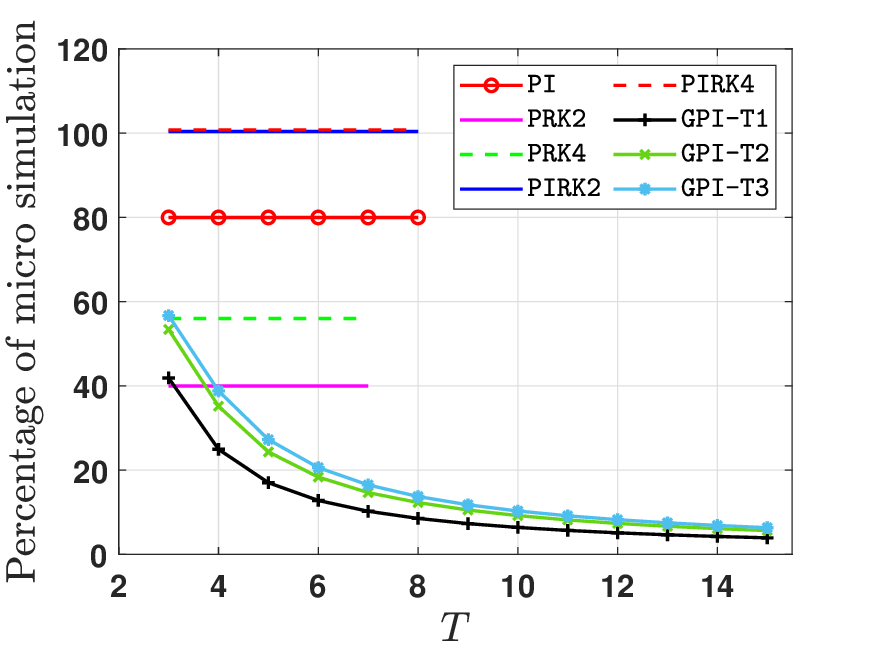}  
		\caption{Percentage of micro-scale simulation over the interval $[2,T]$ as the final time $T$ increases for problem \eqref{eqn:NonLinear_Stiff_ODE_Increase_Stiffness}. Unlike the other methods, the \texttt{GPI} schemes exhibit a decreasing trend, indicating reduced reliance on computationally expensive micro-scale simulations.}
		\label{fig:O1_Comparison_Per_Mic_Sim}
	\end{figure}

	The comparative study of the methods—\texttt{PI}, \texttt{PRK}, \texttt{PIRK}, \texttt{PIG} and the proposed \texttt{GPI} variants—reveals clear and consistent differences in performance across all key metrics, including micro time step count, accuracy, computational time, memory usage and the proportion of micro-scale simulations (as discussed below and illustrated in Figure~\ref{fig:O1_Comparison_Per_Mic_Sim}). Together, these metrics provide a comprehensive assessment of efficiency and scalability.
	
	A primary observation from the results in Figure \ref{fig:O1_Steps} is the markedly different growth behaviour in the number of micro time steps required by each scheme. The \texttt{PI} and \texttt{PRK} methods exhibit a steady increase in the number of time steps as the final time $T$ increases, while the \texttt{PIG} schemes demonstrate a much steeper growth, indicating poor scalability. In contrast, the \texttt{PIRK} methods reduce this growth to some extent; however, the \texttt{GPI} schemes exhibit the slowest growth among all methods. This reduced step count demonstrates that the \texttt{GPI} framework advances the solution more efficiently over long time intervals, thereby lowering the overall computational cost.
	
	In terms of accuracy, as mentioned in Figure \ref{fig:O1_Per_Err}, the \texttt{GPI} schemes, particularly \texttt{GPI-T2} and \texttt{GPI-T3}, consistently achieve lower percentage errors compared to \texttt{PI}, \texttt{PRK} and \texttt{PIRK} methods across the tested range of $T$. Although the \texttt{PIG} schemes also yield competitive accuracy, this comes at the expense of substantially higher computational cost. The \texttt{GPI} methods therefore achieve a superior balance between accuracy and efficiency, making them suitable for reliable long-time integration of multiscale systems.
	
	The computational time results in Figure \ref{fig:O1_Time} further highlight the advantage of the proposed \texttt{GPI} scheme. While \texttt{PI}, \texttt{PRK}, \texttt{PIRK} and \texttt{PIG} methods exhibit rapid growth in computational time as $T$ increases, the \texttt{GPI} schemes show an almost negligible increase. This near-constant computational time indicates that the \texttt{GPI} schemes effectively reduce their reliance on micro-scale simulations as the system evolves, enabling efficient long-time integration. 
	
	A similar trend is observed in memory usage as shown in Figure \ref{fig:O1_Memory}. The existing methods \texttt{PI}, \texttt{PRK}, \texttt{PIRK} and \texttt{PIG} require rapidly increasing memory as $T$ grows, which becomes a limiting factor for large-scale simulations. In contrast, the \texttt{GPI} schemes maintain low and slowly varying memory requirements, reflecting efficient data management and reduced storage of micro-scale information.
	
	Another critical aspect of multiscale methods is the proportion of micro simulations as discussed in Figure \ref{fig:O1_Comparison_Per_Mic_Sim}. The results show that the percentage of micro simulations remains relatively high and often constant for \texttt{PI}, \texttt{PRK} and \texttt{PIRK} methods, while for \texttt{PIG} it grows excessively. Due to this extreme growth, the \texttt{PIG} results are omitted from the figure. For instance, \texttt{PIG2} requires approximately $22\%$, $177\%$ and $1541\%$ micro simulation for $T=3,4,$ and $5$, respectively, while \texttt{PIG4} requires $80\%$, $380\%$ and $3421\%$. In contrast, the \texttt{GPI} schemes demonstrate a decreasing trend in the percentage of micro simulations as $T$ increases. This indicates that the \texttt{GPI} framework leads to significant gains in computational efficiency.
	
	Overall, the \texttt{GPI} framework remains computationally feasible even for large final times, while the other methods become prohibitively expensive in computational cost, rendering them impractical on the available computational platform. The proposed \texttt{GPI} scheme consistently outperforms the existing approaches across all the considered metrics. They combine low growth in step count, high accuracy, near-constant computational time, reduced memory usage and diminishing reliance on micro-scale simulation. These characteristics make the \texttt{GPI} framework well-suited for time-dependent spectra-based multiscale problems for long-time integration.

	\subsubsection{Performance of the \texttt{GPI} scheme is compared with standard stiff \texttt{ODE} solvers}\label{subsubsec:Performance_GPI_compared_to_ode_solvers}
	
	In this subsection, we compare the performance of the proposed \texttt{GPI} scheme with widely used stiff \texttt{ODE} solvers, namely \texttt{ode15s} and \texttt{Radau IIA}. For a fair comparison, all methods are considered at first-order accuracy: the \texttt{GPI} scheme is an explicit first-order method, while \texttt{ode15s} and \texttt{Radau IIA} are treated in their first-order implicit configurations.
	
	The comparison is conducted across three regimes of stiffness intensity:
	\begin{enumerate}
		\item low stiffness, where $f(t)\in[10^2, 10^3]$ for $t\in [2, 3]$,
		\item medium stiffness, where $f(t)\in[10^7, 10^8]$ for $t\in [7, 8]$ and
		\item high stiffness, where $f(t)\in[10^{11}, 10^{12}]$ for $t\in [11, 12]$.
	\end{enumerate}
	This setup enables a systematic assessment of performance as the stiffness increases.
	
	The results are summarised in Table \ref{table:GPIvsOther_Three_Types_Stiffness}. The second column, $\operatorname{N_t}$, denotes the number of macro time steps for the \texttt{GPI} scheme and the total number of time steps for the other solvers. Accuracy is reported in terms of maximum percentage error, while computational time (in seconds) and memory usage (in MB) are provided in the fourth and fifth columns, respectively. 
	
	In this study, the \texttt{GPI-T1} variant is used, with macro-micro time steps and burst lengths chosen as described in Subsection \ref{subsubsec:Non_Linear_Increasing_Stiffness}.
	
	\begin{table}
		\centering
		\caption{Performance comparison of some first order schemes such as \texttt{GPI}, \texttt{ode15s} and \texttt{Radau IIA} across different stiffness regimes.}
		\begin{tabular}{ccccc}%{p{4cm}p{0.5cm}p{2cm}p{2cm}p{2.5cm}}
			\hline
			Solver&$\operatorname{N_t}$&$||\%Error||_\infty$&Time (\texttt{sec})&Memory (\texttt{MB})\\
			\hline
			\hline
			\multicolumn{5}{c}{Low stiffness} \\
			\hline
			\texttt{GPI}&20&0.59&0.017&0.0013\\
			\hline
			\texttt{ode15s}, adaptive&103&0.79&0.073&0.0065\\
			\texttt{ode15s}, non-adaptive&20&0.79&0.080&0.0048\\
			\hline
			\texttt{Radau IIA}, adaptive&413&0.64&0.14&0.016\\
			\texttt{Radau IIA}, non-adaptive&20&0.64&0.15&0.0067\\
			\hline
			\hline
			\multicolumn{5}{c}{Medium stiffness} \\
			\hline
			\texttt{GPI}&150&0.19&0.022&0.0046\\
			\hline
			\texttt{ode15s}, adaptive&121&0.66&0.079&0.0069\\
			\texttt{ode15s}, non-adaptive&150&0.67&0.091&0.0088\\
			\hline
			\texttt{Radau IIA}, adaptive&434&0.24&0.15&0.016\\
			\texttt{Radau IIA}, non-adaptive&150&0.24&0.15&0.011\\
			\hline
			\hline
			\multicolumn{5}{c}{High stiffness} \\
			\hline
			\texttt{GPI}&150&0.19&0.025&0.0047\\
			\hline
			\texttt{ode15s}, adaptive&130&0.68&0.073&0.0072\\
			\texttt{ode15s}, non-adaptive&150&0.69&0.088&0.0088\\
			\hline
		\end{tabular}
		\label{table:GPIvsOther_Three_Types_Stiffness}
	\end{table}	
	
	The \texttt{GPI} scheme is non-adaptive, whereas both \texttt{ode15s} and \texttt{Radau IIA} can be implemented in either adaptive or non-adaptive modes. For a consistent comparison, the non-adaptive configurations in Table \ref{table:GPIvsOther_Three_Types_Stiffness} use $\operatorname{N_t}=20$, $150$ and $150$ for low, medium and high stiffness regimes, respectively. Furthermore, all methods--\texttt{GPI}, \texttt{ode15s} and \texttt{Radau IIA}--are considered at first-order accuracy to ensure a fair comparison. Across all metrics, including step count, accuracy, computational time and memory usage, the \texttt{GPI} scheme demonstrates superior performance.
	
	In the low-stiffness regime, all non-adaptive schemes use $\operatorname{N_t}=20$. The \texttt{GPI} scheme achieves approximately $75\%$ and $90\%$ of the error of the non-adaptive \texttt{ode15s} and \texttt{Radau IIA} methods, respectively. At the same time, it requires less than one-fourth and one-ninth of their computational time and approximately one-fourth and one-fifth of their memory usage. Similar trends are observed when compared with the adaptive variants of \texttt{ode15s} and \texttt{Radau IIA}. Notably, the adaptive solvers require more than five and twenty times the number of time steps, respectively, without achieving any meaningful improvement in performance.
	
	For the medium-stiffness regime with $\operatorname{N_t}=150$, the advantages of the \texttt{GPI} scheme become more pronounced. Compared to the non-adaptive \texttt{ode15s}, \texttt{GPI} achieves approximately one-third of the error while requiring about one-fourth of the computational time and roughly half of the memory. Compared to \texttt{Radau IIA}, \texttt{GPI} attains comparable accuracy (approximately five-sixths of the error) while reducing the computational time to about one-seventh and halving the memory usage.
	
	A similar trend is observed when comparing \texttt{GPI} with the adaptive versions of \texttt{ode15s} and \texttt{Radau IIA}. In particular, the adaptive \texttt{Radau IIA} scheme requires approximately three times more time steps than its non-adaptive counterpart, yet does not yield a corresponding improvement in accuracy or efficiency.

	Similarly, in the high-stiffness regime, the \texttt{GPI} scheme continues to outperform both the adaptive and non-adaptive versions of \texttt{ode15s}. In contrast, the first-order \texttt{Radau IIA} scheme fails to produce results under such extreme stiffness conditions.

	\begin{figure}
		\begin{subfigure}{.5\textwidth}
			\centering
			% include third image
			\includegraphics[width=1\linewidth]{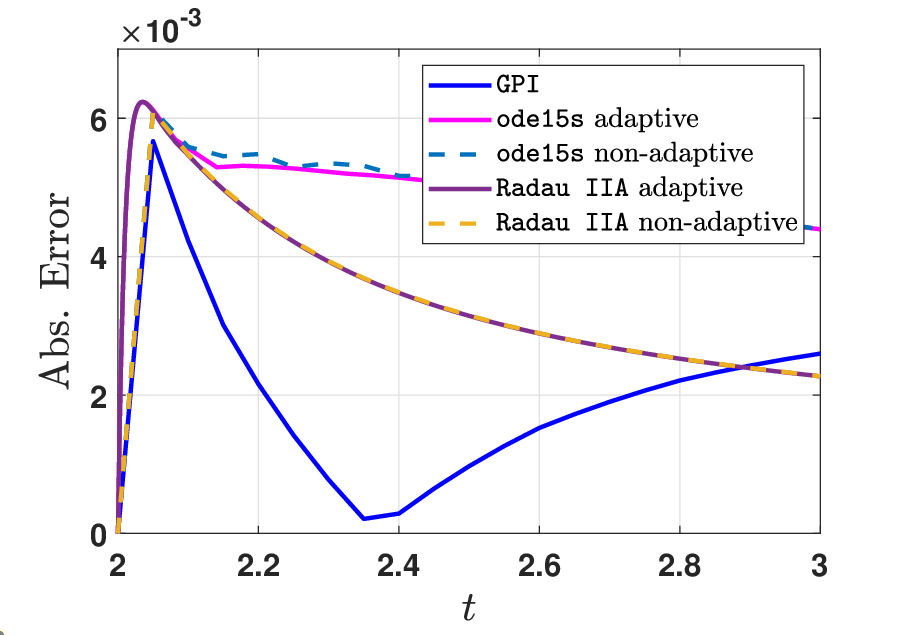}  
			\caption{Low stiffness}
			\label{fig:Low_stiffness_O1}
		\end{subfigure}
		\begin{subfigure}{.5\textwidth}
			\centering
			% include fourth image
			\includegraphics[width=1\linewidth]{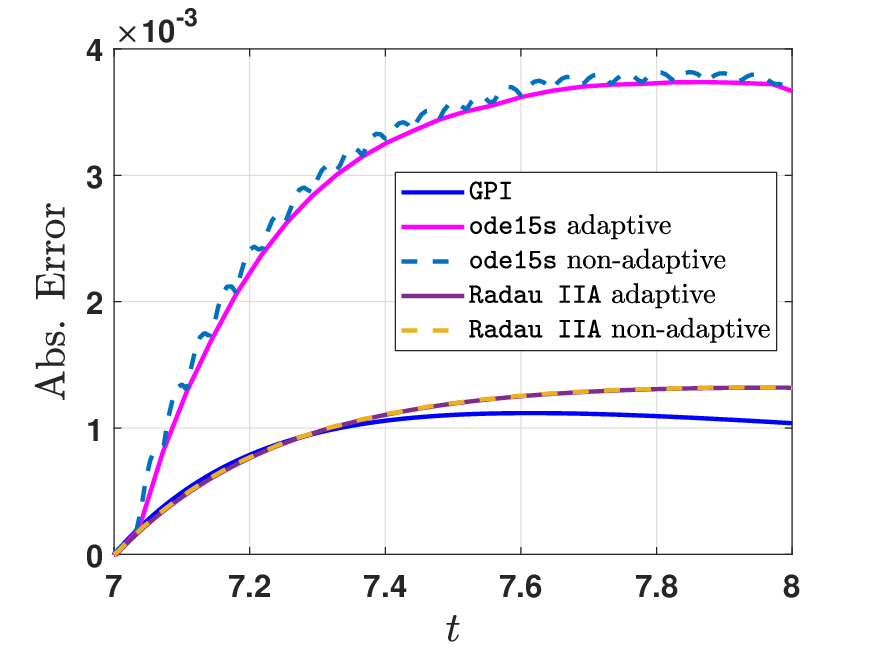}  
			\caption{Medium stiffness}
			\label{fig:Medium_stiffness_O1}
		\end{subfigure}
		
		\centering
		\begin{subfigure}{.5\textwidth}
			\centering
			% include fourth image
			\includegraphics[width=1\linewidth]{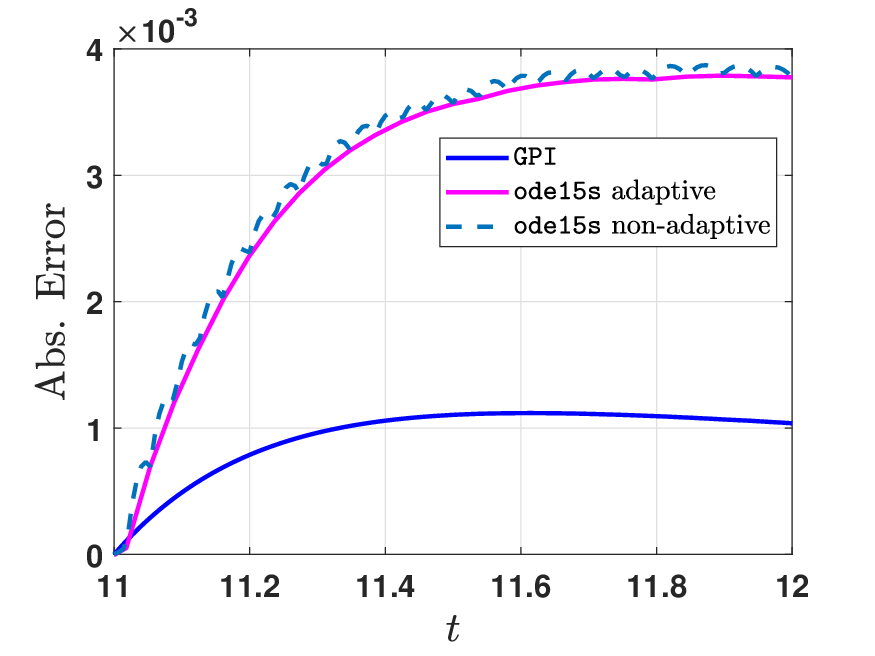}  
			\caption{High stiffness}
			\label{fig:High_stiffness_O1}
		\end{subfigure}
		\caption{A comparison of the absolute errors of the solutions obtained using the \texttt{GPI}, \texttt{ode15s} and \texttt{Radau IIA} schemes for low, medium and high stiffness variations.}
		\label{fig:Three_Types_stiffness_O1}
	\end{figure}
	
	A visual comparison of the absolute errors for the \texttt{GPI}, \texttt{ode15s} and \texttt{Radau IIA} solutions is presented in Figure \ref{fig:Three_Types_stiffness_O1}. Subfigures \ref{fig:Low_stiffness_O1}, \ref{fig:Medium_stiffness_O1} and \ref{fig:High_stiffness_O1} display the absolute errors corresponding to low, medium and high stiffness regimes, respectively. In all three cases, the \texttt{GPI} scheme yields consistently lower absolute errors than both \texttt{ode15s} and \texttt{Radau IIA}.
	
	In the low-stiffness case in Subfigure \ref{fig:Low_stiffness_O1}, a slight increase in the absolute error of the \texttt{GPI} solution is observed near the final time $T=3$; however, the error decreases again for $t>3.5$, indicating stable long-time behaviour.
	
	The results summarised in Table \ref{table:GPIvsOther_Three_Types_Stiffness} and Figure \ref{fig:Three_Types_stiffness_O1} collectively demonstrate that the \texttt{GPI} scheme provides the best overall performance among the considered first-order methods. In particular, it achieves a favourable combination of low step count, high accuracy, minimal computational time and reduced memory usage compared to both adaptive and non-adaptive implementations of \texttt{ode15s} and \texttt{Radau IIA}.

	%%%%%%%%%%%%%%%%%%%%%%%%%%%%%%%%%%%%%%%%%%%%%%%%%%%%%%%%%%%%%%%%%%%%%%%%%%%%%%%%

	\subsection{A pure diffusion equation with time-dependent diffusivity}\label{subsec:Pure_Diff}
	
	We consider the one-dimensional diffusion equation with exponentially time-dependent diffusivity, as discussed by Shampine et al. \cite{1997_Shampine_MATLAB}:
	\begin{equation}\label{eqn:Diff_PDE}
		u_t=\exp(t) u_{xx}, \quad x\in [0,\pi],\quad t\in [0,T],
	\end{equation}
	subject to the Dirichlet boundary conditions  $u(0,t)=0$ and $u(\pi,t)=0$. The initial condition is given by
	\begin{equation}\label{eqn:Diff_IC}
		u(x,0)=\sin(x).
	\end{equation}
	The problem is spatially discretised using a finite difference scheme. A fourth-order central difference scheme is employed at the interior grid points, except for the two points adjacent to each boundary, where a second-order central scheme is used. After spatial discretisation, the \texttt{PDE} $u_t=u_{xx}$ reduces to a system of linear \texttt{ODEs} with a constant coefficient matrix. Consequently, the system has time-independent eigenvalues $\lambda_p<0$ for $p=1,\ldots,\operatorname{N_x}-1$, where $\operatorname{N_x}$ denotes the total number of spatial grids. For the time-dependent problem $u_t=\exp(t)u_{xx}$, the eigenvalues of the resulting system become $\exp(t)\lambda_p<0$ for $p=1,\ldots,\operatorname{N_x}-1$. In contrast to the problem \eqref{eqn:NonLinear_Stiff_ODE_Increase_Stiffness}, where only the fast eigenvalue varies significantly while the slow eigenvalue remains close to zero, here the entire spectrum shifts along the negative real axis as time progresses. 
	
	To solve the problem \eqref{eqn:Diff_PDE}, it is essential to employ variable macro time steps, micro time steps and micro burst lengths. Existing projective integration frameworks based schemes such as \texttt{PI}, \texttt{PRK}, \texttt{PIRK} and \texttt{PIG} do not simultaneously incorporate all three features. Although the \texttt{PIG} scheme allows adaptive macro and micro time steps, however it uses a fixed micro burst length. Due to its high computational cost, it is not well-suited for this class of problems. As demonstrated in Subsection \ref{subsubsec:Non_Linear_Increasing_Stiffness}, even when a single eigenvalue varies in time, these existing schemes perform poorly; in the present case, where the entire spectrum varies, a similar performance is experienced. Therefore, we do not further discuss their failure for this problem. 
	
	We set the spatial grid size as $\Delta x=\frac{\pi}{20}$ over the domain $[0,\pi]$. 
	According to Lafitte et al. \cite{2017_Lafitte_High-order,2016_Lafitte_High-Order}, stable macroscale integration of one-dimensional parabolic \texttt{PDEs} requires that the macro time step must be proportional to the square of the macro spatial step size.
	Accordingly, in the \texttt{GPI} scheme, we consider non-uniform macro time steps of the form 
	\begin{center}
		$\Delta T^n=\omega\exp(-T^n)\Delta x^2$,
	\end{center} 
	where $\omega=0.05$. The micro time steps are chosen according to \eqref{eqn:Variable_Inner_Steps_Formula}, with $\lambda_*(t)=\tilde{\lambda}\exp(t)$, where $\tilde{\lambda}=2500$ is selected to exceed the maximum magnitude of the eigenvalues of the semi-discretised system corresponding to $u_t=u_{xx}$, thereby ensuring sufficiently fine micro time steps. The micro burst length is chosen as
	\begin{center}
		$\frac{\beta}{\lambda_*(T^n)}\log\!\left(\Delta T^{n} \, \lambda_*(T^n)\right)$,
	\end{center}
	with burst control parameter $\beta = 2$. All underlying solvers, including \texttt{ode15s}, \texttt{Radau IIA}, \texttt{ode23s}, \texttt{ode23t} and \texttt{ode23tb}, are implicit methods and both adaptive and non-adaptive time-stepping strategies are considered. Among these, \texttt{GPI}, \texttt{Radau IIA} and \texttt{ode15s} are first-order accurate, whereas \texttt{ode23s}, \texttt{ode23t} and \texttt{ode23tb} are second-order accurate.

	\begin{table}
		\centering
		\caption{A comparison of the accuracy of the \texttt{GPI}, \texttt{Radau IIA}, \texttt{ode15s}, \texttt{ode23s}, \texttt{ode23t} and \texttt{ode23tb} schemes are presented for various time intervals $[0,T]$ for the problem \eqref{eqn:Diff_PDE}.}
		\begin{tabular}{cccc|ccc}%{p{1.5cm}p{1.5cm}p{2cm}p{1.5cm}|p{1.5cm}p{1.5cm}p{1.5cm}}
			\hline
			&\multicolumn{6}{c}{Maximum relative errors} \\
			\hline
			$T$&\texttt{GPI}&\texttt{Radau IIA}&\texttt{ode15s}&\texttt{ode23s}&\texttt{ode23t}&\texttt{ode23tb}\\
			\hline
			0.5&4.3$\mathrm{e}{-5}$&1.4$\mathrm{e}{-3}$&3.6$\mathrm{e}{-3}$&1.2$\mathrm{e}{-4}$&2.5$\mathrm{e}{-4}$&1.4$\mathrm{e}{-4}$\\
			1&1.2$\mathrm{e}{-4}$&5.0$\mathrm{e}{-3}$&1.8$\mathrm{e}{-2}$&2.7$\mathrm{e}{-3}$&1.3$\mathrm{e}{-3}$&1.4$\mathrm{e}{-3}$\\
			1.5&3.5$\mathrm{e}{-4}$&1.2$\mathrm{e}{-2}$&4.3$\mathrm{e}{-2}$&1.0$\mathrm{e}{-2}$&3.5$\mathrm{e}{-3}$&1.5$\mathrm{e}{-3}$\\
			2  &7.9$\mathrm{e}{-4}$&2.6$\mathrm{e}{-2}$&9.0$\mathrm{e}{-2}$&1.9$\mathrm{e}{-2}$&1.1$\mathrm{e}{-2}$&8.8$\mathrm{e}{-3}$\\
			2.5&1.6$\mathrm{e}{-3}$&1.1$\mathrm{e}{-1}$&2.9$\mathrm{e}{-1}$&6.0$\mathrm{e}{-2}$&8.9$\mathrm{e}{-2}$&3.8$\mathrm{e}{-2}$\\
			3  &3.0$\mathrm{e}{-3}$&failure&failure&failure&failure&failure\\
			3.5 &5.3$\mathrm{e}{-3}$&failure&failure&failure&failure&failure\\
			4 &9.2$\mathrm{e}{-3}$&failure&failure&failure&failure&failure\\
			4.5 &1.6$\mathrm{e}{-2}$&failure&failure&failure&failure&failure\\
			5 &2.6$\mathrm{e}{-2}$&failure&failure&failure&failure&failure\\
			\hline
		\end{tabular}
		\label{table:Pure_Diff_Exp}
	\end{table}	
	
	Table \ref{table:Pure_Diff_Exp} presents the accuracy of the aforementioned schemes in terms of maximum relative error. For various final times $T$, we solve the problem \eqref{eqn:Diff_PDE} and examine the performance of each method compared to other methods. The \texttt{GPI} scheme has a better accuracy compared to all other mentioned schemes. For example, over the spatial domain  $[0, \pi]$ and time interval $[0, 2.5]$, the maximum relative error of the \texttt{GPI} solution is approximately $\frac{1}{69}$, $\frac{1}{181}$, $\frac{1}{38}$, $\frac{1}{56}$ and $\frac{1}{24}$ of the maximum relative errors produced by \texttt{Radau IIA}, \texttt{ode15s}, \texttt{ode23s}, \texttt{ode23t} and \texttt{ode23tb}, respectively. 
	
	Despite being a first-order method, $\texttt{GPI}$ achieves higher accuracy than the second-order methods \texttt{ode23s}, \texttt{ode23t}, and \texttt{ode23tb}. Furthermore, these competing schemes fail to maintain both accuracy and stability for final times $T \geq 3$, whereas the proposed \texttt{GPI} scheme continues to produce accurate and stable results. For the adaptive implementations of the existing methods, the observed accuracy is similar to that of the non-adaptive cases reported in Table~\ref{table:Pure_Diff_Exp}. Finally, for this problem \eqref{eqn:Diff_PDE}, the \texttt{IMEX} formulation coincides with fully implicit schemes, since there is no non-stiff component to separate. Therefore, the $\texttt{GPI}$ scheme can also be regarded as outperforming the \texttt{IMEX} approach.
	
	\begin{table}
		\centering 
		\caption{Maximum relative error of \texttt{GPI} for different values of the burst control parameter $\beta$ over the interval $[0,5]$.}
		\begin{tabular}{ccccc}
			\hline
			$\beta$ & 2 & 2.2 & 2.4 & 2.6 \\
			\textbf{Maximum Relative Error} & 2.6$\mathrm{e}{-2}$ & 2.6$\mathrm{e}{-2}$ & 2.6$\mathrm{e}{-2}$ & 2.6$\mathrm{e}{-2}$\\
			\hline
		\end{tabular}
		\label{table:GPI_Diff_Relaxed}
	\end{table}
	
	Table \ref{table:GPI_Diff_Relaxed} evaluates the sensitivity of the \texttt{GPI} scheme to the burst control parameter $\beta$. The results are obtained by varying $\beta$ while keeping all other parameters fixed over the interval $[0,5]$. Increasing the burst length produces no such improvement in accuracy, indicating that the fast dynamics are already sufficiently relaxed onto the slow manifold for $\beta = 2$. This further validates the effectiveness of the proposed burst length selection.
	
	\begin{figure}
		\centering
		\includegraphics[width=0.6\linewidth]{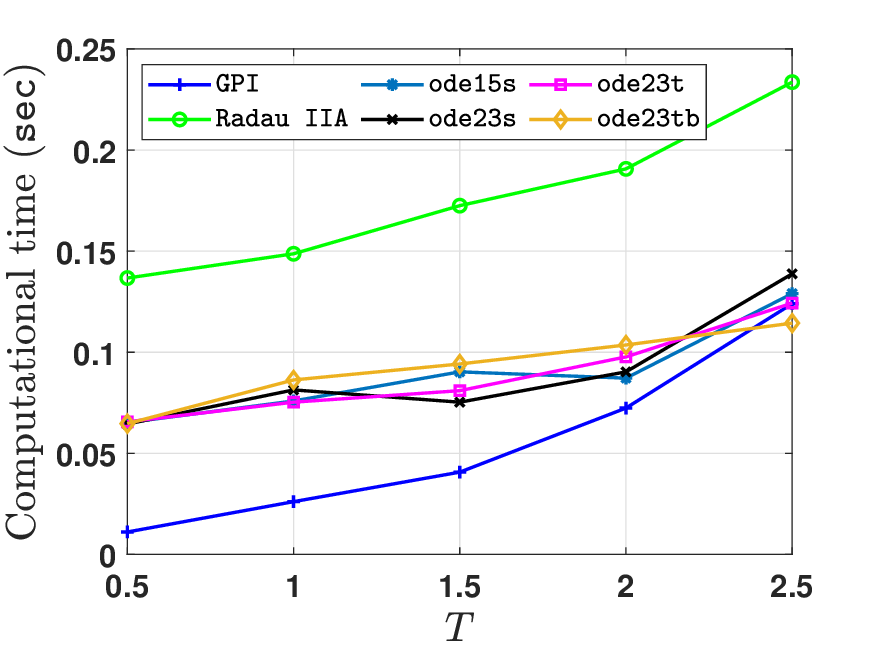} 
		\caption{The computational times required to solve problem \eqref{eqn:Diff_PDE} over the interval $[0,T]$ using the \texttt{GPI}, \texttt{Radau IIA}, \texttt{ode15s}, \texttt{ode23s}, \texttt{ode23t}, and \texttt{ode23tb} schemes are reported.}
		\label{fig:Comparison_Computational_Time_Diff_PDE}
	\end{figure}
	
	In terms of computational time, the \texttt{GPI} scheme requires the least time among the six methods considered, as shown in Figure~\ref{fig:Comparison_Computational_Time_Diff_PDE}. However, at $T = 2.5$, the \texttt{ode23tb} scheme requires slightly less time than \texttt{GPI}. In terms of memory requirements, all schemes exhibit similar usage.
	
	In conclusion, the proposed \texttt{GPI} scheme outperforms the \texttt{Radau IIA}, \texttt{ode15s}, \texttt{ode23s}, \texttt{ode23t}, and \texttt{ode23tb} schemes. These results demonstrate that the \texttt{GPI} scheme is particularly effective for problems with time-dependent spectra.

	\subsection{A highly oscillatory problem: Airy equation}\label{subsec:Airy}
	To demonstrate that the proposed method performs equally well for non-stiff problems, we consider the Airy equation written as a first-order system:
	\begin{equation}\label{eqn:Airy_Eqn}
		\begin{aligned}
			&u_1' = u_2, & u_1(0) = 1, \\
			&u_2' = -t\,u_1, & u_2(0) = 0.
		\end{aligned}
	\end{equation}
	This system is equivalent to the second-order Airy differential equation $u_1'' + tu_1 = 0$. This system is a non-autonomous dynamical system whose instantaneous eigenvalues are given by $\pm i\sqrt{t}$ for $t>0$. Consequently, the solution exhibits purely oscillatory behaviour for positive time, with a time-dependent frequency. In contrast, for $t<0$, the eigenvalues become real, leading to exponentially growing and decaying modes. This transition reflects a qualitative change in the system dynamics. Due to this time-dependent behaviour, the Airy equation serves as a useful benchmark for evaluating the performance of the \texttt{GPI} scheme in handling highly oscillatory solutions with time-varying characteristics.
	
	\begin{table}
		\centering
		\caption{Maximum absolute error in the \texttt{GPI} solution for different numbers of macro time steps when solving the Airy equation \eqref{eqn:Airy_Eqn}.} 
		\begin{tabular}{cc|cc}
			\hline
			\textbf{$\operatorname{N_t}$}& \textbf{$||u_1-u_1^{AS}||_\infty$} & \textbf{$\operatorname{N_t}$}& \textbf{$||u_1-u_1^{AS}||_\infty$} \\
			\hline
			1$\mathrm{e}{6}$& $8.23\mathrm{e}{-2}$& 4$\mathrm{e}{6}$& $1.87\mathrm{e}{-2}$\\ 
			2$\mathrm{e}{6}$& $3.86\mathrm{e}{-2}$& 8$\mathrm{e}{6}$& $9.20\mathrm{e}{-3}$\\ 
			\hline
		\end{tabular}
		\label{table:Airy_Accuracy}
	\end{table}
	
	Table \ref{table:Airy_Accuracy} reports the maximum absolute error in the \texttt{GPI} solution of the problem \eqref{eqn:Airy_Eqn}. The \texttt{GPI} scheme is applied without meso scale (i.e., $l_n=1$), using uniform macro and micro time steps. The total number of macro time steps over the interval $[0,100]$ is listed in the first and third columns of the Table \ref{table:Airy_Accuracy}. The uniform micro time step is chosen as $\delta T=1\mathrm{e}{-8}$. For oscillatory problems without stiffness, there is no intrinsic relaxation time, as the dynamics do not exhibit decay toward a slow manifold but instead remain persistently oscillatory. Therefore, the problem \eqref{eqn:Airy_Eqn} is evolved using only a single micro step (i.e., $M_{n,0}=0$) to compute the slope \eqref{eqn:GPI_Slope}. The second and fourth columns of Table \ref{table:Airy_Accuracy} present the corresponding maximum absolute errors with respect to the analytical solution. For $1\mathrm{e}{6}$, $2\mathrm{e}{6}$, $4\mathrm{e}{6}$ and $8\mathrm{e}{6}$ macro time steps in the entire interval $[0,100]$, the maximum absolute errors in the \texttt{GPI} solutions are $8.23\mathrm{e}{-2}$, $3.86\mathrm{e}{-2}$, $1.87\mathrm{e}{-2}$ and $9.20\mathrm{e}{-3}$, respectively. 
	
	\begin{figure}[h!]
		\centering
		\begin{subfigure}[t]{0.47\textwidth}
			\centering
			% include first image
			\includegraphics[width=1\linewidth]{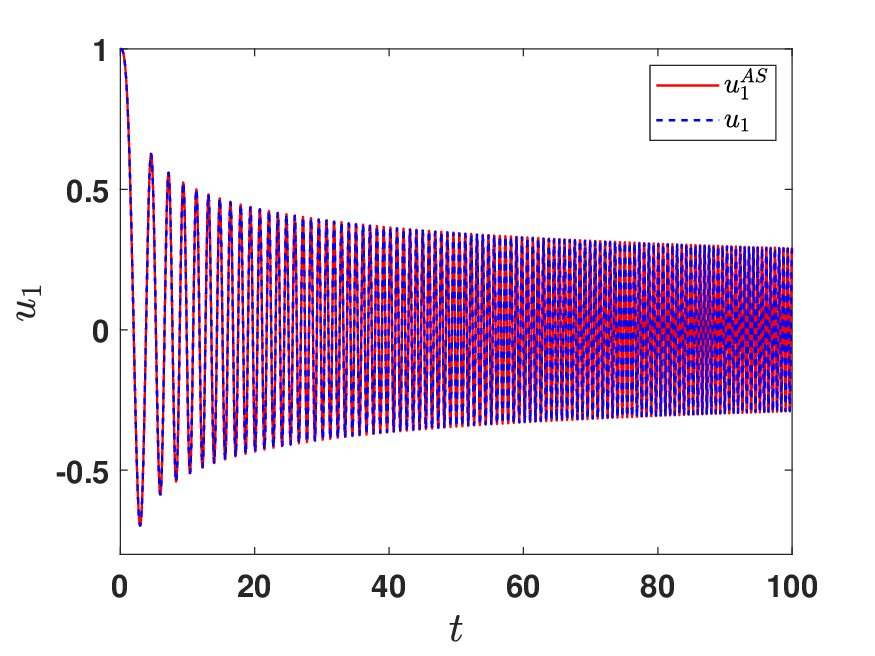}
			\caption{Analytical and \texttt{GPI} solutions of the Airy equation.}
			\label{fig:Airy_Soln_y}
		\end{subfigure}%
		\hspace{0.5cm}
		\begin{subfigure}[t]{0.47\textwidth}
			\centering
			% include second image
			\includegraphics[width=1\linewidth]{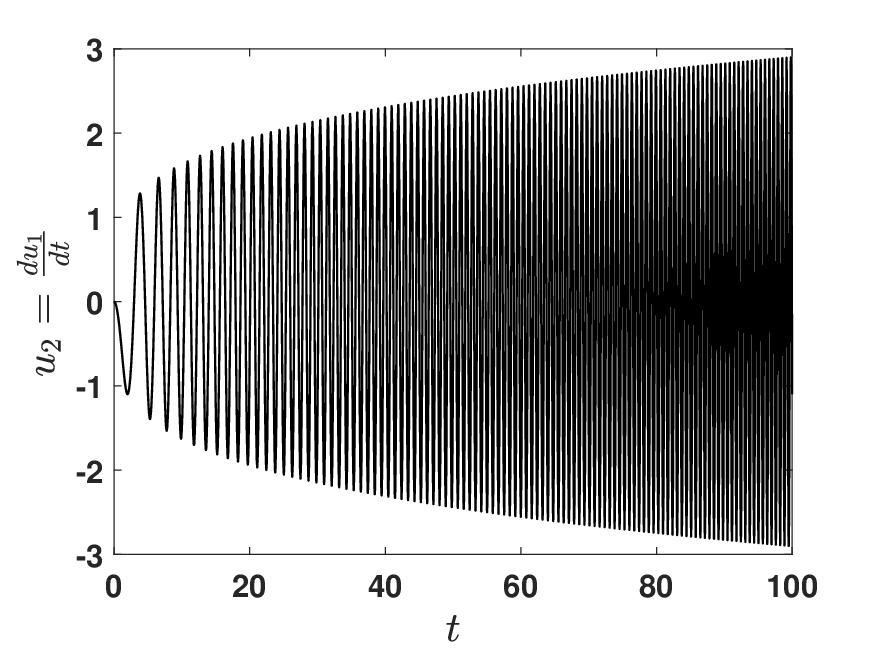}
			\caption{Slope of the Airy solution computed using the \texttt{GPI} scheme.}
			\label{fig:Airy_Soln_Dy}
		\end{subfigure}
		\caption{Solution profiles of the system \eqref{eqn:Airy_Eqn}.}
		\label{fig:Airy_Soln_y_Dy}
	\end{figure}
	
	Figure~\ref{fig:Airy_Soln_y} shows the \texttt{GPI} solution of $u_1$ along with the corresponding analytical solution for $\operatorname{N_t} = 8\mathrm{e}{6}$, demonstrating excellent agreement. Figure~\ref{fig:Airy_Soln_Dy} presents the evolution of $u_2$, which represents the derivative of $u_1$. Both solutions exhibit oscillatory behaviour with time-dependent frequency. It is observed that the amplitude of $u_1$ gradually decreases, while that of $u_2$ increases correspondingly. 
	
	\begin{table}
		\centering 
		\caption{Maximum absolute error of the \texttt{GPI} solution for different numbers of micro time steps $M_{n,0}$ per macro time step over the interval $[0,100]$.}
		\begin{tabular}{ccccc}
			\hline
			$M_{n,0}$ & 0 & 2 & 4 & 6 \\
			\textbf{Maximum Absolute Error} & $8.23\mathrm{e}{-2}$ & $8.23\mathrm{e}{-2}$ & $8.23\mathrm{e}{-2}$ & $8.23\mathrm{e}{-2}$\\
			\hline
		\end{tabular}
		\label{table:GPI_Airy_Relaxed}
	\end{table}
	
	Table~\ref{table:GPI_Airy_Relaxed} compares the accuracy of the \texttt{GPI} solution for different micro burst lengths. The results are obtained by increasing the number of micro steps $M_{n,0}$ within each macro time step while keeping all other parameters fixed over the interval $[0,100]$. It is observed that increasing the burst length does not improve the accuracy. This further confirms that, for oscillatory problems without stiffness, the use of additional micro steps is unnecessary due to the absence of a relaxation time scale.

	The results in Table~\ref{table:Airy_Accuracy} and Figure~\ref{fig:Airy_Soln_y_Dy} demonstrate that the \texttt{GPI} scheme is capable of accurately resolving highly oscillatory problems with time-dependent purely imaginary eigenvalues.

	\section{Conclusion}
	A key limitation of existing equation-free methodologies is that they are primarily designed for systems whose spectral properties remain constant or exhibit only mild temporal variation. Consequently, the associated scale separation is generally assumed to be nearly uniform in time. For systems with evolving spectra and dynamically changing scale separation, however, classical projective integration schemes such as \texttt{PI}, \texttt{PRK2}, \texttt{PRK4}, \texttt{PIRK2}, \texttt{PIRK4}, \texttt{PIG2}, and \texttt{PIG4} lack the flexibility required to adapt to local dynamical features, thereby limiting their effectiveness.
	
	To overcome these limitations, this article proposes a generalised projective integration (\texttt{GPI}) scheme, which introduces enhanced flexibility in the selection of macro-, meso- and micro-time steps, along with an adaptive micro-burst length. As demonstrated, the proposed \texttt{GPI} scheme provides a robust and efficient framework for solving non-autonomous systems with time-dependent spectra as well as scale separation properties.
	
	The main conclusions of this work are summarised as follows:
	
	\begin{enumerate}
		\item The \texttt{GPI} scheme unifies and generalises several existing multiscale approaches, that include the \texttt{PI} scheme, the \texttt{PI} versions of the \texttt{PD} and \texttt{GPD-I} schemes, first-order \texttt{PRK} method, \texttt{HMM}, \texttt{FLAVORS}, \texttt{VSHMM} and the \texttt{BA} strategy within \texttt{HMM}.
		\item Unlike classical projective integration methods that rely on two time scales (micro and macro), the \texttt{GPI} scheme incorporates three distinct time scales--micro, meso and macro--providing a greater flexibility in capturing multiscale dynamics.
		\item A comprehensive stability analysis of the \texttt{GPI} scheme is carried out, offering theoretical insight into its robustness under time-dependent spectral variations.
		\item The splitting of the stability region is analysed in detail. It is shown that for uniform $R^{n,m}$, a single connected stability region is split into at most two disconnected components, whereas for non-uniform $R^{n,m}$, it is split into multiple disconnected components. Moreover, analytical and numerical expressions for the splitting parameter are derived, and excellent agreement is observed between the two approaches.
		\item Problem-dependent guidelines for selecting variable burst lengths are developed, enabling improved adaptability of the method in practical applications.
		\item A new strategy is proposed to select the micro time step in the presence of time-dependent spectra. Within this framework, both mesoscopic and macroscopic time steps can vary adaptively, allowing for non-uniform microsimulation and extrapolation phases.
		\item Extensive numerical experiments demonstrate that the \texttt{GPI} scheme consistently outperforms existing projective integration methods (\texttt{PI}, \texttt{PRK2}, \texttt{PRK4}, \texttt{PIRK2}, \texttt{PIRK4}, \texttt{PIG2}, \texttt{PIG4}) across multiple performance metrics, including micro time step count, accuracy, computational time, memory usage and reduced reliance on microscale simulation. The method achieves a favourable balance between efficiency and accuracy, making it suitable for time-dependent spectra-based multiscale problems for long-time integration.
		
		The \texttt{GPI} scheme also demonstrates superior performance compared to widely used stiff solvers such as \texttt{ode15s} and \texttt{Radau IIA}, achieving lower step count, improved accuracy, minimal computational time and reduced memory under both adaptive and non-adaptive settings.
		\item For the time-dependent diffusion problem, where the entire spectrum evolves dynamically, \texttt{GPI} scheme outperforms the \texttt{Radau IIA}, \texttt{ode15s}, \texttt{ode23s}, \texttt{ode23t} and \texttt{ode23tb} schemes in terms of accuracy and computational time. The \texttt{GPI} scheme maintains accuracy and stability even when classical solvers fail to converge beyond a certain time, highlighting its robustness in challenging scenarios.
		\item The \texttt{GPI} scheme is further shown to be effective in resolving highly oscillatory problems with time-dependent imaginary eigenvalues, demonstrating its applicability beyond dissipative systems.
	\end{enumerate}
	Overall, the proposed \texttt{GPI} framework significantly extends the applicability of projective integration methods to a broader class of non-autonomous and multiscale problems with strongly time-varying spectra and scale separation. In this work, we have confined our attention to the first-order \texttt{GPI} scheme as a foundational step. In future work, we intend to develop a higher-order extension that preserves the key features introduced here, while also enabling a broader range of capabilities and deeper analytical insights.

\bibliographystyle{elsarticle-num}
\bibliography{GPI}
\end{document}